\documentclass[11pt]{article}
\usepackage
[pdftex
,a4paper
,hyperindex
=
true
,hyperfigures
=
true
,pagebackref
=
true
,bookmarks
=
true
,pdftitle
=
{TITLE
},pdfsubject
=
{},pdfauthor
=
{AUTHOR
},pdfkeywords={KEYWORDS},pdfpagemode=UseOutlines,menubordercolor=1 1 1,colorlinks=true,urlcolor=blue,citecolor=red,linkcolor=blue]{hyperref}

\usepackage{amsfonts,amssymb,epsfig}

\newcommand{\al}{\alpha}
\newcommand{\BE}{\begin{equation}}
\newcommand{\EE}{\end{equation}}
\newcommand{\la}{\lambda}
\newcommand{\csla}{c^s_\lambda} 
\newcommand{\RR}{{\mathbb R}}
\newcommand{\QQ}{{\mathbb Q}}
\newcommand{\ZZ}{{\mathbb Z}}
\newcommand{\NN}{{\mathbb N}}
\newcommand{\PP}{{\mathbb P}}
\newcommand{\EEE}{{\mathbb E}}

\newcommand{\ep}{\varepsilon}
\newcommand{\gr}{{\cal G}r}
\newcommand{\La}{\Lambda}
\newcommand{\EMH}{E_f (h)}
\newcommand{\cjk}{c^{i}_{j,k}}
\newcommand{\Hmin}{{h^{min}_f}}
 
\newcommand{\cla}{c_{\lambda}}
\newcommand{\pla}{\psi_{\lambda}}
\newcommand{\Laj}{\Lambda_j}
\newtheorem{lem}{Lemma}
\newtheorem{coro}{Corollary}
\newtheorem{Theo}{Theorem}
\newtheorem{prop}{Proposition}
\newtheorem{defi}{Definition}
\newcommand{\BP}{\begin{prop}}
\newcommand{\EP}{\end{prop}}
\newcommand{\BC}{\begin{coro}}
\newcommand{\EC}{\end{coro}}
\newcommand{\BL}{\begin{lem}}
\newcommand{\EL}{\end{lem}}
\newcommand{\BD}{\begin{defi}}
\newcommand{\ED}{\end{defi}}
\newcommand{\BT}{\begin{Theo}}
\newcommand{\ET}{\end{Theo}}

\def\E{{\hbox{I\kern-.2em\hbox{E}}}}

\author{ P. Abry\thanks{ Signal, Systems and Physics, Physics Dept., CNRS UMR 5672, Ecole Normale 
 Sup\'erieure de Lyon, Lyon, France 
{\tt patrice.abry@ens-lyon.fr}
},   S. Jaffard\thanks{ 
       Universit\'e Paris-Est, LAMA (UMR 8050),  UPEC, Cr\'eteil, France
{\tt jaffard@u-pec.fr}
  }, H. Wendt\thanks{CNRS, IRIT UMR 5505, University of Toulouse, France 
{\tt herwig.wendt@irit.fr}} 
}

 \date{\today } 
\title{  Irregularities and Scaling in Signal and Image Processing: Multifractal Analysis 
 }

\begin{document}
 
 \maketitle

{ \bf Abstract:} 

B. Mandelbrot gave a new birth to the notions of scale invariance, selfsimilarity and non-integer dimensions, gathering them as the founding corner-stones used to build up  \emph{fractal geometry}. 
The first purpose of the present contribution is to review and relate together these key notions,  explore their interplay and show that they are different facets of a same intuition.
Second, it will explain how these notions lead to the derivation of the mathematical tools underlying multifractal analysis.
Third,  it will reformulate these theoretical tools into a wavelet  framework, hence enabling their better theoretical understanding as well as their efficient practical implementation. 
B. Mandelbrot used his concept of  \emph{fractal geometry} to analyze  real-world applications of very different natures. 
As a tribute to his  work,  applications  of various origins, and where multifractal analysis proved fruitful, are revisited to illustrate the theoretical developments proposed here. \\



\clearpage
\newpage
\tableofcontents

\clearpage
\newpage

\section{Introduction }

\label{intro}

As many students, before and after us, we first met Beno\^{\i}t Mandelbrot through his books 
\cite{Mandelbrot1977,Man1}; 
 indeed, they were scientific bestsellers, and could be found in evidence in most scientific bookshops. A quick look was sufficient to realize that Beno\^{\i}t was an unorthodox scientist; these books radically differed from what we were used to, and were trespassing several taboos of the time: First,  they could not be straightforwardly associated with any of the standardly labelled fields of science, which had been taught to us as separate subjects.
This could already be inferred from their titles: { \em The fractal geometry of nature} \cite{Man1} constitutes a statement that mathematics and natural sciences are intrinsically mixed, and that the purpose will not be to disentangle them artificially, but rather to explore their interactions.  
The contents of these books  confirmed this first feeling: Though they contained a large proportion of very serious mathematics, they did not follow the old \emph{definition-theorem-proof} articulation  we were used to. 
The focus was directly on the observation of figures and on their pertinence for modeling;  the mathematics were there to help  and comfort the deep geometric intuition that Beno\^{\i}t had developed.  
This does not mean that he considered mathematics as an ancillary discipline. 
To the contrary, many testimonies confirm  how enthusiastic he was when one of his mathematical intuitions was proven correct, and how he immediately incorporated the new mathematical methods in his own toolbox, which then helped him to sharpen his intuition and go further.  

Another revolutionary point  was that  these books reversed the classical focusses: 
While one did find familiar mathematical objects (the nowhere differentiable Weierstrass functions, or the  \emph{devil's staircase} had been encountered before), the role he made them play was radically different:  
So far, they  had only been counterexamples, or pathological objects pedagogically used to shed light on possible pitfalls;  
instead of this restrictive role, Beno\^{\i}t used them as the key examples and cornerstones to found a new geometry based on roughness and selfsimilarity.  
These books were  opening a subject, and drawing tracks through a  new territory; their purpose was not to give the usual  polished, final description of  a well understood area, but to open windows  on how science advances,  and they had a strong influence on our vocations to research.  

Later, personal meetings played a decisive role; let us only mention one such occasion:  At the end of his PhD, one of us (S. J.) met Beno\^{\i}t who was  visiting  Ecole Polytechnique.  
Beno\^{\i}t had heard of wavelets, which, at the time, were a  new tool still in  infancy, and immediately envisaged possible interactions between wavelets and fractal analysis.
A  key feature of fractals is scaling invariance, and, because wavelets can be deduced one from the others by translations and dilations of a given function, their particular algorithmic structure should reveal  the scaling invariance properties of the analyzed objects. 
Driven by his insatiable curiosity, Beno\^{\i}t wanted to learn everything about wavelets (which was actually not so much at the time!) and then, with his characteristic generosity,  he shared  with this student he barely knew his intuitions on the subject.  
Thus, the present contribution can be read as a tribute to some of the scientific lines Beno\^{\i}t had dreamed about at the time: It will show how and  why wavelet techniques  yield powerful tools for analyzing {fractal} properties of  signals and images, and thus supply new collections of tools for their classification  and modeling.  
 
\subsection{Scale invariance and self similarity}

\noindent {\bf Exact (or deterministic) self similarity.} \quad  Fractal analysis can roughly be thought of as a way to characterize and investigate the properties of  irregular sets. 
Note that \emph{irregularity}  is a negatively defined word; the description of the scientific domain it covers is therefore not clearly feasible,  and had not been tried  before B. Mandelbrot  developed the  concepts of fractal geometry. 
One of his leading ideas was a notion of \emph{regularity in irregularity or selfsimilarity}: The object is irregular, but there is some invariance in its behavior across scales. 
For some sets, the notion of \emph{selfsimilarity} replaces the notion of smoothness. 
For instance, the triadic Cantor set  is selfsimilar, in the sense that it is made of two pieces  which are the same as  the whole set is shrunk by a factor 3.  The purpose of multifractal analysis is to transpose this idea from the context of geometrical sets to that of functions.   
Sometimes, this transposition can be directly achieved, when the geometrical set consisting of the graph of the function also displays some selfsimilarity. 
Consider the example of  the \emph{Weierstrass-Mandelbrot functions} (a slight variant of the famous  Weierstrass functions, see \cite{Man1} pp. 388-390): 
\BE \label{weier}  { \cal W}_{a, H} (x) = \sum_{n = -\infty}^{+ \infty} \frac{\sin ( a^n x)}{a^{H n}}  \hspace{15mm} \mbox{for} \hspace{6mm}  a >1   \hspace{6mm} \mbox{and } \hspace{6mm}  H \in (0,1) \EE 
(note that the series converges  when $n \rightarrow + \infty$ because $a >1$ and when $n \rightarrow - \infty$ because 
$ | \sin ( a^n x)  | \leq |a^n x | $ and $H <1$).    Those functions clearly satisfy the exact selfsimilarity relationship
\BE  \label{selfsim} \forall x \in \RR, \hspace{6mm} { \cal W}_{a, H} (  a x) =  a^H { \cal W}_{a, H} (x) . \EE
Note that  selfsimilar functions associated with a  selfsimilarity exponent $H >1$ can be obtained using slight variants:  For instance, the function
\[  \sum_{n = -\infty}^{+ \infty} \frac{ \cos ( a^n x) -1}{a^{H n}}   \hspace{15mm} \mbox{for} \hspace{6mm}  a >1   \hspace{6mm} \mbox{and } \hspace{6mm}  H \in (0,2) \]
is selfsimilar  with a  selfsimilarity exponent $H $, which can go up to  2;   and  
the function
\[  \sum_{n = -\infty}^{+ \infty} \frac{ \sin ( a^n x) - a^n x  }{a^{H n}}   \hspace{15mm} \mbox{for} \hspace{6mm}  a >1   \hspace{6mm} \mbox{and } \hspace{6mm}  H \in (1,3) \]
is selfsimilar  with a  selfsimilarity exponent $H $ which can go up to 3. Note  that this ``renormalization technique'' can be pushed further to deal with higher and higher values of $H$, yielding smoother and smoother functions: One thus obtains functions that are everywhere  $C^{[H ]} $ and nowhere $C^{[H ] +1} $.    \\  

\noindent {\bf Random (or statistical) self similarity.} \quad  Functions satisfying (\ref{selfsim}) are very rare; 
therefore the notion of exact selfsimilarity can be considered too restrictive for real world data modeling.  
Fruitful generalizations of this concept can be developed into several directions.  
A first possibility lies in  weakening the exact  deterministic relationship (\ref{selfsim}) into a probabilistic one: 
The (random) functions $a^H f( ax)$ do not coincide sample path by sample path, but  share the same statistical \emph{laws}. 
A stochastic process $X_t$  is said to be { \em selfsimilar}, with selfsimilarity exponent $H$ iff 
\BE  \label{autosimstoch}  \forall a >0, \;\;\; \{ X_{a t},\;  t \in \RR \}Ê\;  \stackrel{{\cal L} }{=} \{ a^H X_t ,\;  t \in \RR \} \EE
(see Definition \ref{defeuqalaw} where the precise definition of equality in law of processes is recalled). 

The first example of such selfsimilar processes is that of fractional Brownian motion (hereafter referred to as fBm), introduced by Kolmogorov \cite{Kolmo40}. Its importance for the modeling of scale invariance and fractal properties in data was made explicit and clear by Mandelbrot  and Van Ness in \cite{mvn68}. 
This article is characteristic of Mandelbrot's genius, which could perceive similarities between   very distant disciplines: Motivations simultaneously rose  from hydrology (cumulated water flow), economic time series and fluctuations in solids. 
Notably, B. Mandelbrot explained that fBm (with $0.5 < H <1$) was well-suited to model long term dependencies, or long memory effects, in data, a property he liked to refer to as the \emph{Joseph effect} in a colorful and self-speaking metaphoric comparison to the character of the Bible (cf. e.g., \cite{H10}). 
Selfsimilar processes have since been the subject of numerous exhaustive studies, see for instance the books  \cite{Abry2002a,CohIst,EmbrechtsMaejima2001,Samorodnitsky1994}.   Note also that B. Mandelbrot put selfsimilarity (which he actually more accurately called { \em selfaffinity}) as  a central notion  in the study stochastic processes in views of applications, see \cite{Mandelbrot2002} for a panorama of his works in this direction. 
 
\subsection{Selfsimilarity and Wavelets} 

\label{secwavmot}

Let us now briefly come back to B. Mandelbrot's intuition and see  how  wavelet  analysis  can reveal relationships such as  (\ref{selfsim})  or  (\ref{autosimstoch}),  as well as other important probabilistic properties of stochastic processes. 
In its simplest form, an orthonormal wavelet basis of $L^2 (\RR)$ is defined as the collection of functions
\BE \label{wav} 2^{j/2}\psi (2^jx-k), \;\; \mbox{ where} \;\; 
j, k\in\ZZ , 
\EE
and the wavelet $\psi$ is a smooth, well localized function. Wavelet coefficients are further defined by 
\[ {    c_{j,k} = 2^{j} \int f(x) \;   \psi  (2^j x-k) \; dx .  }  \]
A $L^1$-normalization is used (rather than the classical $L^2$-normalization), as it is better suited to express scale invariance relationships, as shown in (\ref{samelaw})  below. 
 
Let $f$ denote a function satisfying (\ref{selfsim})  with $a =2$, i.e.  such that
 \[ \forall x \in \RR, \hspace{6mm}  f(2x) = 2^H  f(x); \]
 then, clearly, its wavelet coefficients satisfy 
 \[ \forall j,k, \hspace{8mm}  c_{j,k}= 2^H c_{j+1,k}.\]
 Similarly, in the probabilistic setting,  if (\ref{autosimstoch}) holds, then the  (infinite) sequences  of random vectors  $ {\cal C}_j \equiv (2^{H j} c_{j,k})_{k \in \ZZ} $   satisfy the following  property:
 \BE     
 \label{samelaw} \{ {\cal C}_j\}_{j \in \ZZ}   \hspace{4mm} \mbox{ share the same law.}  
 \EE  
Other probabilistic properties have simple wavelet counterparts; let us mention a few of them.

 Let  $X_t$ be  a {  \em Markov process}, i.e. a stochastic process  satisfying the property
\[ \mbox{  $ \forall s >0$, $ X_{ t+s} -X_t$ is independent of the $(X_u)_{ u \leq t}$ } ; \]   
typical examples of Markov processes are supplied by  { \em Brownian motion}, or, more generally, by { \em L\'evy processes}, i.e. processes  with stationary and independent increments. 
Because wavelets have vanishing integrals, the Markov property implies that, if  the supports of the $\psi (2^{j} x-k )$ are disjoint, then the corresponding coefficients $ c_{j, k}  $ are independent random variables.  
  
Recall that a stochastic process $X_t$ has { \em stationary increments} iff 
    \[ \mbox{  $ \forall s >0$, $ X_{ t+s} -X_t$ has the same law as  $X_s$ };  \] 
      typical examples are   { \em fractional Brownian motions}, or { \em L\'evy processes}. 
If $X_t$ has {stationary increments}, then
\BE \label{statinc} \forall j, \;\; \mbox{ the random variables } \; (c_{j,k})_{k \in \ZZ} \;  \mbox{ share the same law. } \EE 
Note that this property holds for each given $j$, but implies no relationship between wavelet coefficients at different scales; however, if $X_t$ is a selfsimilar process of exponent $H$ with stationary increments, then it follows from (\ref{samelaw}) and (\ref{statinc}) that 
 the random variables $ 2^{Hj}c_{j,k} $  all share the same law (for all $j$ and $k$). 

A random variable $X$ has a  {\em  stable distribution} iff 
\[ \forall a, b \geq 0, \;\; \exists c >0\; \mbox{ and} \;  d\in \RR \;\; \; \mbox{ such that } \; 
 \;\; \;  \;   aX_1 + b X_2  \stackrel{{\cal L} }{=}  cX + d \] 
 where $X_1$ and $X_2$ are independent copies of $X$.  
 A random process is stable if all linear combinations of the $X (t_i)$ are stable.  
This property implies that the wavelet coefficients of a stable process are stable random variables.
Further, it can actually be shown that any finite collection of wavelet coefficients forms a stable vector (cf. e.g., \cite{abpt00,PIPIRAS:2003:A,PIPIRAS:2007:A,STOEV:2002}).
 
Finally, if $X_t$ is a { \em Gaussian process}, then any finite collection of its wavelet coefficients  is a Gaussian vector.  
Note however a common pitfall: 
Unless some additional stationarity property holds, the sole Gaussianity hypothesis for the marginal distribution of the process does { \em not} imply that the empirical distributions of the coefficients $(c_{j,k})_{k \in \ZZ}$ are close to Gaussian: It can only be said that they will consist of  \emph{Gaussian mixtures} (which can strongly depart from a Gaussian distribution). 
The statistical properties of wavelet coefficients of fBm were studied in depth in the seminal works of P. Flandrin  \cite{f89a,Flandrin1992} that significantly contributed to the popularity of wavelet transforms in the analysis of scale invariance. 

Among the results that we mentioned above, it is important to draw a difference between those that are a consequence of the sole fact that wavelet coefficients are defined as a linear mapping, and therefore would hold for any other basis (Gaussianity and stability) and those that are the consequence of the particular algorithmic nature of wavelets (selfsimilarity and stationarity).  These (rather straightforward) results are either implicit or explicit in the literature concerning wavelet analysis of stochastic processes; however, in order to give a flavor of the ideas and methods involved, we give a  proof of (\ref{samelaw}) at the end of  Section \ref{secwavedis}. 

Further, let us mention that  all these wavelet properties do not constitute exact characterizations of the probabilistic properties of the corresponding processes, but only implications, because of a deficiency of wavelet expansions  that will be explained in Section \ref{secwavedis}.  \\
  
 \subsection{Beyond selfsimilarity: Multifractal analysis}
 \label{sec-ss2mf}
 
\noindent {\bf Beyond selfsimilarity.} \quad  An obvious drawback of using {\em exact} probabilistic properties, such as selfsimilarity,  or independence of increments is that, though they can relevantly model physical properties of the underlying data, they usually do not hold exactly for real-life signals.
Two limitations of selfsimilarity can classically be pointed out. 
First, the scaling properties may not hold for all scales (as required by the mathematical definition of selfsimilarity). 
This may  be due to noise corruption or it may also result from the fact that the \emph{physical} phenomena implying those properties intrinsically act  over a finite range of scales only.
Second, the appealing property of selfsimilarity (the sole selfsimilarity parameter $H$ contains the  essence of the model) may also constitute its limitation: Can one really  expect that the complexity  of real-world data can be accurately encapsulated in  one single parameter only? 
Let us examine such issues.\\

\noindent {\bf Multiplicative Cascades.} \quad As mentioned above, B.~Mandelbrot coined fBm as one of the major candidates to model scale invariance. However, combining selfsimilarity and stationary increments implies that all finite moments of the increments of fBm (and of any process whose definition gathers these two properties)  verify, for all $p$s such that $ \E( |X(1)|^{p}) < \infty$ \cite{Samorodnitsky1994,EmbrechtsMaejima2001},
\begin{equation}
\forall \tau \in \RR, \hspace{6mm}  \E (|X(t+\tau) -X(t)|^{p} ) = \E (|X(1)|^{p} ) \; |  \tau|^{pH}. 
\end{equation}
Though very powerful, this result can also be regarded as a severe limitation for applications, where empirical estimates of moments,  $\int |X(t+\tau) -X(t)|^p dt $, computed on the data assuming ergodicity, may instead behave as: $\int |X(t+\tau) -X(t)|^p dt = \tau^{\zeta(p)} $, for a range of $\tau $, and where $\zeta(p)$ is, by construction, a concave function, but is not necessarily linear. 
This is notably the case when analyzing velocity or dissipation fields in the study of hydrodynamic turbulence  where strictly concave scaling functions were obtained (cf. e.g., \cite{Frisch1995} for a review of experimental results).
This was justified by a celebrated argument due to Kolmogorov and Obukhov in 1962 \cite{kolmogorov62,obukhov62,Yaglom1966}: 
Navier-Stokes equation,  governing fluid motions, implies that the (third power of the) gradient of the velocity $u(x+r)-u(x) $ is proportional to the mean of the local energy $\epsilon_r$ dissipated into a bulk of size $r$: $\E (u(x+r)-u(x))^3 \propto \E (\epsilon_r) \cdot  r$. 
If $\epsilon_r$ is constant, then any power scales as  $\E (u(x+r)-u(x))^p \propto \epsilon_r^p  r^{p/3}$.
But, in turbulence, $\epsilon_r$ is not constant over space and should rather be regarded as a random variable, hence implying:    $\E (u(x+r)-u(x))^p \propto \E (\epsilon_r^p ) \;  r^{p/3} $ which naturally implies a general scaling of the form $  \sim \E (u(x+r)-u(x))^p \propto r^{\zeta(p)}$, with $\zeta(p)$  which is  likely to depart from  the linear behavior $pH$ which is implied by  the  selfsimilarity  hypothesis.  \\ 

To account for this form of observed scale invariance,  new stochastic models  have been proposed: Multiplicative cascades. 
Numerous declinations of cascades were proposed in the literature of hydrodynamic turbulence to model data, see e.g., \cite{Frisch1995} (on the mathematical side, see  the book  \cite{Panor} for a detailed  review,  and  also the contribution of J. Barral and J. Peyri\`ere \cite{BarPeyr} in the present volume, for a lighter introduction).
For the sake of simplicity, the simplest (1D) formulation is presented here;  however, it  contains the key ingredients present in all other more general cascade models. 
Starting from an interval (arbitrarily chosen as the unit interval), a cascade is constructed as the repetition of the following procedure ($c$ being an arbitrary integer): 
At iteration $j$, each of the $c^{j-1}$ sub-intervals $\{ I_{j-1,k}, k=0, \ldots c^{j-1}-1 \} $,  obtained at  iteration $j-1$, is split into $c$ intervals of the same length $c^{-j}$;  a mean-$1$ positive random variable $w_{j,k} $ (drawn  independently and with identical distribution) is associated with each $\{ I_{j,k}, k=0, \ldots c^{j}-1 \} $. 
After an infinite number of iterations, the \emph{cascade} $W$ is defined   as the limit (in the sense of measures) of the finite products  
\[  W_J (x) dx = \displaystyle\prod_{j \leq J, \; x \in I_{j,k} }  w_{j,k} \; dx . \]  
Mathematically, this sequence of measures converges  under mild  conditions  on the distribution of the multipliers $w$ (cf. \cite{KaPe}).
Such constructions and variations have been massively used to model the Richardson energy cascade 
implied by the Navier-Stokes equation  \cite{richardson22} (cf.  \cite{Frisch1995} for a review of such models).
One of the major contributions of B.~Mandelbrot  to the field of hydrodynamic turbulence was proposed in \cite{Mandelbrot1974}, where he presented  the themes and variations on multiplicative cascades in a single  framework and studied their common properties, notably in terms of scale invariance. 
This celebrated contribution paved the way towards the creation of the word \emph{multifractal}, to the notions of multifractal processes and multifractal analysis, where \emph{multi} here refers to the fact that instead of a single scaling exponent $H$, a whole family of exponents $\zeta(p)$ are needed for a better description of the data. 
B.~Mandelbrot also discussed the importance of the distribution of the multipliers and how multiplicative cascades may significantly depart from log-normal statistics that could be obtained from a simplistic argument: 
\[  \Pi w_{j,k} = \Pi e^{\log (w_{j,k})} = e^{\sum \log w_{j,k}},  \]  where $\sum \log w_{j,k}$  tends to  a normal distribution, hence $ e^{\sum \log w_{j,k}}$ should tend to  a log-normal distribution (cf. \cite{Mandelbrot1990,N14}, and also \cite{Yaglom1966}).
This opened many research tracks proposing alternatives to the log-normal  model, the most celebrated one being the log-Poisson model \cite{sheleveque01}. 
This will be further discussed in Section \ref{sec-turb}.
Much later, multiplicative cascades were connected to Compound Poisson cascades and further to Infinitely Divisible Distributions (cf. e.g., \cite{novikov94}). 
Mandelbrot himself contributed to these latter developments with the earliest vision of Compound Poisson Cascades \cite{BaM02} and  also with his celebrated fBm in multifractal time \cite{Mandelbrot1999},  which  will be considered in Section \ref{sec-finance}.  
This gave birth to multiple constructions of processes with \emph{multifractal} properties (cf. e.g., \cite{Bacry2001,Bacry2002,Bacry2003,BaSe3,BaSe2,Chainais2005,Chainais2005a}).

Compared with the concept of selfsimilarity (and in particular fBm), which is naturally related to random walks and hence to additive constructions, cascades are, by  construction, tied to a multiplicative structure. 
This  deep difference  in nature explains why practitioners, for real-world applications, are often eager to distinguish between these  two classes of models. 
This distinction is however often misleadingly confused with the opposition of mono- versus multi-fractal processes. 
This will be further discussed in Section~\ref{sec-mfMF}. \\

\noindent {\bf Scaling functions and scaling exponents.} \quad 
The empirical observations of strictly concave scaling exponents, made on hydrodynamic turbulence data, led to the design of a new tool for analysis:  \emph{Scaling functions}.
 Let $f: \RR^d \rightarrow \RR$. The { \em Kolmogorov scaling function}  of $f$  (see \cite{Kol41}) is the function $\eta_f (p)$ implicitly defined by the relationship
\begin{equation} \label{kolmo} \forall p \geq 1, \hspace{8mm} 
\int |f(x+h) - f(x)|^p dx \quad\sim\quad   |h|^{\eta_f(p)},
\end{equation}
the symbol $\sim$ meaning that 
\BE  \label{scalKolmo1} \eta_f (p) =\lim_{|h |  \rightarrow 0} \frac{\log \left(\displaystyle\int |f(x+h) - f(x)|^p dx\right) }{\log |h| }  . \EE
Note that this limit  does not necessarily    always exist; when it does, it reflects the fact that   $f$ shows clear scaling properties (and  it  is needed to hold in order to numerically  compute scaling functions). However, we  are not aware of any simple and general mathematical assumption implying that  it does; 
therefore, it is important to formulate a mathematical and slightly more general  definition of the scaling function,  which  is well defined for any function $f\in L^p $: 
\BE  \label{scalKolmo} \forall p \geq 1, \hspace{6mm} \eta_f (p) =\liminf_{|h |  \rightarrow 0} \frac{\log \left(\displaystyle\int |f(x+h) - f(x)|^p dx\right) }{\log |h| }.  \EE 

It is important to note that, if data are smooth (i.e., if one obtains that $\eta_f (p) \geq p$), then one has to use differences of order 2 or more in (\ref{scalKolmo1}) and (\ref{scalKolmo})  in order to define correctly the scaling function.

For application to real-world data, which are only available with a finite time or space   resolution, (\ref{scalKolmo1})  can be interpreted as a linear regression  of the  log of the left hand side of (\ref{kolmo})  vs. $\log (|h|)$. 
The same remark holds  for the other scaling functions that we will meet.
 \\

\noindent {\bf Multifractal analysis.} \quad In essence, multifractal analysis consists in the determination of scaling functions (variants to the original proposition of Kolmogorov $\eta_f$ will be considered later). 
Such scaling functions can then be involved into classification or model selection procedures. 
Scaling functions constitute a much more flexible analysis tool than the strict relation (\ref{selfsim}).
For many different real-world data stemming from applications of various natures, relationships such as (\ref{kolmo}) are found to hold, while this is usually not  the case for (\ref{selfsim}), see Section~\ref{sec-App} for illustrations on real world data.

An obvious key advantage of the use of the scaling function $\eta_f(p)$ is that its dependence in $p$ can take a large variety of forms, hence providing versatility in adjustment to data. 
Therefore, multifractal analysis, being based on a whole function (the scaling function $\eta_f(p)$, or some of its variants introduced below), rather than on the single selfsimilarity exponent, yields much richer tools for classification or model selection. 
The scaling function  however satisfies a few constraints, for example, it has to be a concave non-decreasing function (cf. e.g., \cite{Jaffard2006,WENDT:2007:E}).

Another fundamental difference between the selfsimilarity exponent and the scaling function  lies in the fact that, for stochastic processes, the selfsimilarity exponent is, by definition,  a deterministic quantity, whereas the scaling function is constructed for  each sample path, and therefore is  a priori defined  as a random quantity. 
And, indeed, even if  the limit in (\ref{scalKolmo})  turns out to be deterministic for many examples  of  stochastic processes (such as fBm, L\'evy processes, or multiplicative cascades, for instance), it is not always  the case, as shown, for instance, by the example of fairly general classes of Markov processes, such as the ones  studied in \cite{BaJaFoSe}. 
Note that, if a stationarity assumption can be  reasonably assumed, then  one can define a deterministic quantity by taking  an expectation instead of a space average in (\ref{scalKolmo}), which  leads to a definition of the scaling function through moments of increments. 
Under appropriate ergodicity assumptions,  time (or space) averages can be shown to coincide with ensemble averages, and, in this case, both approaches lead to the same scaling functions, see \cite{Riedi2003}  for a discussion  (checking such ergodicity assumptions for experimental data does, however, not seem feasible and can become an intricate issue, cf. \cite{lac04}). \\


\subsection{Data analysis and Signal Processing} 

In its early ages, the concept of { \em scale invariance}  has been deeply tied to that of $1/f$-processes, that is, to second order stationary random processes whose { \em power spectral density} behave as a power law with respect to frequency over a broad range of frequencies: $\Gamma_X(f) \sim C |f|^{-\beta} $, where the frequency $f$ satisfies $ f_m \leq f \leq f_M$, with $ f_M/f_m \gg 1 $. 
From that starting point, spectrum analysis was regarded as the natural tool to analyze scale invariance and estimate the corresponding scaling parameter $\beta$. 
In that context, the contribution of B. Mandelbrot and J.  Van Ness \cite{mvn68} can be regarded as seminal since, elaborating on \cite{Hurst1951}, it put forward fBm, and its increment process, fractional Gaussian noise (fGn),  as a model that extends and encompasses $1/f$-processes: fGn is a $1/f$-process, with $\beta = 2H-1$, fBm is a non-stationary process that shares the same scale invariance intuition. 
This change of paradigm immediately raised the need for new and generic estimation tools for measuring the  parameter $H$ from real world data, which were no longer based on classical spectrum estimation. 
Relying strongly on the intuitions beyond selfsimilarity and long memory (a concept deeply tied with selfsimilarity, cf. e.g. \cite{Beran1994}),  the R/S estimator (for Range-Scale Statistics) has been the first such tool proposed in the literature, to the study of which B. Mandelbrot also contributed \cite{Mandelbrot1998,H10,H25}. The R/S estimator has notably been used in finance and economy, in particular by D. Strauss-Kahn and his collaborators,  (see e.g., \cite{HMSK78, Peters1994}). 
Later, in the 90s, with the advent of wavelet transforms, seminal contributions studied the properties of the wavelet coefficients of fBm \cite{f89a,Flandrin1992,tk92} and opened avenues towards accurate and robust wavelet-based estimations of the  selfsimilarity parameter $H$ \cite{agf95,Abry1998,Veitch1999,Veitch2001}.

To some extent, the present contribution can be read as a tribute to ideas and tracks originally addressed by B. Mandelbrot (and others) aiming at estimating the parameters characterizing scale invariance in real world data: The wavelet based formulations of multifractal analysis (i.e., of the estimation from data of scaling functions and multifractal spectra) devised in Sections~\ref{secwav} and \ref{secmulti} can be read as continuations and extensions of these pioneering works towards richer models and  refined variations  on the  characterization of  scale invariance. \\

\subsection{Goals and Outline} 

In Section~\ref{sec-F2MF}, it is explained how various paradigms pertinent in fractal geometry,  such as selfsimilarity or fractal dimension, can be shifted to the setting of functions,  thus leading to the introductions of several { \em scaling functions}. 
In Section~\ref{sec-GeomAnalysis}, it is shown how such scaling functions can be reinterpreted in terms of function space regularity: On one hand, it allows to derive some of their properties; on other hand, it paves the way towards equivalent  formulations that will be derived later. 
Section \ref{secwav} introduces  wavelet techniques, which allow to reformulate these former notions in a more efficient setting, from both a theoretical and  pratical point of view.  
Section \ref{secmulti} introduces the { \em multifractal formalism}, which offers a new interpretation for scaling functions in terms of fractal dimensions of the sets of points where the functions considered has a given pointwise smoothness.  
Applications of  the tools and techniques developed  here will be given in Section~\ref{sec-App}, aiming at illustrating the various aspects of multifractal analysis addressed in the course of this contribution. 
Such applications are chosen either because B. Mandelbrot significantly contributed to their analysis and/or because they illustrate particularly well some concepts which will be developed here. 
The selected applications stem from widely different backgrounds, ranging from Turbulence to Internet, Heart Beat Variability and Finance (in 1D), and from natural textures to paintings by famous masters (in 2D), thus showing the extremely rich possibilities of these new techniques. 
  
\section{From fractals to multifractals}
\label{sec-F2MF}

\subsection{Scaling and Fractals}

\noindent {\bf Dynamics versus Geometry.} \quad One of the major contributions of B. Mandelbrot  has been to understand the benefits of using mathematical tools and objects already defined in earlier works but regarded as incidental, such as fractal dimensions, in order to explore the  world of irregular sets, to classify them and to study their properties.   
A key  example is supplied by the { \em Von Koch curve}  (cf. e.g., \cite{ps88,Schroeder1991} and  Fig.~\ref{fig-BV}). 
The important  role that B. Mandelbrot  gave to this curve   is typical of  his way of diverting  a mathematical object from its initial purpose and of making it fit to his own views; indeed, it was initially introduced  as an example of the graph of a continuous, nowhere differentiable function (see  the book by G. Edgar
 \cite{Edgar}, which made available many classic  articles of historical importance in fractal analysis).  
B. Mandelbrot shifted the status of this  curve from a peripherical example of strange set to a central element in the construction of his theory, and he baptized it the  { \em Von Koch snowflake}, a typical example of the poetic accuracy that  B. Mandelbrot  displayed for coining new names (other examples, being the { \em devil's staircase},  { \em L\'evy flights},  {\em  Noah  or Joseph effects},...). 
This example will be revisited below in Section \ref{secwavcoe} to validate a numerical method for multifractal analysis.

Two features --- dynamics and geometry --- are intimately tied into the Von Koch snowflake: It is constructed from (and hence invariant under) a dynamical system consisting of the iterative repetition of a \emph{split and reproduce} procedure (highly reminiscent of the \emph{split and multiply} procedure entering the definition of multiplicative cascades, cf. Section~\ref{sec-ss2mf} earlier); the resulting geometrical set has well defined geometrical properties  that can 
be put into light through the determination of its  \emph{ box dimension}, whose definition we now recall.\\

\noindent {\bf Dimensions.} \quad 
\BD  \label{defbox} Let   $A$ be a bounded subset of $ \RR^d$; if 
 $\ep>0$,  let $N_\ep (A)$  be the smallest number  such that there exists a covering of  $A$ by $N_\ep (A)$ balls of radius $\ep$. 

The upper and lower box dimension of $A$  are respectively given by 
\BE \label{bd}  \overline{dim}_B (A)  =\limsup_{\ep \rightarrow 0} \frac{\log N_\ep (A)}{-\log \ep}  , \hspace{5mm} \mbox{and} \hspace{8mm}  \underline{dim}_B (A)  =\liminf_{\ep \rightarrow 0} \frac{\log N_\ep (A)}{-\log \ep}  .\EE
\ED


When both limits coincide (as it is the case for the Von Koch snowflake), they are referred to as the { \em  box dimension }  of the set $A$: 
\BE 
\label{dimboi} {
dim}_B (A)  =\lim_{\ep \rightarrow 0} \frac{\log N_\ep (A)}{-\log \ep}. 
\EE 
This implies that $N_\ep (A)$ displays an approximate  power-law behavior with respect to scales: 
\begin{equation}
\label{equ-scalingBoxdim}
  N_\ep (A) \sim \ep^{-{dim}_B (A)},
  \end{equation}
a major feature that makes the box dimension useful for practical applications, as it can be computed through a linear regression in a log-log plot.

The strong similarity (in a geometric setting) between the box dimension and the Kolmogorov scaling function (\ref{scalKolmo}) can now clearly be seen. 
The latter is also defined as a linear regression in log-log plots (log of the $p$-th moments vs. the log of the scales): 
Both objects hence fall under the general heading of { \em scale invariance}. 

This relationship between geometric and analytic quantities is not only formal: On one hand,  the determination of scaling functions has important consequences in the numerical determination of the fractal dimensions of graphs, and related quantities, such  as the $p$-variation, see Section \ref{bdsec} ; on other hand,  the { \em multifractal formalism}, exposed in Section \ref{secmulti}, allows to draw relationships  between scaling functions and  the dimensions of other types of sets, defined in a different way:  
They are the sets of points where the function considered has  a given pointwise smoothness (see Definition \ref{defholdreg}).  
 Note however that the notion of dimension used  in the  context of the multifractal formalism differs from the box dimension: One has to use the { \em Hausdorff dimension} which is now defined (see \cite{Falc1}). 

\BD  
\label{defmeshaus} 
Let   $A$ be a subset of $ \RR^d$. 
If $\ep>0$  and  $\delta \in [0,d]$,  let
\[ M^\delta_{\ep} =\inf_R \;  \left( \sum_{  i} | A_i |^\delta \right) ,\]
where $R$ is an  $\ep$-covering of   $A$, i.e. a covering of  $A$  by bounded sets  $\{ A_i\}_{i \in
\NN}$
 of diameters  $| A_i | \leq \ep$.
(The infimum is therefore taken on all  $\ep$-coverings.)

For any $\delta \in [0,d]$, the 
$\delta$-dimensional  Hausdorff measure of 
$A$ is 
\[  mes_\delta (A) = \displaystyle\lim_{\ep\rightarrow 0}  M^\delta_{\ep}. \]
One can show that there exists
 $\delta_0 \in [0,d]$ such that 
\[ \begin{array}{l} \forall \delta < \delta_0, \;\;\; 
 mes_\delta (A) = + \infty  \\ \forall \delta > \delta_0, \;\;\; 
 mes_\delta (A) = 0
. \end{array}\]
This critical  $\delta_0$ is called the  Hausdorff dimension of  $A$, and is  denoted by $\dim (A)$. 
\ED  

An important convention, in view of the use of these dimensions in the context supplied by the multifractal formalism,  is that, if   $A$ is empty,  then $\dim \,(\emptyset )= -\infty $.

In contradistinction with the box dimension, the Hausdorff dimension is not obtained through a regression  on log-log plots of  quantities  that are computable on real-life data; therefore it is fitted to theoretical questions and can not be used directly  in signal or image processing; however, we will see that the multifractal formalism allows to bypass this problem by supplying an indirect way to compute  such dimensions. \\

\noindent {\bf Scaling, dimensions and dynamics.} \quad The discussions in this section aimed at illustrating the interplay between three different notions: 
The signal processing notion of scaling, or scale invariance, the geometrical notion of dimension and the dynamic notion of constructive split/multiply or split/reproduce procedures.  Note that the relationships between fractals and dynamics have been  investigated by B. Mandelbrot, see \cite{Mandelbrot2004} for a panorama. 
The dynamic constructive procedures of geometrical sets (e.g., Von Koch snowflake) or of functions (e.g., multiplicative cascades) result in geometrical properties of the constructed object.
These geometrical properties can be theoretically characterized using dimensions. 
These dimensions can be measured practically if they additionally imply the existence of power law behaviours (or scaling) (such as those reported in (\ref{equ-scalingBoxdim}) or  (\ref{kolmo})).
When the  constructive procedure underlying the geometrical sets or functions analyzed are known and parametrized, the measure of geometrical dimension may enable to recover the corresponding dynamic parameters.
However this is  certainly no more  true  outside of a parametric model: Many different dynamic constructions may yield the same dimensions.
Therefore, the inverse problem of identifying the nature of the dynamic or even its existence from the sole measure of dimensions is ill-posed and can not be achieved in general. 
Specifically, a concave scaling function, that hence departs from linearity in $p$, does not prove the existence of an underlying cascade mechanism in data. 

\subsection{From graph box-dimensions to $p$-variation: The oscillation scaling function}

\label{bdsec}

\noindent {\bf From graphs to functions.} \quad For data analysis purposes, an obvious difference between the determination of the fractal dimension of a set $A$ and the scaling function $\eta_f (p)$ of  a function $f$ has already been pointed out: While the former reduces to a single number,  the latter consists of a whole function of the variable $p$, hence yields a much richer tool for classification and identification. 
However, the gap can  easily be  bridged in the case where the geometrical set $A$ consists of the graph, denoted $\gr (f)$, of the function $f$.
Let us examine how a slight  extension of the notion of box dimension of a graph allows to introduce a dependency on a variable $p$, and hence a family of exponents whose definition bears some similarity with the Kolmogorov scaling function. \\

\noindent {\bf Oscillations.} \quad Let $f: [0,1]^d \rightarrow \RR$ be a continuous function  (the set $ [0,1]^d$ plays no specific role and is taken here just to fix ideas).  
If $A \subset [0,1]^d$, let $ Os_f ( A )$ denote the \emph{oscillation} of $f$ on the set $A$, i.e. 
\[ Os_f ( A ) = \sup_A f - \inf_A f.  \]
Note that this notion defines only first order  oscillation;   if $A$ is a convex set, second order oscillation (which has to be used for smooth functions) is  defined  as  
\[ Os_f ( A ) = \sup_{x, y \in A} \left|  f(x)  -2 f\left( \frac{x+y}{2}\right) +   f(y)  \right|   \]
(and so on for higher order oscillation). 
Let $ \la_j (x_0)$ denote the dyadic cube  of width $2^{-j} $ that contains $x_0$; $\La_j$ will denote the set of such dyadic cubes  $\la\subset \RR^d$  of width $2^{-j}$.
Let now $p$ be given. 
The $p$-oscillation at scale $j$ of $f$ is defined as 
\BE 
\label{posc}  
R_f(p, j)  = 2^{-dj} \sum_{ \la \in \La_j}   (Os_f (\la))^p. 
\EE
This quantity allows to define the { \em oscillation scaling function} of $f$ as
\BE 
\label{osf} O_f (p) =  \displaystyle \liminf_{j \rightarrow + \infty} \;\; \frac{\log \left(   R_f (p, j) \right) }{\log (2^{-j})}. 
\EE

\noindent {\bf Box dimension of graphs versus the oscillation scaling function.} \quad  Since $f$ is continuous,  the graph of $f$ is a bounded subset of $\RR^{d+1}$. 
First, note that, in the definition  of the upper box dimension (\ref{bd}),  one can rather consider  coverings by dyadic cubes of side $2^{-j}$; indeed each optimal ball  of width $\ep$ used in the covering is included in at most $4^{2d}$ dyadic cubes of width $2^{-j}= [ \log_2 (\ep )]$; and, in the other way, a dyadic cube is contained in the ball of same center and width $\sqrt{d}\cdot  2^{-j}$. 
Therefore using dyadic cubes (instead of optimally positioned balls) only changes prefactors  in  $N_\ep (\gr (f) )$ and   not scaling exponents.  

Let us now consider the dyadic cubes used in the covering of the graph of $f$  and which  stand  above a dyadic cube  $\la$ of dimension $d$ and width $2^{-j}$ included in $[0,1]^d $. 
Their number $N( \la )$ clearly satisfies:
\[ 2^j \left( \sup_\la f -\inf_\la f \right) \leq N( \la ) \leq 2^j \left( \sup_\la f -\inf_\la f \right) +2. \]
This, together with the definition of the oscillation scaling function, implies that for $p=1$, 
\[ \overline{dim}_B (\gr (f))  = \max ( d, d+1 -O_f (1)). \]

However, the oscillation scaling function can be used for other (positive) values of $p$, and therefore  yields a classification tool with similar potentialities as those offered by the Kolmogorov scaling function. We will see an example of its use in finance, in order to estimate quadratic variations, see Section \ref{sec-finance}.
An extension of the oscillation scaling function will be defined later, defined through wavelet quantities: The { \em leader scaling function}, cf. Section \ref{secleade}. \\

\noindent {\bf Box dimension and selfsimilarity exponent.} \quad For some selfsimilar functions and processes, a  direct relationship between the box dimension of their graphs and their selfsimilarity exponent can be drawn.  
For instance,  one-dimensional fBm of selfsimilarity index $H$ has  a graph of box dimension $2-H$, and  
the same result  holds in a deterministic setting for the graphs of the Weierstrass-Mandelbrot functions, which have box dimension $2 - H$ (see \cite{Fal93, Falc97} where similar results can also be found for Hausdorff dimensions of graphs). 
Note that B. Mandelbrot himself had already drawn   connections between selfsimilarity and fractal dimension (cf. e.g., \cite{H22,H23,H21}).

\section{Geometry vs. Analysis }
\label{sec-GeomAnalysis}

The shift from fractal geometry to multifractal analysis essentially consists in replacing the study of selfsimilarity, and of  geometric properties of sets, by the study of analytic properties of functions.  
This shift has several facets:  First, in Section \ref{secpoin}, we will explore   relationships between selfsimilarity and pointwise  regularity, two concepts which are often confused, and actually   have subtle interplays which are far from being fully understood today.
Second,  in Section~\ref {sec-sfvsfs},   we will give  the interpretation of scaling functions in terms of function space regularity, which allows to use the full apparatus of functional analysis. 
Finally, in Section~\ref{secmulti},  the multifractal formalism is introduced: It allows to  interpret the (Legendre transform of the) scaling function  in terms of dimensions of sets of points where $f$ has  a given pointwise H\"older regularity.     

\subsection{Pointwise H\"older regularity}

\label{secpoin}

\noindent {\bf H\"older exponent.} \quad  Selfsimilarity of either deterministic functions or stochastic processes is intimately related with   their local regularity.  
Before explaining these relationships,  let us start by recalling the notion of { \em pointwise  H\"older regularity} (alternative notions of pointwise  regularity,  have also been considered, allowing to deal with functions that are not locally bounded, see \cite{JaffCies,JaffToul,Bergou,JaffMel}). 

 \BD  \label{defholdreg}  Let $f: \; \RR^d  \rightarrow  \RR $  be a { locally   bounded function}, $x_0 \in   \RR^d$ and  let $\gamma \geq 0$; $f $ belongs to $ C^\gamma (x_0)$ if there exist $C>0$, $R >0$  and a polynomial  $P$ of degree less than  $\gamma$ such that  
\[ \mbox{ if} \;\; |x-x_0| \leq R,  \;\; \mbox{ then} \hspace{8mm}   |f(x) -P(x-x_0) | \leq C | x-x_0|^\gamma . \]  
The H\"older exponent of    $f$ at $x_0
$ is  
\[  h_f  (x_0) =\sup \left\{ \gamma : \;\; f \;\; \mbox{ is  }   \;\;   C^{\gamma
} (x_0) \right\} . \] 
\ED

When H\"older exponents lie between $0$ and $1$, the Taylor polynomial  $P(x-x_0)$ boils down to $f(x_0)$ and the definition of the H\"older exponent heuristically means that,
\[  |f(x) -f(x_0) | \sim | x-x_0|^{ h_f(x_0) }. \]
For these cases, a relationship can be drawn with the local oscillation considered in Section \ref{bdsec}. Indeed, let 
 $ \la_j (x_0)$ denote the dyadic cube  of width $2^{-j} $ that contains $x_0$, and $3\la_j (x_0)$ the  cube of same center and 3 times wider (it therefore contains $3^d$ dyadic cubes of width $2^{-j}$). 
One easily checks that, if $ 0 < h_f (x_0 ) <1$, then 
\BE  \label{holoscill} h_f  (x_0) = \displaystyle \liminf_{j \rightarrow + \infty} \;\; \frac{\log \left( Os_f ( 3\la_j (x_0) ) \right) }{\log (2^{-j})}. \EE 
This formula puts into light the fact that the H\"older exponent can theoretically be obtained through linear regressions on log-log plots of the oscillations $ Os_f ( 3\la)$ versus the log of the scale.  
However, in practice, this formula  has found little direct applications for two reasons: 
One is  that alternative formulas based on { \em wavelet leaders}  (cf. Section \ref{secleade}) are numerically more stable, and the other is that many  types of signals display an extremely erratic H\"older exponent, so that  the instability of its determination is intrinsic, no matter which method is used.   
We will see in Section \ref{secleade} how to turn this deeper obstruction. \\

\noindent {\bf Implications of selfsimilarity on  regularity and dimensions.} \quad  The H\"older exponent sometimes coincides with the selfsimilarity exponent, hence resulting in a confusing identification of the two quantities.
This is the case for two frequently used examples:   Weierstrass-Mandelbrot functions and fBm sample paths  have everywhere the same and constant H\"older exponent, $\forall x:\,h_f(x) = H$.   
The theoretical definition of such examples essentially relies on a single parameter, and therefore it comes as no surprise that, for such simple functions with a constant H\"older exponent,  both their selfsimilarity and their local regularity are governed by this sole parameter (and coincide).  Similarly, their scaling functions are linear in $p$, and therefore are governed by  one parameter only, the slope of this linear function, which, for these two examples,   is  $H$.  
 
However, one should not infer hasty conclusions from  these  examples:  Selfsimilar processes with stationary increments do not always possess a single H\"older exponent that corresponds to the selfsimilarity exponent, as  
 shown by the example of { \em stable L\'evy processes}: 
These processes are  selfsimilar  of exponent $H =1/\alpha$, where   $\alpha $ is the stability parameter. 
In addition,  their scaling function is  linear in $p$ only for $p$ small enough:
 \[  
 \begin{array}{rll} 
 \eta_f (p)  & =\displaystyle H{p}  &  \mbox{ if}  \;\; 0 \leq p \leq\frac{1}{H}  \\ 
 &  \\  
 & = 1 &  \mbox{ if}  \;\;  p > \frac{1}{H} . 
\end{array} 
\] This is intimately related with the fact that their H\"older exponent is not constant: 
 Sample paths of a stable L\'evy process   have, for each $h \in [0, H ]$, a dense set of points where their H\"older exponent takes the value $h$ (cf. \cite{JaffLev} for L\'evy processes, and \cite{DurJaff} for { \em L\'evy sheets} indexed by $\RR^d$). 
Thus,  their H\"older exponents  jump in a very erratic manner from point to point, as a consequence of the high variability of the increments: 
Indeed, their distributions only have power-law decay, and therefore (in strong contradistinction with Gaussian variables) can take extremely large values with a relatively large probability; the importance of this  phenomenon in applications  was largely undervalued in the past, until B. Mandelbrot pointed out  its  possibly dramatic consequences, referring to them in a  colorful way as \emph{Noah's effect}. 

Let us now be more explicit concerning the relationships between selfsimilarity and H\"older regularity. A first result relates selfsimilarity with uniform regularity of sample paths. 

\BP  \label{propChents} Assume that $X_t$ is  a selfsimilar process of exponent $H$, with  stationary
increments, and  such that 
$ \EEE ( | X_1|^\alpha ) < \infty $ for some $\al >0$.
Then,  
\BE \label{unifregselfsi} \forall h  < H-{1}/{\alpha},  \;\; \mbox{ a.s. } \;   \exists C, \;\;  \forall  t, s , \hspace{1cm}  | X_t -X_s |   \leq C | t-s|^{h} \EE
\EP

The proof  follows directly from the Kolmogorov-Chentsov theorem: Recall that this theorem  asserts that, if $X_t$ is  a random process  satisfying 
\[\forall  t, s , \hspace{1cm}  \EEE ( | X_t -X_s |^\al  ) \leq C | t-s|^{1+\beta} \]
then, 
\[ \forall h  <\beta /{\alpha},  \; 
 \mbox{ a.s. } \;   \exists C, \;\;  \forall  t, s , \hspace{5mm}  | X_t -X_s |   \leq C | t-s|^{h} \] 
 (see \cite{CohIst}). But, if $X$ is selfsimilar, then 
\BE  \EEE ( | X_t -X_s |^\al  ) = | t-s|^{H \alpha}  \EEE( | X_1|^\alpha ) , \EE
and the proposition follows.  \\Ê

 We will reinterpret this result in Section \ref{secunif}   when the uniform H\"older exponent will be introduced.
The following  result  goes in the opposite direction: It relates selfsimilarity with irregularity of  sample paths.

\BP \label{propirreg}  Let  $X_t$  be   a selfsimilar process of exponent $H$, with  stationary
increments,  and which   satisfies $\PP ( \{ X_1 =0 \} ) = 0 $, i.e. the law of $X$  at time $t=1$ (or at any other time $t >0$) does not contain a Dirac mass at the origin. 
Then 
\BE \label{irregselfsi} \mbox{ a.s. } \;\; \mbox{ a.e. } \hspace{5mm}  h_X (t) \leq H . \EE 
\EP 

{ \bf Proof:}  We first prove the result for $t =0$.  The assumption $\PP ( \{ X_1 =0 \}Ê ) = 0 $   implies the existence of  a decreasing sequence $(\delta_n)_{n \in \NN} $ which converges to 0 and satisfies 
\[  \PP ( |X_1|  \leq \delta_n  )  \leq \frac{1}{n^2} . \]
Let $l_n = \exp (-1/\delta_n )$. Then $ X_{l_n} \;  \stackrel{{\cal L} }{=}  \; (l_n)^H X_1 $. Therefore
\[ \PP \left( |X_{l_n}|  \leq \frac{(l_n)^H }{\log (1/l_n) }  \right)  = \PP ( |X_1|  \leq \delta_n  ) \leq \frac{1}{n^2} . \]
The Borel-Cantelli lemma implies that a.s.  a finite number of these events happens, and therefore $h_X (0) \leq H$. 
By the stationarity of increments, this argument holds the same at each $t$, i.e. 
\[ \forall t \;\; \mbox{ a.s.} \hspace{5mm}  h_X (t) \leq H ; \] and Fubini's theorem  allows to conclude that (\ref{irregselfsi}) holds. \\

 In particular, assume that $X_t$ satisfies the assumptions of Proposition \ref{propirreg}, and also of Proposition \ref{propChents} for any $\al >0$. Then, both results put together show  that, in this case, (\ref{unifregselfsi}) is optimal. 
Relationships between   selfsimilarity  and   H\"older exponents for selfsimilar processes will be  further investigated in Section~\ref{sec-mfMF}.

Note that, without additional assumptions,  one can not expect results  more precise  than those given by Propositions  \ref{propChents} and \ref{propirreg}, as shown by the subcase of stable, selfsimilar processes with  stationary
increments, see \cite{MaejKon,Samorodnitsky1994}: Let  $H$ denote the selfsimilarity index of such a process, and $\al $ its stability index. 
First, note that there is no direct relationship between $\al$, $H$ and H\"older regularity:  
They always satisfy 
\[ 0 < \al <2 \;\; \mbox{  and } \;\; H \leq \max (1, 1/ \al); \]  the case of   fBm is misleadingly simple since it satisfies 
\[ \al =2  \;\; \mbox{  and } \;\; \forall x, \;\;  h_{B_H} (x) = H \in (0, 1);  \]
however L\' evy processes satisfy $H = 1/ \al$ and $ h_f (x) $ takes all values in $ [0, H] $. 
Finally, L\'evy fractional stable processes   have nowhere locally bounded stable sample paths if $H < 1/ \al$, and  continuous sample paths if  $H > 1/ \al$ (see Section \ref{secunif} for a more precise result). 
Sample paths of selfsimilar processes (with stationary increments) together with their scaling functions (and multifractal spectra) are illustrated in Fig.~\ref{fig-wlmfSS} in Section~\ref{secspec}. \\Ê

Let us now mention a few results concerning the relationship between  selfsimilarity and fractional dimensions. We start with fBm which again is  the simplest  and best understood case. Such results are of two types:
\begin{itemize}  \item Dimensions of the graph of  sample paths:  Box and Hausdorff dimensions coincide and take the value $2-H$, see  \cite{Fal93, Falc97}. 
\item Relationship between the Hausdorff dimension of a set $A$ and of its image $X (A)$ (see see \cite{ShiehXiao} and references therein)
\BE \label{trandormdim} \mbox{ a.s.,    for any Borel set, } \;  A \subset \RR,  \hspace{6mm} \dim X( A) = \min (1, \dim (A)/H  ). \EE
\end{itemize}
We will see that this  second result has important consequences for the multifractal analysis of some finance models proposed by  L. Calvet, A. Fisher, and B. Mandelbrot, see (\ref{equ-fbmmf}). 
Partial results of these two types also hold for other selfsimilar processes, see \cite{ShiehXiao}.   \\

\noindent {\bf Regularity versus dimension: The mass distribution principle.} \quad The example of stable L\'evy processes  suggests that the most fruitful connection  to be drawn does not involve the notion of selfsimilarity, but rather should connect scaling functions with pointwise regularity. 
 This is precisely what the \emph{multifractal formalism} performs, as it establishes a relationship between scaling functions and the Hausdorff  dimensions of the 
 { \em isoh\"older sets}  of $f$ defined by  
\BE \label{iso}  E_f(h) = \{ x_0: \hspace{5mm}  h_f(x_0) = h\} .  \EE 
The multifractal formalism will be detailed in Section~\ref{secmulti}.
A clue which  shows that such relationships  are natural can however be immediately given, via the { \em mass distribution principle} which goes back to Besicovitch (see Proposition \ref{distmas} below).  
This principle establishes a connection between the dimension of a set $A$ and the regularity of the measures that this set can carry. 
The motivation for it stems from the following remark.  
Mathematically, it is easy to obtain upper bounds for the Hausdorff dimension of a set $A$: 
Indeed, any  well chosen sequence of particular   $\ep$-coverings (for a sequence  $\ep_n \rightarrow 0$) yields an upper bound.  
To the opposite, obtaining a lower bound by a direct application of the definition is more difficult, since it requires to consider { \em all possible } $\ep$-coverings of $A$.   
The advantage of the {mass distribution principle} is that it allows to obtain lower bounds of Hausdorff dimensions by just constructing one measure { \em carried}  by the set $A$. 
Recall that a  measure  $\mu$ is {carried by} $A$ if $\mu (A^c)= 0$. 
Note, however, that a set $A$ satisfying  such a definition  is not unique; In particular, this notion should not be mistaken with the more usual one of { \em support} of the measure, which is uniquely defined as the intersection of all closed sets carrying $\mu$.  For example, consider the measure 
\[ \mu = \sum_{p/q \in [0,1], p\wedge q =1} \frac{1}{q^3}\delta_{p/q} \]
(where $\delta_a$ denotes the Dirac mass at the point $a$). The support  of this measure 
is the  interval  $[0,1]$, but it is (also) carried by the (much smaller) set  $A= \QQ \cap [0,1]$. 

\BP \label{distmas} Let  $\mu $  be a probability measure carried by  $A$.  Assume that there exists  $\delta \in [0, d]$,  $C>0$  and  $\ep >0$  such that, for any ball 
$B$  of diameter less than  $\ep$,  $\mu $ satisfies the following uniform regularity assumption
\[  \mu (B) \leq C |B|^\delta  . \]
 Then the $\delta$-dimensional measure of $A$ satisfies:   $mes_{\delta }  (A) \geq {1}/{C} $, and therefore $\dim (A) \geq \delta$.
\EP

{ \bf Proof:}   Let  $\{ B_i\}_{i \in \NN}$ be an arbitrary  $\ep$-covering of  $A$.  Then
\[1 = \mu (A) =  \mu  \left(\bigcup B_i\right) \leq \sum  \mu  (B_i) \leq C \sum |B_i|^\delta  . \]
The result follows by passing to the limit  when  $\ep\rightarrow 0$. 

\subsection{Scaling functions vs. function spaces} 
\label{sec-sfvsfs}

\noindent {\bf From scaling functions to function spaces.} \quad  Kolmogorov scaling function $\eta_f(p)$ 
(as defined by (\ref{scalKolmo}))  receives a
function space interpretation,  which  is important for several reasons.  
On one hand,  it allows to derive a number of its  mathematical properties, and on  other hand, it points towards variants and extensions of this scaling function, yielding sharper  information on the singularities existing in data, and therefore offering a deeper understanding of the multifractal formalism.  \\
 


\noindent {\bf Function space interpretation.} \quad  The  function space interpretation of the Kolmogorov scaling function can be  obtained  through the use of  the spaces $\mbox{Lip}  ( s , L^p ) $ defined as follows. 
 \BD \label{esplip}
 Let $s  \in (0,1)$, and  $p \in [1, \infty]$; 
$ f \in \mbox{Lip}  ( s , L^p  (\RR^d)) $ if $\int | f(x) |^p dx < \infty $ and 
\begin{equation} \label{nicol1}  \exists C>0, \;\; \forall h  >0,  \hspace{6mm} 
\int |  f(x+h) -f(x) |^p  dx  \leq C |h|^{sp}  .  \end{equation}
 \ED
Note that  larger smoothness indices $s$  are  reached by replacing in the definition the first-order difference 
 $|  f(x+h) -f(x) |$ by higher and higher order differences:  $|  f(x+2h) -2   f(x+h) +f(x) | $,... (which is coherent with the remark just after the introduction of  Kolmogorov scaling function (\ref{scalKolmo}), where, there too, higher order differences have to be introduced in order to deal with  smooth functions $f$).  
 
It follows from (\ref{kolmo}) and  (\ref{nicol1}) that, 
\begin{equation} \label{nicol}  \forall p \geq 1, \hspace{6mm}  \eta_f (p) = \sup \{ s: f \in \mbox{Lip}  ( s/p , L^p  (\RR^d)) \} .
\end{equation}
In other words, the scaling function allows to determine within which spaces $\mbox{Lip}  ( s , L^p ) $ data are contained, for   $p \in [1, \infty]$.   
The reformulation supplied by (\ref{nicol}) has several advantages:
It will lead to  a better numerical  implementation, based on wavelet coefficients, and  it    will have the additional advantage of  yielding a natural extension of these function spaces to  $p \in (0,1)$; this, in turn, will  lead to a scaling function also defined for $p \in (0,1)$,  therefore supplying  a  more powerful tool  for classification, see Section \ref{secwavcoe}. Note that the reason why Kolmogorov scaling function cannot be directly extended to $p <1$ is that  $L^p$ spaces and the spaces which are derived from them (such as  $\mbox{Lip}  ( s , L^p ) $)  do not make sense for  $p <1$; indeed Definition \ref{esplip} for $p <1$ leads to mathematical inconsistencies (spaces of functions thus defined would not necessarily be included in spaces of distributions).

Another advantage of the function space interpretation of the scaling function is that it allows to draw relationships with areas of modeling where the a priori assumptions which are made are function space regularity assumptions.  
Let us mention two famous examples, which will be further used as illustrations  for the wavelet methods developed in the next section: { \em Bounded Variations}  (BV) modeling in image processing (cf. Section~\ref{secwavcoe}), and bounded quadratic variation hypothesis in finance (cf. Sections \ref{secleade} and \ref{sec-finance}, this latter case being related with the determination of the oscillation scaling functions (\ref{osf})). \\

\noindent {\bf Bounded variations and image processing.} \quad A function $f$ belongs to the space BV, i.e., has { \em bounded variation}, if its gradient, taken in the sense of distributions,  is a  finite (signed) measure. 
A standard assumption in image processing is that real-world images  can be modeled  as the sum of a function $u \in BV$ which models the \emph{cartoon part}, and  another term $v$ which accounts for the noise and texture parts (for instance, the first ``$u+v$ model'', introduced by Rudin, Osher and Fatemi in 1992 (\cite{ro92}) assume that  $v \in L^2$, see also \cite{H22,H23}).  
The BV model is motivated by the fact that if an image is composed of smooth parts separated by {contours}, which are  piecewise smooth curves, then its gradient will be the sum of  a smooth function (the gradient of the image inside the smooth parts)  and Dirac masses along the edges, which are typical finite measures.  On the opposite, characteristic functions of domains with fractal boundaries usually do not belong to BV, see Fig.~\ref{fig-BV}  for an illustration).  Therefore, a natural question in order to validate such models is to determine whether an image (or a portion an image)   actually belongs to the space BV or not. 
We will see in Section  \ref{secwavcoe} how this question can be given a sharp answer using  a direct  extension of  Kolmogorov scaling function. \\
 
\noindent {\bf Bounded quadratic variations.} \quad  Another motivation for function space modeling  is supplied  by the \emph{finite quadratic variation} hypothesis in finance.  
In contradistinction with the previous image processing example, this hypothesis is not deduced from a plausible intrinsic property of financial data, but rather constitutes an ad hoc regularity hypothesis which (considering the actual state of the art in stochastic analysis) is required  in order to use the tools of stochastic integration. 
There exist several slightly different formulations of this hypothesis.  
One on them, which we consider here, is that the $2$-oscillation is uniformly bounded at all scales.  
This means that $f$, assumed to be defined on $[0,1]^d$, satisfies
\[ \exists C, \;\; \forall j  \geq 0, \hspace{1cm}  2^{dj} R_f (2, j)  = \sum_{ \la \in \La_j}   (Os_f (\la))^2 \leq C   \]
(where $R_f (p, j)$ was defined by (\ref{posc})).
The definition of the oscillation scaling function (\ref{osf}) shows that, if $O_f (2) >d$, then $f$ has finite quadratic variation, and if $O_f (2) <d$, then this hypothesis is not valid.  
This will be further discussed in Section~\ref{secwavcoe} and illustrated in Section~\ref{sec-finance} and Fig.~\ref{fig-finance}.

\section{Wavelets: A natural tool to measure global scale invariance} 
\label{secwav}

The general considerations developed in Section \ref{sec-sfvsfs} emphasized the importance of scaling functions as a data analysis tool. 
To further develop this program, reformulations of scaling functions which are numerically more robust than the Kolomogorov and oscillation scaling functions introduced so far are of both practical and theoretical interests.  Motivations for introducing wavelet techniques for this purpose were already mentioned in Section \ref {secwavmot}. 
We now introduce these  techniques in a more detailed way, and develop and explain the benefits  of rewriting scaling functions in terms of wavelet coefficients. 

\subsection{Wavelet bases}
\label{secwavedis}

Orthonormal wavelet bases on $\RR^d$ exist under two slightly different forms.  
First, one can use $2^d-1$ functions $\psi^{(i)}$ with the following properties: The functions 
\BE
\label{wavbasbad}  
2^{dj/2}\psi^{(i)}(2^jx-k), \hspace{5mm} \mbox{ for} \hspace{3mm} i = 1, \cdots 2^d -1, \hspace{3mm}   j\in \ZZ \hspace{3mm} \mbox{ and} \hspace{3mm}
k\in\ZZ^d 
\EE
 form an orthonormal basis of $L^2 (\RR^d)$.
 Therefore, $\forall f \in L^2$, 
\begin{equation} \label{ecrit22} f(x) =    \sum_{j\in \ZZ} \sum_{k\in \ZZ^d}
\sum_{i} \cjk \psi^{(i)} (2^jx-k), \end{equation}  where the $\cjk$ are the wavelet coefficients of 
$f$, and are given by  
\begin{equation} \label{cjk}  \cjk = 2^{dj} \int_{\RR^d} f(x)  \psi^{(i)} (2^j x -k) dx . \end{equation} 
Note that, in practice, discrete wavelet transform coefficients are generally not computed through the integrals (\ref{cjk}), but instead using the discrete time convolution filter bank based pyramidal and recursive algorithm (cf. \cite{Mallat1998}), referred to as the  {\em  Fast Wavelet Transform} (FWT). 

However,  the use of such  bases rises several difficulties as soon as one does not have the a priori information that $f \in L^2$. 
For instance, if the only assumption available is that $f$ is a continuous and bounded function, then one can still compute wavelet coefficients of $f$ using (\ref{cjk}) (indeed, these integrals still make sense), however, the series at the right-hand side of (\ref{ecrit22})  may converge to a function which  differs from $f$. 
This is due to the fact that  certain \emph{special} functions (which do not belong to $L^2$) have all their wavelet coefficients that vanish, and therefore, such possible components of functions will always be missed in the reconstruction formula (\ref{ecrit22}).  
These  \emph{special}  functions always include polynomials, however, in the case of Daubechies compactly supported wavelets, other \emph{fractal-type} functions also have this property, see \cite{Lemar}.  
This explains why, at the end of Section \ref{intro},  we mentioned that the  properties of wavelet coefficients of certain processes that  we obtained were not characterizations: Indeed information on these special functions, as possible components of the processes, is systematically  missed, whatever the information on the wavelet coefficients is.  
However, this drawback  can be turned by using a slightly different basis, which we now describe.   \\Ê

The alternative wavelet bases that will systematically be used from now on are of the following form: 
There exists a function  $\varphi (x)$ and  $2^d-1$ functions $\psi^{(i)}$ with the following properties: The functions 
  \BE \label{seizebis} \mbox{ $\varphi (x-k)$ (for $k\in \ZZ^d ) \hspace{5mm} $  and   $\hspace{5mm} 2^{dj/2}\psi^{(i)}(2^jx-k)$
( for $k\in\ZZ^d,$  and $j\geq 0 $)} \EE
 form an orthonormal basis of $L^2 (\RR^d)$.  
In practice, we will use Daubechies compactly supported wavelets \cite{Daubechies1988}, which can be chosen arbitrarily smooth. 
Since these functions form an orthonormal basis of $L^2$, it follows, as previously, that
\begin{equation} \label{ecrit2}\forall f \in L^2, \hspace{6mm}  f(x) =  \sum_{k\in \ZZ^d}  C_k \varphi (x-k) +  \sum_{j=0}^{\infty} \sum_{k\in \ZZ^d}
\sum_{i} \cjk \psi^{(i)} (2^jx-k); \end{equation} the $\cjk$ are still given  by (\ref{cjk}) and
\begin{equation} \label{ck} C_k = \int_{\RR^d} f(x)  \varphi ( x -k) dx .\end{equation} 

This choice of basis enables to answer the following important data analysis questions (which was not permitted by the  previous choice (\ref{wavbasbad})):  In real-world data, the a priori assumption that $f \in L^2$ does not need to hold (for instance, the sample paths of  standard models, such as Brownian motion,  do not belong to $L^2 (\RR) $).  
As already mentioned, the sole information that the series of coefficients of a function $f$ on an orthonormal basis is square summable, does not imply that $f \in L^2$ (consider the example of $f=1$ and the basis (\ref{wavbasbad})). 
However, for the alternative basis given by (\ref{seizebis}),  the  following property  now holds: Assume that  $f$ is a function in $L^1_{loc}$ with slow growth, i.e. satisfies 
 \BE \label{lunsg} \exists  C, A >0: \hspace{2cm} \int_{B(0,R)} | f(x) | dx \leq C(1+R)^A,\EE
where $B(x, R)$ denotes the ball of center $x$ and radius $R$;  
then the wavelet expansion (\ref{ecrit2})  of $f$ converges to $f$ almost everywhere.  
Additionally, if the wavelet coefficients of $f$ are square summable, then $f\in L^2$
(in contradistinction with what happens when using the basis (\ref{wavbasbad})). 
 Note that the slow growth assumption is a very mild one, since it is satisfied by all standard models in signal processing, and actually is necessary, in order  to  remain in the mathematical  setting supplied by  tempered distributions.

Note that  one can even go further:  (\ref{cjk}) and (\ref{ck}) make sense even if $f$ is not  a function; indeed,  if one uses smooth enough wavelets, these formulas 
can be interpreted as a duality product between  smooth functions (the wavelets) and tempered  distributions. 
This rises the problem of determining if the series (\ref{ecrit2}) converges in other  functional settings, and it is indeed one of the most remarkable properties of wavelet expansions that it is very often the case. Wavelet characterizations of   function spaces   are detailed in \cite{Mey90I}. Let us  mention   a property   which is particularly  useful, and shows that  convergence properties of wavelet expansions are ``as good as possible'': Suppose that $f$ is  continuous with slow growth, i.e. that it satisfies 
\BE \label{sg1}   \exists C, N >0, \;\;  \forall x \in \RR^d, \hspace{1cm} 
\left|\,f(x)\,\right| \leq C (1+  |x |)^N;
\EE 
then the wavelet expansion of $f$ converges  uniformly on every compact set. 

However, while the differences that we pointed out between the two types of wavelet bases are important for reconstruction, they are not  for the type of analysis that we will perform: Indeed, properties will be deduced from  the inspection of asymptotic behaviors  of wavelet coefficients at small scale, where both bases (\ref{wavbasbad}) and (\ref{seizebis}) coincide. \\

 As an example of  the fruitful  interplay between  the algorithmic structure of wavelet bases and the concept of selfsimilarity, we now prove (\ref{samelaw}).  We start by recalling the definition of equality in law of stochastic processes.  
 
 \BD \label{defeuqalaw}  
 Two vectors of
$\RR^l$: 
 $ (X_{1}, \cdots X_{l})$  and $ 
(Y_{1}, \cdots Y_{l}) $
 share the same law if, for any Borel set $A\subset \RR^l$, $\PP (\{ X \in A\} )  = \PP (\{ Y \in A\} )$. 
 
Two processes
$X_t$ and $Y_t$  share the same law if,   $\forall l \geq 1$, for
any   finite set of time points $t_1, \cdots t_l$,  the vectors of
$\RR^l$
 \[ (X_{t_1}, \cdots X_{t_l}) \hspace{1cm} \mbox{ and} \hspace{1cm}  
(Y_{t_1}, \cdots Y_{t_l}) \]
 share the same law. 
 \ED
 
 We apply this definition to the two  processes $X_{at}$  and 
$a^H X_t$  (such as stated in (\ref{autosimstoch})), which are assumed to share the same law.  
 Taking for $A$ the domain between two parallel  hyperplanes, this, in turns, implies that any finite linear
combinations
  \[ \sum_{i=1}^l \al_i X_{t_i} \hspace{6mm} \mbox{ and} \hspace{6mm}  
\sum_{i=1}^l \al_i Y_{t_i} \]
   share the same law.  Therefore, if $X_t$ is a selfsimilar process of
exponent $H$, then (\ref{autosimstoch}) means that 
   \BE \label{samelawbis} \forall t_1, \cdots t_l,  \hspace{1cm}  \sum_{i=1}^l \al_i X_{ a
t_i} \hspace{6mm} \mbox{ and} \hspace{6mm}  \sum_{i=1}^l \al_i a^H X_{ 
t_i}  \hspace{1cm}  \mbox{ share the same law.} \EE
      Assume now that for almost every $\omega \in \Omega$, the sample
path of $X_t$ is Riemann integrable. The
    coefficient
    \[ c_{j,k} = \int 2^j  X_t \psi (2^j t-k) dt \]
    is a limit of the Riemann sums
    \BE \label{riemann} \sum_i ( t_{i+1}-t_i)  2^j X_{t_i} \psi (2^j t_i
-k) ;\EE
    using (\ref{samelawbis}),  it follows that  (\ref{riemann}) has the same
law as
    \[ \sum_i ( t_{i+1}-t_i)  2^{-H} 2^j X_{ 2 t_i} \psi (2^j t_i -k)  \]
    which, writing $u_i = 2 t_i$, can be written
    \[ \sum_i ( u_{i+1}-u_i)  2^{-H} 2^{j-1} X_{  u_i} \psi (2^{j-1} u_i
-k) , \]
    which is a Riemann sum of the integral  \[  \int 2^{-H} 2^{j-1}    X_u
\psi (2^{j-1} u-k) du  = 2^{-H} c_{j-1, k} .\]
    Passing to the limit when the  largest spacing between sampling points $t_i$ (and therefore $u_i$) 
tends to 0, we obtain that
  \[   c_{j,k}  \stackrel{{\cal L} }{=}  2^{-H}  c_{j-1,k} .  \]
The argument that we developed  for one coefficient can be reproduced similarly for an arbitrary vector of coefficients, hence (\ref{samelaw}) holds.

\subsection{Wavelet Scaling Function}
\label{secwavcoe}

We now investigate properties which arise as consequences of the  wavelet characterizations of function spaces. \\Ê

\noindent {\bf Notations.} \quad Compact notations are used for  indexing  wavelets.   
Instead of the three indices $(i,j,k)$,  dyadic cubes are used:
$ \la\; (= \la (j,k)) \; = \displaystyle\frac{k}{2^j} + \left[ 0, \displaystyle\frac{1}{2^{j}}\right)^d$. 
In order to keep notations as simple as possible, the index $i$ is dropped in formulas and notations, with the implicit convention that sums or suprema over indices are also taken on the index $i$.
Accordingly, $\cla = \cjk$ and $\pla (x)= \psi^{(i)}(2^jx-k) $.
The wavelet $\pla$  is  localized near the cube $\la$: 
Since we use Daubechies compactly supported  wavelets,  
$ \exists C >0 $ such that $ \forall i,j,k,$ $ supp \left( \pla \right)  \subset  C\cdot \la $ 
(where $C \cdot \la$ denotes the cube of same center as $\la$ and $C$ times wider). 
Finally,  let $\Laj$  denote the set of dyadic cubes  $\la$  which index a wavelet of scale  $j$, i.e. wavelets 
 of the form  $\pla (x)= \psi^{(i)}(2^jx-k) $  (note that  $\Laj$ is a subset of the dyadic cubes of size 
 $2^{j+1}$), and $\La$ will denote the union of the $\La_j$ for $j \geq 0$.  \\Ê
 
 \noindent {\bf Wavelet scaling function.} \quad 
A key property of wavelet expansions  is that many function spaces have a simple characterization by conditions bearing on  wavelet coefficients.  
This property has a direct consequence on the determination of the Kolmogorov scaling function. 
Let
 \begin{equation}
 \label{equ-WSF}
 S_f (p,j) = 2^{-dj} \displaystyle\sum_{ \la \in \La_j }  | c_\la |^p.   
 \end{equation}
Classical embeddings between function spaces imply  that, if the wavelets are smooth enough, then    the Kolmogorov scaling function can be re-expressed as (cf. \cite{jmf2})
  \BE \label{defscalond}  \forall p \geq 1, \hspace{6mm} 
  \eta_f (p) =   \displaystyle \liminf_{j \rightarrow + \infty} \;\; \frac{\log \left( S_f (p,j)  \right) }{\log (2^{-j})}, \EE 
This formulation, which, again, yields the scaling function through linear regressions in log-log plots,  has several advantages when compared to the earlier version (\ref{scalKolmo}). 
First, it allows to extend the Kolmogorov scaling function  to the range  $0 < p \leq 1$; it will be shown in Section~\ref{secleade}  how to go even further and define scaling functions for negative $p$s (see also \cite{JAFFARD:2006:A}  for an extension of $\eta_f (p)$ to $p <0$). 
We will call this extension of the Kolmogorov scaling function  the \emph{wavelet scaling function}, 
but  we will keep the same notation.    By construction, $\eta_f (p) $ is a concave increasing function, of derivative at most $d$, see   \cite{Jaffard2006}. \\Ê

\noindent {\bf Wavelet regularity.} \quad Let us say a few words about the required smoothness of wavelets. The rule of thumb is that wavelets should be smoother and have more vanishing moments than  the regularity index appearing in the definition of the function space.  In signal processing, one is not necessarily aware beforehand of the regularity of the data. In practice, one uses smoother and smoother wavelets, and when the results   do not vary, this means that smooth enough wavelets are used.  (Note that this is reminiscent  of the definition of the Kolmogorov scaling function where, in theory, one has  to use higher and higher order differences and stop when the results are numerically stabilized.) In the following we will never specify the required smoothness, and always assume that { \em smooth enough wavelets} are used (the minimal regularity required being always  easy to infer). 

Note that, in the wavelet setting, one should also draw a difference  between  (\ref{defscalond}), which allows to define the  wavelet scaling function of any function $f$,  and its use in data analysis, where scaling properties need to hold to enable measurements based on linear regressions in log-log plots.  \\
 
\noindent {\bf Function space interpretation.} \quad Let us examine how the wavelet scaling function can be used in order to validate the function space assumptions  discussed in Section~\ref{sec-sfvsfs}.
For instance, regarding the $u+v$ model, where $u\in BV $ and $v \in L^2$, the values taken by wavelet scaling function at $p=1$ and $p= 2 $ allow practitioners to determine if data belong to $ BV $ or $ L^2$:
 \begin{itemize}
\item  If {  $\eta_f (1) >1$}, then $f \in BV$, and   if {  $\eta_f (1) <1$}, then $f \notin BV$
\item  If  { $\eta_f (2) >0$}, then $f \in L^2$ and  if  { $\eta_f (2) <0$}, then $f \notin L^2$.
\end{itemize}
Thus wavelet techniques  allow  to check if the function space assumptions which are made, e.g., in certain denoising algorithms relying on then $u+v$ model, are valid. 

\begin{figure}[h]
\centering
\includegraphics[width=0.8\linewidth]{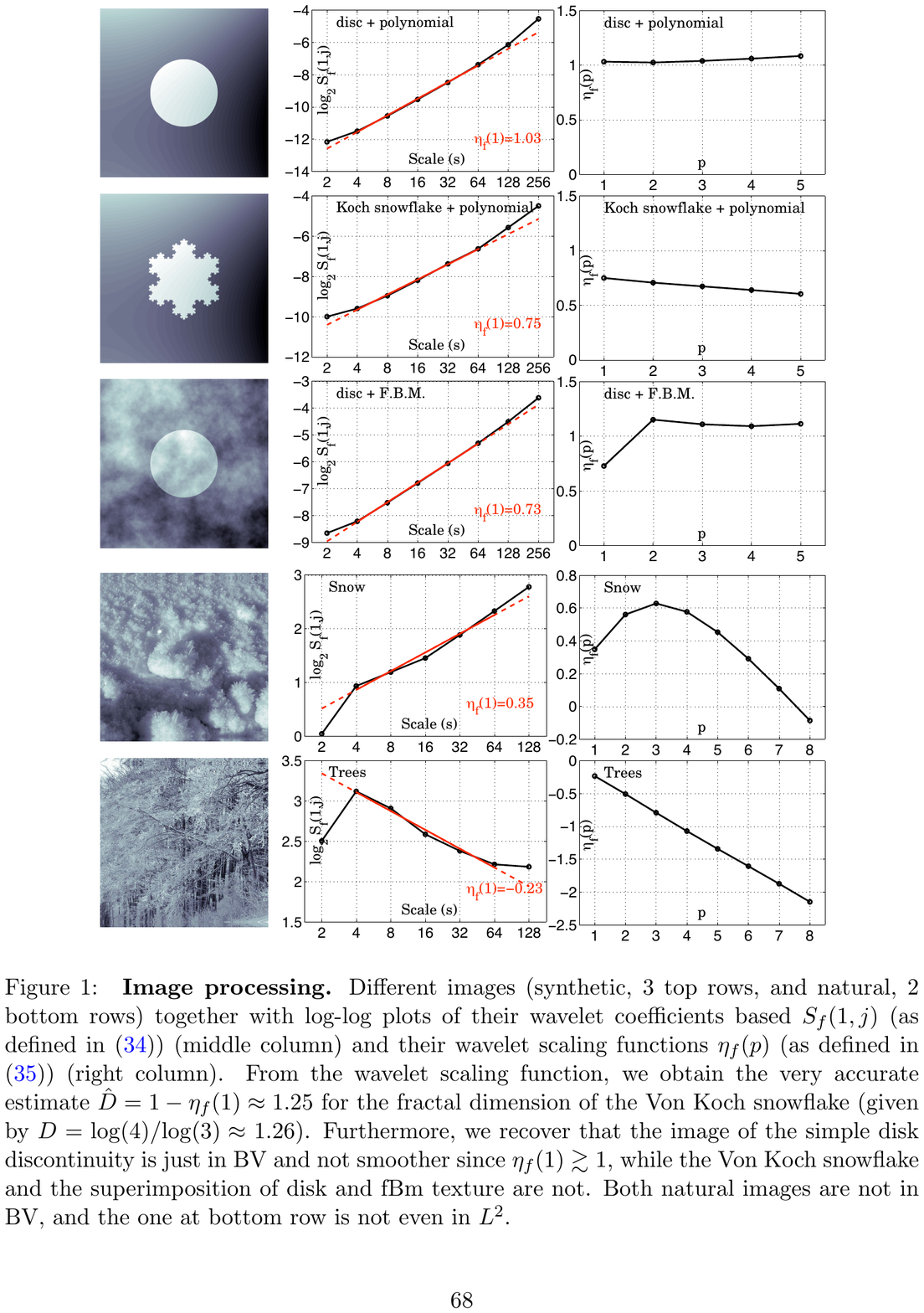}
\caption{\label{fig-BV} {\bf Image processing.} Different images (synthetic, $3$ top rows, and natural, $2$ bottom rows) together with log-log plots of their wavelet coefficients based $S_f(1,j)$ (as defined in (\ref{equ-WSF})) (middle column) and their wavelet scaling functions $\eta_f(p)$ (as defined in (\ref{defscalond})) (right column). From the wavelet scaling function, we obtain the very accurate estimate $\hat D = 1-\eta_f(1) \approx 1.25$ for the fractal dimension of the Von Koch snowflake (given by $D={\log (4)}/{\log (3)}\approx1.26$). 
Furthermore, we recover that the image of the simple disk discontinuity is just in BV  and not smoother since $\eta_f(1)\gtrsim1$, while the Von Koch snowflake and the superimposition of disk and fBm texture are not. 
Both natural images are not in BV, and the one at bottom row is not even in $L^2$.
}
\end{figure}

Examples of synthetic and natural images are shown in Fig.~\ref{fig-BV}, together with the corresponding measures of $\eta_f (1)$ and $\eta_f (2)$. 
The image consisting of a simple discontinuity  along a circle and no texture, (i.e.,  a typical  \emph{cartoon} part of the image in the $u+v$ decomposition)  is in BV, while the image of textures or discontinuities existing on a complicated support (such as the Von Koch snowflake) are not.  Y. Gousseau and J.-M. Morel were the first authors to rise the question of finding statistical tests to verify if natural images 
belong to  BV \cite{GoussMor}. Our results confirm that  the BV requirement is seldom met;   actually, images  often turn out not to be even in  $L^2$, as illustrated in Fig.~\ref{fig-BV} (bottom row). 

\subsection{Uniform H\"older exponent}
\label{secunif} 

 It has been shown in the previous section that the Kolmogorov scaling function can be rewritten and extended as a wavelet coefficient based version. 
It is hence natural to ask whether the same can be performed for the oscillation scaling function (\ref{osf}). 
Because oscillations are defined only if $f$ is locally bounded, this condition needs to be checked in practice.
This test can be performed using  the \emph{uniform H\"older exponent}, which constitutes the subject of this subsection. 
However, before entering technicalities, the following general question concerning function space modeling needs to be addressed.  \\Ê
 
\noindent {\bf Function space modeling.} \quad The issue of determining whether experimental data, acquired by a physical device, can be modeled by a bounded function, or not, may appear meaningless at first sight.
Indeed, data can be collected only with finite size and resolution, and hence consist of finite, and thus  bounded, sequences. 
For images, a common pitfall is to conclude that images are  grey-levels, so that   they necessarily take values in $[0,1]$, and, therefore, by construction,   are appropriately modeled by bounded functions. 
The issue is even more general regarding the whole strategy of functional space modeling:   Given any finite resolution and size sequence of data practically available, they can always be extrapolated by a $C^\infty$ function, that thus belongs to any possible function space. 
The resolution of this paradox requires the use of the notion of scale invariance.   
Let us consider again the example of image modeling and address the issue of deciding whether a digital image can be modeled  by a bounded function or not:   
It is clear that if the image is analyzed only at its finest scale, then the answer is positive. 
However, following the point of view that we developed until now, the problem can be reinterpreted in terms of  analysis of  scaling properties, over a range  of available scales,
and  will therefore amount to determine whether an observed scaling exponent fits with  those that are compatible with  bounded functions, or not. 
The uniform H\"older exponent, which will be considered in this section, answers this question.

More generally, numerical results, obtained from finite resolution and size data, leading to conclude that data belong to certain function spaces and not to others, should actually be understood in the following way: 
The function space regularity thus determined holds under the hypothesis that scaling that are observed on the range of  available scales would continue to hold at finer scales if they were available. \\
 
\noindent {\bf Uniform H\"older exponent: Definition.} \quad The mathematical idea behind testing whether data are locally bounded or not is to consider boundedness as a particular case in  the whole range of Lipschitz spaces, which corresponds to the regularity exponent $s=0$. 
Let us be more specific. 
  
The Lipschitz spaces $C^s_{sg} (\RR^d ) $ are defined for $0 < s <1$  by the conditions : 
 \[  \exists C, N,   \;\; \;\; \forall x,y \in \RR^d, \hspace{1cm} | f(x) | \leq C(1+ |x|)^N. \]
 and 
\[  \exists C , N, \;\;  \;\; \forall x,y \in \RR^d,  \hspace{1cm} | f(x)-f(y) | \leq C |x-y|^s(1+ |x| + |y|)^N. \]
If $s >1$, they are then defined by recursion on $[s]$ by the condition:
$ f \in C^s_{sg} (\RR^d )$ if $f\in L^\infty$ and  if all  its partial derivatives (taken in the sense of distributions) ${\partial f}/{\partial x_i} $ (for $i =1, \cdots d$) belong to 
$C^{s-1}_{sg} (\RR^d )$. If $s <0$, then the $C^s_{sg}$ spaces are composed of distributions,  also defined by recursion on $[s]$  as follows:
$ f \in C^s _{sg} (\RR^d )$ if $f$ is a finite sum of partial derivatives (in the sense of distributions) of order 1 of elements of 
$C^{s+1}_{sg} (\RR^d )$. 
This allows to define the $C^s_{sg}$ spaces  for any $s \notin \ZZ$ (note that a consistent definition using the Zygmund classes can also be supplied for $s\in \ZZ$, see \cite{Mey90I},  however we will not need to consider  these specific values in the following).  
We can now define the parameter that we will use for determining boundedness. 

\BD The  {    uniform H\"older exponent}  of  a tempered distribution $f$ is 
\BE \label{caracbeswav2hol}  \Hmin = \sup \{ s : \; f \in C^{s }_{sg} (\RR^d) \}   . \EE
\ED
This definition does not make any a priori assumption: 
The uniform H\"older exponent  is defined for any tempered distribution, and it can be positive and negative. Note that we have already met this notion: Proposition \ref{propChents}  can be reinterpreted as stating that, if   $ 0 <  H-{1}/{\alpha} <1$, then $h^{min}_X \geq   H-{1}/{\alpha}$. In particular the remark following the proof of Proposition \ref{propirreg} allows to recover the fact that, for fBm, $\Hmin = H$; additionally, it yields that the same result holds for  the Rosenblatt process (see \cite{AbrPip,Pip} for the definition  and the properties of the wavelet expansion of this process).  

The  values taken by the uniform H\"older exponent  have the following interpretation:   
 \begin{itemize}
\item  If {  $\Hmin > 0$}, then $f $ is  a locally bounded function,
\item  if {$\Hmin <0$}, then $f$ is not a locally bounded function.
\end{itemize}
Typical examples are supplied by   L\'evy fractional stable processes which  already appeared in Section \ref{secpoin}; they    satisfy $\Hmin =  H- 1/ \al$, and therefore (as already mentioned) they have nowhere locally bounded stable sample paths if $H < 1/ \al$, and  continuous sample paths if  $H > 1/ \al$. 

The importance of this exponent lies in the fact that it does not only play the role of a prerequisite (i.e.,  a parameter whose value has to be positive if one wants to determine the oscillation scaling function), but also that it turns out to be a very efficient parameter for discrimination.
As will be illustrated later (cf. Section~\ref{sec-App}), the classification of several types of data can actually be based on the determination of their uniform H\"older exponent.  
\\

\noindent {\bf Uniform H\"older exponent and wavelets.} \quad In practice, $ \Hmin $ can be derived  directly from the wavelet coefficients of $f$ through a simple regression in a log-log plot; 
indeed, it follows from the wavelet characterization of the spaces $C^s_{sg}$  that:
\BE 
\label{caracbeswav3hol}  
\Hmin =\liminf_{j \rightarrow + \infty} \;\;  
 \;\;  \frac{ \log 
\left(  \displaystyle  \sup_{\la \in \La_j }  | \cla |   \right) }{\log (2^{-j})}. 
\EE
This estimation procedure has been studied in more detail in \cite{Bergou}. 

The exponent  $\Hmin$  can also be derived from the wavelet scaling function,  \cite{Bergou},
\[   \Hmin = \lim_{p \rightarrow + \infty}  \frac{\eta_f(p)}{p} = \lim_{p \rightarrow + \infty} {{\eta'}_f(p)}. \]     \\

\noindent {\bf Uniform H\"older exponent and applications.} \quad In Section~\ref{sec-App}, which  reviews a number of real-world applications,  we will see that 
 $ \Hmin$ can be found either positive or negative depending on the nature of the application.
For instance, velocity turbulence data (cf. \ref{sec-turb} and Fig.~\ref{fig-TurbFig}, middle row, left plot) and price time series in finance (cf. \ref{sec-finance} and Fig.~\ref{fig-finance}, 2nd and 3rd row row, second plots and bottom row, first plot) are found to always have $ \Hmin>0$, while aggregated count Internet traffic time series always have $ \Hmin<0$ (cf. \ref{sec-Internet} and Figs.~\ref{fig-TraffFigb} and \ref{fig-TraffFigc}, top row, left plots). 
For biomedical applications (cf. e.g., fetal hear rate variability as in Section~\ref{sec-fhbv} and Fig.~\ref{fig-fECG}, 3rd row), as well as for image processing (cf. Section~\ref{sec-VG} and Fig.~\ref{fig-VG}, 4th column), $ \Hmin $ can commonly be  found either positive or negative. 
For these last two examples,  $ \Hmin $ is found to be a  highly relevant parameter for classication purposes (cf. Fig.~\ref{fig-fECG}, bottom row, and Fig.~\ref{fig-VG}, bottom row, respectively). \\

\noindent {\bf Functions  that are not locally bounded and fractional integration.} \quad 
The examples listed above indicate that it is not uncommon in real-world application to find  that $ \Hmin <0$. 
However, this does not imply that further analyses which require that $\Hmin >0$ (and  will be described in the following) are impossible and should be discarded.
Indeed, \emph{fractional integration} offers a one to one way to associate  with any distribution $f$ (with potentially negative $ \Hmin$) a function $f^{(-s)}$, whose exponent  $ h^{min}_{f^{(-s)}}$ is positive, and therefore classification can then be operated on this new function $f^{(-s)}$ rather than on the initial data $f$.  
We  now expose this procedure, both from a theoretical and  a practical point of view. 
Note that there exist many variants of the definition of fractional integration, which would however all lead  to the same definition for the exponents used here, see \cite{PipTaq} for a discussion of fractional integration, especially in the context of stochastic processes. 

\BD 
Let  $f$ be a tempered distribution.
Fractional integral of order $s$ of $f$, denoted by  $f^{(-s)}$, is defined as convolution operator $(Id-\Delta )^{-s/2}$ which,  in the Fourier domain, is the multiplication by the function   $(1+|\xi |^2)^{-s/2}$.
\ED
Numerically, fractional integration is intricate to perform (because of finite size and border effect issues). 
However, for our purpose here,  these difficulties can be circumvented by instead using  \emph{pseudo-fractional integration}, defined as follows.
 
 \BD
 Let $s >0$, let $\pla$ be a  wavelet basis  whose elements belong to  $C^r$ with $r >s+1$;    let $f$ be a function, or a distribution, with wavelet coefficients $\cla$. The {  pseudo-fractional  integral} (in the basis $\pla$) of $f$ of order $s$, denoted by ${I}^s (f) $,  is the function whose wavelet coefficients (on the same wavelet basis) are 
\[ \csla  = 2^{-sj} \cla. \]
\ED 
This operation is numerically straightforward, and the following result shows that it retains the same pointwise and multifractal properties as exact fractional integration, see \cite{Bergou,Zuhai,Jaffard2006,AJL05}.
 
\BP \label{theopseudobis}  Let $f$ be a tempered distribution; then
\[  \forall s \in \RR, 
\left\{ \begin{array}{l}
 \forall p  >0, \hspace{3mm}  \eta_{ f^{(-s)} } (p)  =   \eta_{ {I}^s (f)  }  =  \eta_{ f } (p) +sp, \\Ê \\Ê
  \displaystyle h^{min}_{f^{-s}} = h^{min}_{I^s(f)} =  \Hmin  +s. \end{array} \right. \]
Additionally, the pointwise H\"older exponents satisfy
\[  \forall s > -\Hmin, \;\; \forall x \in \RR^d, \hspace{6mm} h_{f^{(-s)}} (x) = h_{ {I}^s (f)  } (x) .\]
\EP

This last property also holds for the leader  scaling functions which  will be introduced  in Section \ref{secleade} (under the same  condition $s > -\Hmin$, so that those quantities are well defined).
Therefore, in practice, for multifractal analysis, pseudo-fractional integration is preferred to exact fractional integration. 
These results point toward a  natural  strategy in order to perform the multifractal analysis of a function which is not locally bounded (or of a distribution): 
First determine  its exponent $\Hmin$; 
then, if  $\Hmin <0$, perform a pseudo-fractional integration  of order $s > - \Hmin$; 
it follows that the uniform regularity exponent of $I^s(f)$ is  positive, and therefore that it is a bounded function, to which multifractal analysis can be applied. 
Yet, this strategy rises an important problem: 
There is no canonical choice for the fractional integration order, and the shape of the  wavelet leader scaling function  (see Definition \ref{defilsf}  below) obtained after fractional integration may depend on this choice. 
In practice, when  multifractal properties of a collection of data have to be studied,  one follows the following strategy:  One first determines the exponent $\Hmin$ of each signal in the database; 
if some of them have negative $\Hmin $,  one picks  an exponent $s$  such that $ \Hmin +s $ is positive for all signals. Under these constraints,  $s$ should be picked as small as possible, so that data are modified by  a fractional integration of order as small as possible.   
The rule of thumb practically used is to select $s$ as the smallest multiple of $0.5$ ensuring that $ \Hmin +s $ is positive for all signals and to perform a pseudo-fractional integration { \em of the same order} $s$ for all  signals. 

\subsection{Wavelet leaders} 
\label{secleade}

\noindent {\bf From oscillations to wavelet leaders.} \quad  Our  purpose is now to obtain a wavelet reformulation of the oscillation scaling function (\ref{osf}) (in the same way as the 
wavelet scaling function is the wavelet  reformulation of the Kolmogorov scaling function).
This requires the introduction of  \emph{wavelet leaders}, which are local suprema of wavelet coefficients.
Let us start with a heuristic argument showing that they are natural candidates to estimate oscillation. 
Recall that, if $\lambda$ denotes a dyadic cube, then  $3 \lambda$ denotes the cube with same center and  three times wider. 

\BD 
Let   $f$ be  a    {locally  bounded function} satisfying (\ref{sg1}).   The   {   wavelet leaders  } of $f$ are the coefficients 
\BE  
\label{deflead} 
d_\lambda = \sup_{\lambda ' \subset 3 \lambda} | c_{ \lambda '} |.
\EE
\ED 
The fact that the supremum  in (\ref{deflead}) is finite is a consequence of the boundedness assumption for $f$.   
Therefore, checking that ${\Hmin}>0$ is  an absolute prerequisite prior to computing leaders. 

Let us now  sketch the reason why  wavelet leaders allow to estimate the oscillation. Consider a wavelet coefficient $ c_\la$. Since we use compactly supported wavelets, 
\[ c_\la = 2^{dj}\int _{ C \la } f(x) \psi^{(i)} (2^jx-k) dx \]
Since wavelets have vanishing integral,
\[\begin{array}{rl}  | c_\la |   &  = \left| 2^{dj}\displaystyle\int _{ C \la }  (f(x)- f(k 2^{-j} ) \psi^{(i)} (2^jx-k) dx \right| \\ & \\Ê
 &  \leq  2^{dj}\displaystyle\int _{ C \la }   \left| f(x)- f(k 2^{-j})  \right| |  \psi^{(i)} (2^jx-k) |  dx   \\ & \\ 
&  \leq  2^{dj} Os_f (C \la) \displaystyle\int _{ C \la }   \left|  \psi^{(i)} (2^jx-k)   dx  \right| \\ & \\Ê
& = C  \cdot Os_f (C \la) .  \end{array} \]
Consider now a given cube $K$;  
 this argument works for any wavelet whose support is included in $K$, so that 
each wavelet coefficient is bounded by  $C Os_f (C \la) \leq  C Os_f (K)$. This explains why  wavelet leaders  are bounded by the local oscillation.  

Converse estimates follow form the following argument:  From (\ref{ecrit22}), we get
\[  f(x) -f(y) =  \sum_{j,k} c^{(i)}_{j,k}\left(  \psi^{(i)} (2^jx-k)  -  \psi^{(i)} (2^jy-k)  \right)  \]
Let $j_0$ be such that $2^{-j_0} \sim | x-y|$. 
The low frequency terms ($j \leq j_0$) are estimated using the  smoothness of the wavelets, which makes the difference $ \psi^{(i)} (2^jx-k)  -  \psi^{(i)} (2^jy-k)$ small. 
The terms such that  $j_0 \leq j \leq A j_0$ (where $A$ has to be chosen correctly) are estimated using the bound of the wavelet coefficients supplied by the wavelet leader. 
Finally, high frequency terms ($j \geq  A j_0$) are estimated using  the uniform regularity of $f$.  This allows to obtain the required converse estimates. Note that subtle relationships between oscillations and wavelet coefficients have recently been obtained by M. Clausel and her collaborators, see e.g. \cite{ClauNic} and references therein.  \\ 
 
\noindent {\bf Leader scaling function.} \quad  The previous heuristic argument motivates the introduction of a new scaling function, constructed on the model of the wavelet scaling function, but replacing wavelet coefficients by wavelet leaders. 

\BD \label{defilsf}  Let   $f$ be  a    {locally  bounded function}  with slow growth (i.e. satisfying (\ref{sg1})), and let
\begin{equation}
\label{equ-WLSF}
 T_f (p, j) =  2^{-dj} \displaystyle\sum_{ \lambda \in \Lambda_j}  | d_\lambda |^p . 
 \end{equation}
The  {leader scaling function  } $\zeta_f (p)$  is given by
  \begin{equation} \label{vingub}  \forall p \in \RR, \hspace{6mm} 
  \zeta_f (p) =   \displaystyle\liminf_{j \rightarrow + \infty} \;\;   \frac{\log (T_f (p, j))}{\log (2^{-j})} \;  . \end{equation}
  \ED
The relationship between the leader scaling function and the oscillation scaling function is very similar to the one mentioned to exist between the wavelet scaling function and the  Kolmogorov scaling function: 
They coincide when $p \geq 1$;
therefore  the former  extends the later  for $p <1$, see \cite{Jaffard2004}. 

Let us discuss two consequences of this result:  
First, the fact that they  coincide  for $p =1$  implies that, if $\Hmin >0$,  then the upper box dimension of the graph of $f$  can be derived from  the leader scaling function, see \cite{jaffbox}: 
 \[  \overline{dim}_B (\gr (f)) =  \max ( d, d+1 - \zeta_f (1) ). \]
Second, when $\Hmin >0$, one can determine whether $f$  has bounded quadratic variation or not by inspecting its leader scaling function  for $p=2$, indeed:
\begin{itemize}
\item  If  { $\zeta_f (2) >d$}, then $f$ has bounded quadratic variation,
\item   if  { $\zeta_f (2) <d$}, then 
  the  quadratic variation of $f$ is unbounded.
\end{itemize}

\noindent {\bf Discussion of the condition $\Hmin >0$.} \quad
One may wonder if the condition $\Hmin >0$ is really necessary in order to obtain such  results. 
It is actually the case, and one can not obtain estimates of the quadratic variation (or of any $p$-variation) using wavelet methods for functions which are locally bounded but do not have some uniform smoothness for the following reason:
Recall that the Hilbert transform is defined as the convolution with $p.v. (1/x)$ (or alternatively as the Fourier multiplier by the function $sign( \xi )$) and consider as an example the sawtooth function $x-[x]$ which is bounded and has discontinuities at the integers.  
It obviously has bounded $p$-variation for any $p$. 
However, after applying the Hilbert transform to it, discontinuities are transformed into logarithmic singularities, and therefore the property of bounded $p$-variation is lost for any value of $p$. 
Since the Hilbert transform does not modify  the wavelet-based quantities  that we have considered, such as $\Hmin, \zeta_f (p)$, $\eta_f (p)$, it is clear that wavelet methods can not estimate $p$-variations of functions that have discontinuities, and therefore some uniform regularity assumption is indeed required.  

In finance, for instance, checking whether $\Hmin >0$ is positive or not is of double importance: 
First, it suggests to reject jump models for price; 
second, it enables to probe the finiteness of quadratic variations, a requested property to validate the use of stochastic integration on such data. 
This is further discussed in Section~\ref{sec-finance} and illustrated in Fig.~\ref{fig-finance}.  \\

\noindent {\bf The bonus of negative $p$s.} \quad  The leader scaling function can also be given a function space interpretation for $p >0$ in the framework of  \emph{oscillation spaces},  studied in \cite{Jaf05,JafOscill}. In particular, embeddings between these spaces and other, more classical, function spaces (such as Besov spaces) allow to derive an important relationship between the two scaling functions constructed through wavelet coefficients, see \cite{Jaffard2004}.

\BP 
  Let $f$ be a function satisfying  $\Hmin >0$. Then 
\[ \forall p >0,\hspace{1cm}  \eta_f (p) \geq  \zeta_f (p). \]
 Let  $p_c$ be the ``critical exponent'' defined by the condition $\eta_f (p_c) = d$, then 
\[ \forall p \geq p_c, \hspace{1cm}  \eta_f (p) = \zeta_f (p). \]
\EP

An important property of the leader scaling function is that it is \emph{well defined} also for $p <0$, although it can no longer  receive a function space interpretation.  
By well defined, it is meant that it has the following robustness properties if the wavelets belong to the Schwartz class (partial results still hold otherwise), see \cite{Zuhai,Jaffard2006}: 
\begin{itemize}
\item   $\zeta_f (p)  $ is  independent of the wavelet basis.
\item   $\zeta_f (p) $ is invariant under the addition of a $C^\infty$ perturbation.  
\item   $\zeta_f (p) $ is invariant under a $C^\infty$ change of variable.  
 \end{itemize}
The invariance under a smooth change of variable is a mandatory property for texture classification:  
Indeed, it is needed that natural textures can be recognized even after the deformation produced by setting them on smooth surfaces.  
 
The possibility of involving negative $p$s in analysis can turn crucial: A celebrated example is that of hydrodynamic turbulence where two multiplicative cascade models are classically in competition. 
Using positive $p$s only, the Kolomogorov and the  wavelet scaling functions computed from experimental data fail to discriminate between the two models. 
However, the leader scaling function, enabling the use of negative $p$s,  permits to show that one model is unambiguously preferred by data. 
This is discussed in Section~\ref{sec-turb} and illustrated in Fig.~\ref{fig-TurbFig} (bottom row plots).

\section{Beyond functional analysis: Multifractal analysis}
\label{secmulti}

The fact that the leader scaling function extends the analysis to negative $p$s shows that multifractal analysis allows to go beyond the scope of functional analysis.
Indeed,  the scaling function no longer has a  function space interpretation if $p <0$. 
Note that  this extension has some canonical  character, as a consequence of its independence of the chosen wavelet basis, and also because of its robustness properties. 
The point made in the present  section is to show that even more  is true by exhibiting a natural interpretation of  scaling functions, which requires both positive and negative values of $p$. 
This interpretation is in terms of the pointwise singularities of the function, and will be supplied by the multifractal formalism, which we now describe.

\subsection{Multifractal formalism for oscillation} 
\label{secrelat1}

\noindent {\bf Multifractal spectrum versus multifractal formalism.} \quad
From now on, it is assumed that $\Hmin >0$, 
which implies that $f$ is a continuous function, and therefore that its oscillations are well defined. 
Recall that the notion of oscillation came up in two occurrences until now: 
The definition of the oscillation scaling function (\ref{osf}), and the characterization of the pointwise H\"older exponent (\ref{holoscill}). 
This points towards a possible relationship between both quantities, which can be made explicit via \emph{multifractal formalisms}:
Specifically, a multifractal formalism establishes a connection between the scaling functions and the spectrum of singularities (or multifractal spectrum) of $f$, which is defined as follows.

\BD 
Let  $f$ be a   locally  bounded function. 
The  { spectrum of singularities} of $f$ is the function  
\[  D_f (h)= \dim (E_f(h) ), \] where we recall that $\dim$ denotes the Hausdorff dimension, and the $E_f (h) $ are the isoh\"older sets 
$  E_f(h) = \{ x_0: \hspace{3mm}  h_f(x_0) = h\} . $
\ED

This definition justifies the denomination of { \em multifractal functions}:   On the theoretical  side,  one typically considers  functions $f$ that have non-empty isoh\"older sets for $h$ taking all values in an interval $[ \Hmin , h^{max}_f ]$, and therefore one has to  deal with potentially  an infinite (and even uncountable) number of fractal sets $E_f (h) $ whose Hausdorff dimensions have to be determined.

Recall that $\dim (\emptyset) = -\infty$ so that  $ D_f (h) =-\infty $  if $h$ does not belong to the range of $h_f$.   We will call { \em support} of $D_f$  the set 
\[ supp (D_f) = \{ h : \;\; D_f(h) \geq 0 \} = \{ h: E_f (h) \neq \emptyset \} . \]
(Note  that the support of the spectrum needs not be a closed set, in contradistinction with the usual definition of the support of a function). 

The initial formulation of the multifractal  formalism was due to the physicists G. Parisi and U. Frisch and was  based on Kolmogorov scaling function  \cite{ParFri85}. 
It was further  grounded mathematically in \cite{Bmp92}, see also \cite{HJKPS}.
Here, the focus is on a formalism based on the oscillation scaling function, which we describe first since it paves the way towards a more satisfactory leader scaling function formulation. \\

\noindent {\bf Oscillation multifractal formalism.} \quad In this paragraph, it is assumed that $h_f (x) $ takes on values between 0 and 1, so that (\ref{holoscill}) is valid. 
Since  pointwise H\"older  exponents can thus be derived from the quantities $Os_f (3 \la) $, it is natural to consider scaling functions based on these quantities, i.e., the oscillation scaling function  of $f$  defined  by (\ref{osf}).

Let us now show how the spectrum of singularities can be obtained from the scaling function.
The definition of the scaling function roughly means that the quantity $R_f (p, j)$ defined by (\ref{posc}) satisfies  $R_f  (p, j) \sim 2^{- Os_f (p)j}$. 
Let us  estimate the contribution to $R_f  (p, j)$ of the cubes $ \la$ that cover the points of $\EMH$;
(\ref{holoscill}) asserts that they satisfy
$Os_f (3 \la)  \sim 2^{-hj}$; since we need about $2^{-D_f (h) j}$ such cubes to cover $\EMH$, the corresponding
contribution roughly  is 
\[ 2^{-dj}2^{D_f (h) j}2^{-hpj} = 2^{-(d-D_f (h) +hp)j}.\]
 When $j \rightarrow +
\infty$, the dominant contribution stems from the smallest exponent  (the others bring contributions which are exponentially smaller), so that 
\begin{equation} 
\label{etamuleg} 
Os_f (p) =  \inf_h (d- D_f (h)+ hp).
\end{equation}
It follows from the definition of  the scaling function $Os_f (p)$  that it is  a concave function on $\RR$. 
This is in agreement with the fact that the right-hand side of (\ref{etamuleg}) necessarily is a
concave function (as an infimum of a family of linear functions) no matter whether $D_f (h)$ is concave or not.
However, if the spectrum also is   a concave function, then the Legendre transform in (\ref{etamuleg})  can be inverted (as a consequence of a general result on the duality of convex functions), which justifies the following definition.

 \BD 
 Let $f$ be a  function which satisfies $\Hmin >0$; $f$ follows the multifractal formalism if its spectrum of singularities satisfies 
\begin{equation} 
\label{formult1}
D_f (h) = \inf_{p } (d-Os_f (p) + hp).
\end{equation} 
\ED

\noindent {\bf Validity of the multifractal formalism.} \quad The derivation exposed above is not a strict mathematical proof, but rather consists of a heuristic explanation of the relation between $D_f(h)$ and $Os_f (p)$.
The determination of the range of validity of (\ref{formult1}), or of any of its variants, is one of the main open mathematical problems in multifractal analysis.   
Nonetheless, let us emphasize that the cornerstone of this derivation relies heavily on (\ref{holoscill}), i.e., on the fact that the H\"older exponent of a function can be estimated from the set of values that its oscillation takes on the ``enlarged'' dyadic cubes $3 \la$.
In particular, the variant of this multifractal formalism based on increments (i.e., on the Kolmogorov scaling function) has a narrower range of validity, see \cite{Ciocco1,Zuhai,Jaffard2006} for a discussion. 

An important remark is that, in applications, one should not focus too much on the problem of the validity of the multifractal formalism. 
On one side, because one can not compute directly a spectrum of singularities of experimental data, the multifractal formalism can be checked only for mathematically defined functions or processes. 
On other side, the  purpose of multifractal analysis often is to use a spectrum as a tool for classification or model selection, and, with this respect, using a scaling function (or it Legendre transform) is  a perfectly valid approach, independently of the fact that the multifractal formalism holds or not. 
However, Section  \ref{sec-mfMF} below details an example where the multifractal formalism allows to determine effectively the multifractal spectrum: It can be the case when it is degenerate, i.e., when $h_f(x)$ is a constant function (in which case $D_f$ is supported by  only one point). \\

\noindent {\bf Further comments.} \quad  The formulation of the multifractal formalism given by (\ref{formult1}) has the following  advantage: 
The scaling function  is invariant under  bi-Lipschitz deformations of the function, which is a natural requirement since the spectrum of singularities possesses this invariance property (as  a consequence of the assumption that the range of $h_f$ is included in $(0,1)$). 

To this Legendre transform based multifractal formalism, initially proposed in \cite{ParFri85} and thoroughly theoretically grounded in \cite{Jaffard2004,Riedi2003}, has been opposed an alternative formulation relying on the \emph{large deviation principle} (cf. e.g., \cite{Riedi2003} for an early introduction and \cite{BarGon} for a recent version involving oscillations of order 1). 
However, we will not follow this track for the following reasons:  
\begin{itemize}
\item  Using oscillations restricts the possible range of H\"older exponents to $h \in [0,1]$;
\item the  method can not be adapted to data that are not locally bounded (which is often the case in signal and image processing).
\end{itemize}
This second drawback, as we saw, can be circumvented easily by a pseudo-fractional integration when using wavelet techniques, and this is the reason why we now turn towards these methods.  

The limitation to the range $h \in [0,1]$ implied by either the use of increments of order 1 (Kolmogorov scaling function) or oscillations of order 1 (oscillation scaling function) has long been recognized and led to the pioneering contributions of A. Arneodo and collaborators promoting the use of wavelet based methods in multifractal analysis (cf. e.g., \cite{AADMV,ABM95,muzyetal91,muzyetal94} and references therein). 

\subsection{A fractal interpretation of the leader scaling function: The leader multifractal formalism}
\label{secspec}
 
The introduction of wavelet leaders was motivated by the fact that they allow to estimate local oscillations, and therefore come as a natural ingredient in a wavelet reformulation of the oscillation scaling function and of pointwise H\"older regularity. 

Let $x_0 \in \RR^d$, and recall that  $\lambda_j (x_0)$  denotes the dyadic cube of width  $2^{-j} $  which contains $x_0$. If $\Hmin >0$,  then the heuristic relationship that we gave in Section \ref{secleade} and which related  oscillations and wavelet leaders actually leads to the following result (see \cite{Jaffard2004} and references therein):
\begin{equation} \label{repun}   h_f(x_0) = \liminf_{j \rightarrow + \infty}\;  \; \frac{ { \log 
 \left( d_{\lambda_j (x_0)} \right)  }}{\log (2^{-j})}. \end{equation}
The heuristic arguments  that were developed above for the derivation  of the oscillation multifractal formalism can be reproduced in the setting of wavelet leaders.
They lead to the following reformulation of the multifractal formalism. 

\BD  Let $f$ be a function statisfying $\Hmin >0$; 
the { leader spectrum}  of $f$ is defined through a Legendre  transform of  the leader scaling function, i.e. 
\[ 
L_f(h)  =    \inf_{p \in \RR}  \left( d+hp - \zeta_f  (p)  \right).      
\]
The { wavelet   leader multifractal formalism} holds if 
 \BE 
\label{wlmf}    \forall h \in \RR ,  \hspace{3mm}   D_f (h) =  L_f(h).  \EE
 \ED

This formalism holds  for large classes of  models used in signal and image processing, such as fBm,  lacunary and random wavelet series \cite{AJ02}, or multiplicative cascades.  
In this wavelet leader setting, the justification of the derivation  of (\ref{wlmf})  relies heavily on (\ref{repun}).
In particular, multifractal formalisms based on other quantities (such as the wavelet  scaling function) have a narrower range of validity, see \cite{Jaffard2004} for a discussion. 

The multifractal formalism does not hold in all generality, nonetheless the following upper bound does, yielding a general relation between the leader scaling function and the spectrum of singularities.

\begin{Theo}  
\label{theomaj} 
If $\Hmin >0$,    then, 
\[  \forall h \in \RR ,  \hspace{3mm}   d_f (h)  \leq  L_f(h).  \] 
\end{Theo}

\begin{figure}[h]
\centering
\includegraphics[width=1\linewidth]{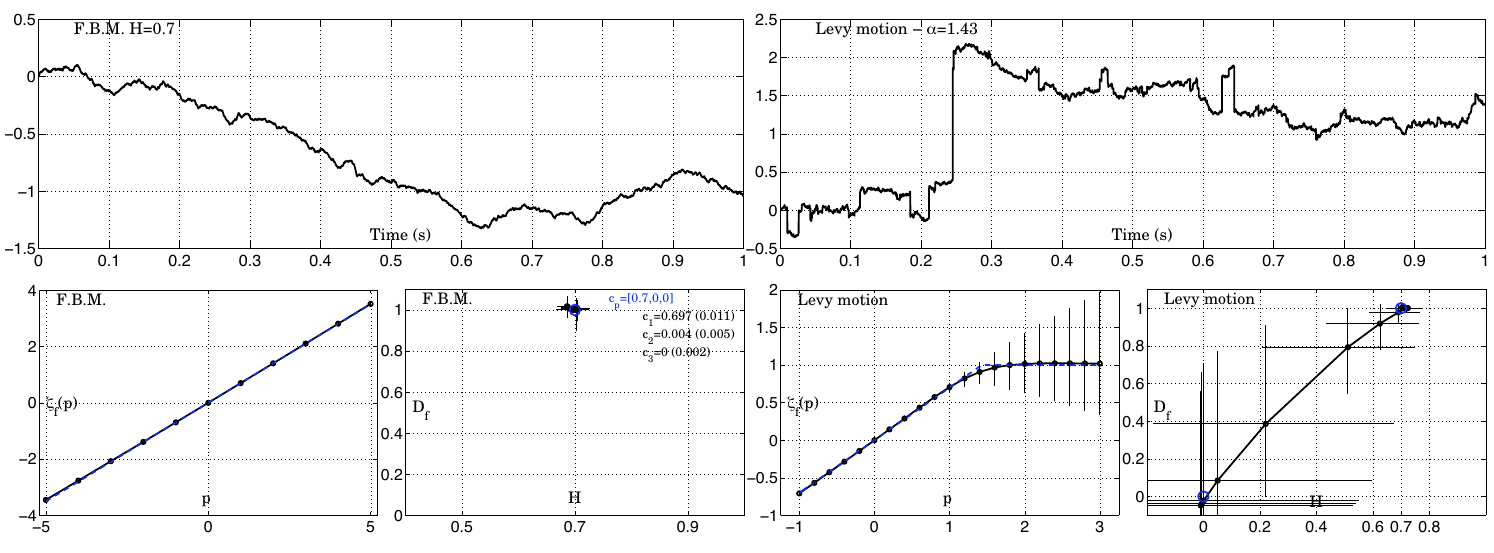}
\caption{\label{fig-wlmfSS} {\bf Selfsimilar processes.} Top row: sample paths of fBm with $H = 0.7$  (left) and of  stable L\'evy motion with $H = 0.7 $ (right). 
Bottom row: estimated wavelet leader scaling functions and Legendre spectra.  
Both processes have leader scaling functions that behave linearly in $p$ for $p$ close to $0$: $\zeta_f(p) = pH$ for fBm,  $\zeta_f(p) = H  p$ (if $p$ is close to 0) for  stable L\'evy motion. These self similar processes with stationary increments  are thus characterized by $c_2 = 0$. 
While fBm shows a linear in        $p$ leader scaling function and is mono-H\"older,   stable L\'evy motion is multifractal and its leader scaling function is only piecewise linear (cf. discussion in Section~\ref{sec-mfMF}). Confidence intervals are constructed based on the non parametric bootstrap procedure described in \protect \cite{WENDT:2007:E}.}
\end{figure}

\begin{figure}[h]
\centering
\includegraphics[width=1\linewidth]{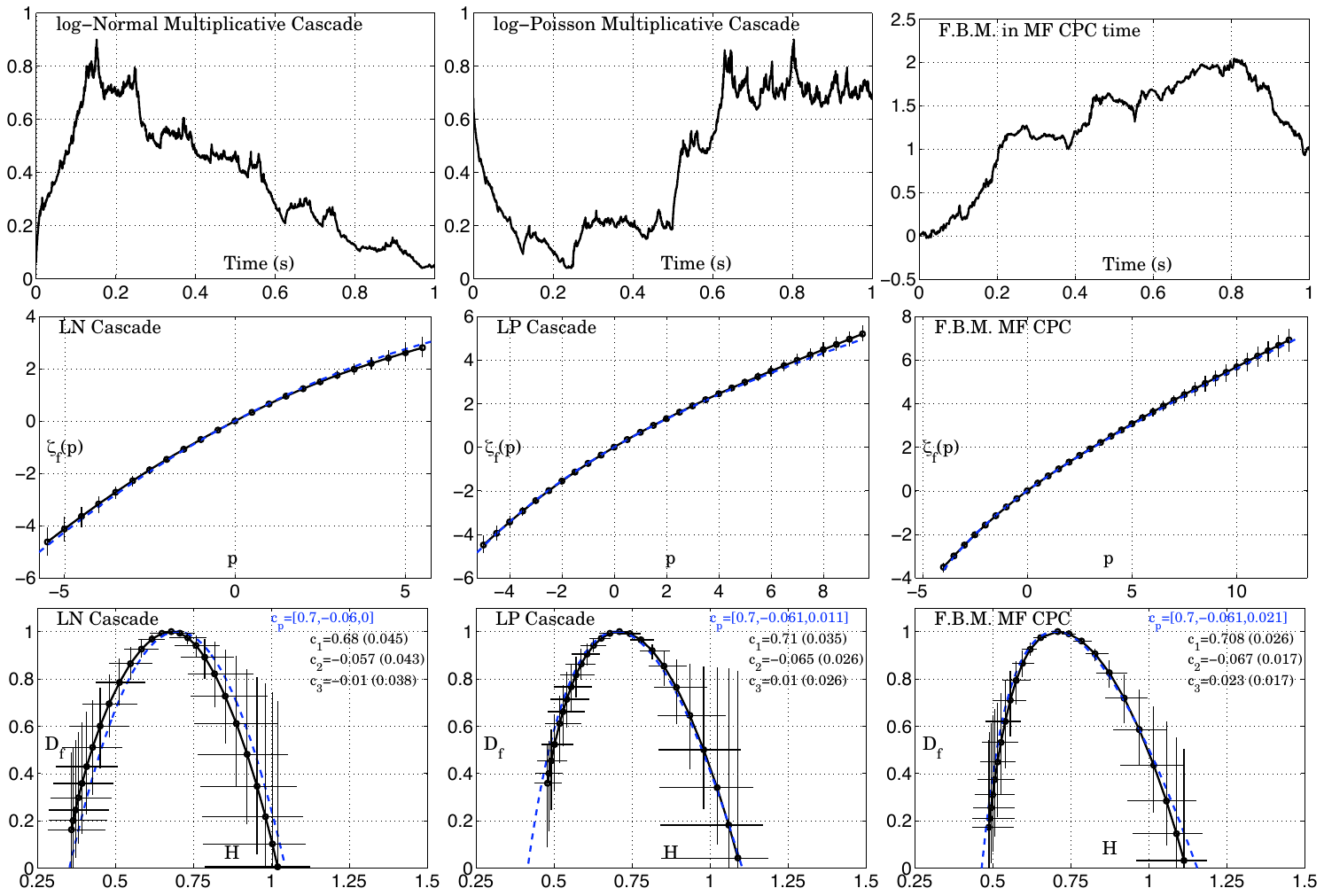}
\caption{\label{fig-wlmf1D} {\bf 1D multifractal processes.} Top row: sample paths of log-normal cascade (left), log-Poisson cascade (middle), and a fBm in multifractal time (as defined in (\ref{equ-fbmmf})) (right).  
2nd and 3rd rows show the superimposition of the theoretical (blue) and estimated (black) corresponding leader scaling functions and multifractal spectra. Theory and estimation are strikingly in agreement both for positive and negative $p$s and across the entire multifractal spectra (increasing and decreasing parts). Confidence intervals are obtained by a non parametric bootstrap procedure described in \protect \cite{WENDT:2007:E}.}
\end{figure}

\begin{figure}[h]
\centering
\includegraphics[width=1\linewidth]{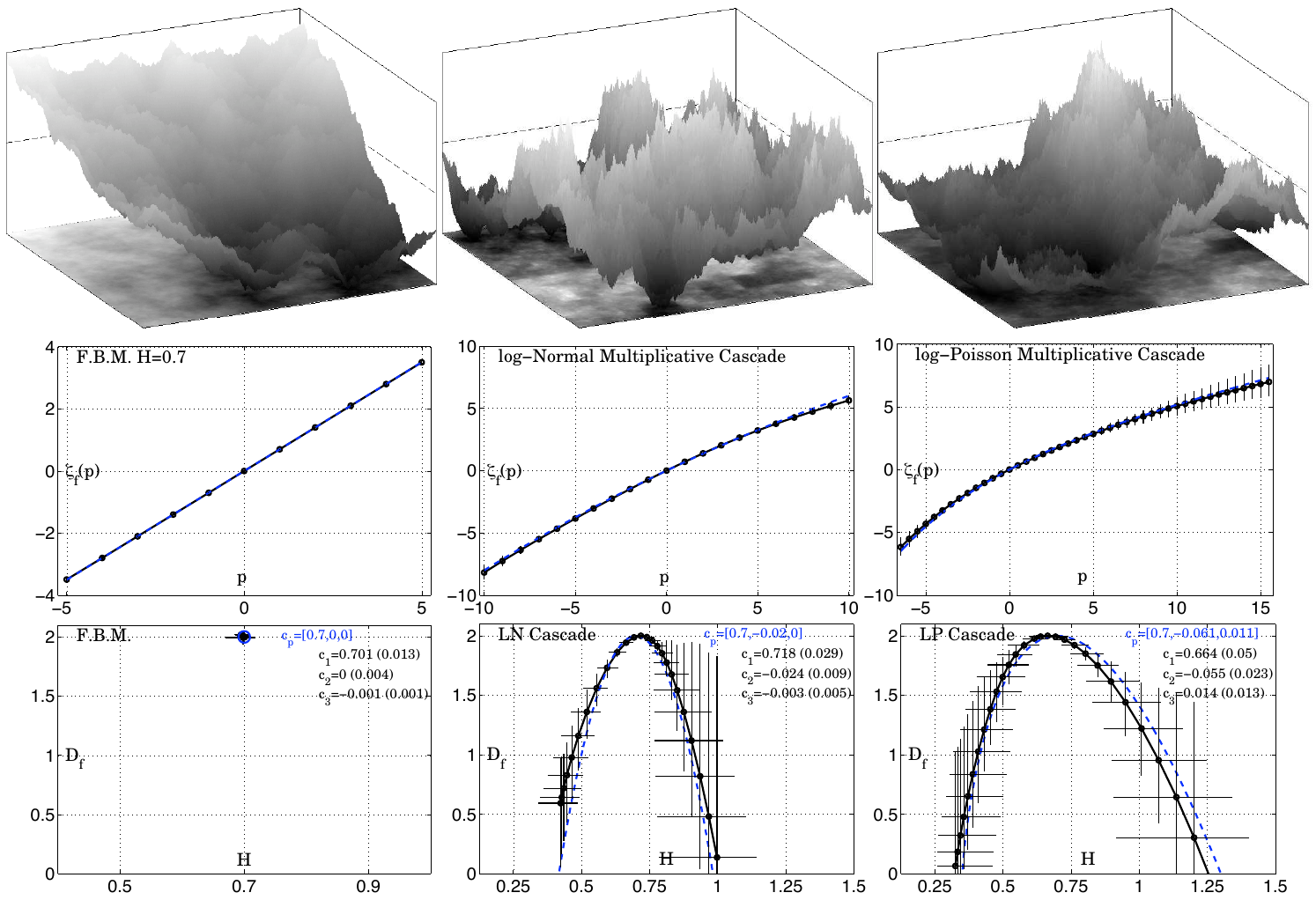}
\caption{\label{fig-wlmf2D} {\bf 2D multifractal processes.} Top row: 2D realizations of fBm (left), log-normal cascade (middle), log-Poisson cascade (right).  
2nd and 3rd rows show the superimposition of the theoretical (blue) and estimated (black) corresponding leader scaling functions and multifractal spectra. As in the 1D case (cf. Figs. \ref{fig-wlmfSS} and \ref{fig-wlmf1D}), theory and estimates are strikingly in agreement both for positive and negative $p$s and across the entire multifractal spectra (increasing and decreasing parts). Confidence interval are constructed using a non parametric bootstrap procedure described in \protect \cite{WENDT:2007:E}.}
\end{figure}

Examples of leader multifractal spectra $L_f(h)$ computed from a single realization of various processes and fields are illustrated in Figs.~\ref{fig-wlmfSS}, \ref{fig-wlmf1D} and \ref{fig-wlmf2D}, and compared with the theoretical multifractal spectra $D_f(h)$. 
It can be observed that, for all these reference processes, the leader multifractal spectra computed from a single realization very satisfactorily match the theoretical multifractal spectra $D_f(h)$, hence suggesting that the multifractal formalism,  $L_f(h)=D_f(h)$, holds in all cases. 

The following relationship holds for most  standard models used in multifractal analysis: 
\BE \label{egalhmin}   \Hmin = \inf   ( supp (D_f) ) . \EE
Note however that this relationship is not true in all generality; for instance the ``chirp'' 
\[ F_{ \alpha, \beta} (x) = | x|^\al \sin \left ( \frac{1}{| x|^\beta} \right)  \]
(for $\al$ and $\beta$ positive) satisfies 
\[ h^{min}_{F_{ \alpha, \beta}} = \frac{\al }{ \beta +1} \hspace{3mm}  \mbox{and } \hspace{3mm} \inf \left( supp \left(D_{F_{ \alpha, \beta}}\right) \right)   = \al.  \] 
However, we will mention a partial result confirming (\ref{egalhmin}) in Section \ref{homog}.

\subsection{Local spectrum and homogeneity}
\label{homog}

So far, little attention has been paid to the region where multifractal analysis is performed, implicitly assuming that it is conducted on the entire set of data available.  
However, one may wonder whether multifractal analysis would yield different results when applied to different domains.  
The answer is often negative:  Scaling functions, and  spectra of singularities, when  measured on a sub-part of the data (say a non-empty open subset, in order to dismiss boundary problems) do not depend on the particular open set which is chosen. 
When such is the case, the function analyzed  is said to have \emph{homogeneous multifractal properties}.  
It is observed that large categories of  experimental data  have homogeneous multifractal properties: 
It is for instance the case for fully developed turbulence, where  analyses performed on different parts of the signal always yield the same results. However,
this might  not systematically be the case. 

A simple cause of non-homogeneity can be that data consist of a juxtaposition of different  pieces, concatenated together. 
This situation is commonly met in image processing:
Indeed, because of occlusion effects (i.e. objects in the front of the picture partly mask objects behind), images usually appear as a patchwork of  homogeneous textures. 
In this case, it is reasonable to isolate the different components of the data, and perform multifractal analysis over each piece independently. 
Scaling functions computed on the whole image simply result as the minimum of the \emph{local} scaling functions, and conversely the global multifractal spectrum consists of the maximum of the local spectra. 
Note that such cases can lead to simple examples of (global) non-validity of the multifractal formalism.
Indeed, the scaling function of the whole image remains, by construction, a concave function, whereas the spectrum, as a maximum of concave function needs no more be concave, even if each component is.
Therefore, the multifractal formalism will fail, since, by construction, the right hand side of (\ref{formult1}) always yields concave functions. 

However, more intricate situations may exist, where multifractal characteristics can evolve with time (or the location, for images) in a way which is more involved than piecewise constant.  
One can expect, for instance, that stock market data display different statistics during different market phases, such as crises or bubbles. Different multifractal characteristics could stem
from such differences. 
One may also expect that the whole conditions of economy evolve with time, as a consequence of  the ever increasing complexity of financial operations,  of new financial products, and (more subtly) of mathematical models which themselves influence the way that operators react to the market; all this pleads for possible smooth (or not so smooth) modifications of multifractal characteristics with time. 
Such possible mechanisms may help to interpret the empirical observations that  multifractal properties of market price fluctuations,  such as those reported in Fig.~\ref{fig-finance} (bottom row plots) for the hourly Euro-USD rate, are significantly changing along time (cf. Section~\ref{sec-finance} for discussion).  
 
On the mathematical side, such situations have been barely studied so far (see a contrario the results concerning heterogeneous ubiquity, by J. Barral, A. Durand and S. Seuret,  where sophisticated mathematical tools have been developed to deal with the theoretical challenges of this situation \cite{BaSe2,BaSe3,BaSe1,Dur3}).  
A typical  example, studied in \cite{BaJaFoSe}, is supplied by Markov processes (i.e. processes with independent increment)  which do not have  stationary increments (and thus are not L\'evy processes).  
Roughly speaking, such processes have jumps of random amplitude, but each jump brings the process in a new environment,  and the local multifractal properties will depend on this environment  which is randomly chosen. One therefore obtains a non-homogeneous process, with a random spectrum, and random scaling functions. 
To the opposite, homogeneous multifractals have some additional regularity properties. For instance, the pitfall mentioned in the previous subsection does not occur: An homogeneous multifractal function satisfies (\ref{egalhmin}). 
An interesting phenomenon of homogeneity breaking will be discussed in the next subsection: 
We will consider a  homogeneous monoh\"older stochastic process, whose square is not homogeneous, and no longer monoh\"older. 

Note also that  a test for the constancy along time of the multifractal properties measured using the leader scaling function and relying on a non parametric bootstrap procedure was proposed in \cite{ICASSP}. 

\subsection{Mono vs. Multifractality?}
\label{sec-mfMF}

\noindent {\bf Disentangling issues.} \quad For practical purposes, it is often of interest to decide whether data are mono- or multifractal.
However, despite apparent simplicity, as such the question is ill-posed and can be rephrased in different ways: 

\begin{enumerate}
\item  Are data characterized by a single H\"older exponent that takes the same and unique value everywhere or do there exist a variety of H\"older exponents in data? 
\item Are data better described by a selfsimilar process (often implicitly meaning fBm) or by a cascade based process? 
In this later formulation, the important underlying question is: Is there an additive (as in random walk and hence selfsimilar models) or a multiplicative (as in cascades) structure underlying data and hence revealing significant differences in the nature of the data.
\item Are the scaling exponents (i.e., the scaling function) measured on data linear in $p$ or not?
\end{enumerate}
 Relying on the fBm (Gaussian based) intuition, these three different formulations are often thought to be equivalent. 
This section aims at discussing some of the subtleties underlying the relations between these three questions. \\

\noindent {\bf Log-cumulants.} \quad  Recall that, if they exist, the cumulants of a random variable $X$ are the coefficients of  the Taylor expansion of the second characteristic  function of $X$, i.e. the coefficients  $c_m$ (if they exist)  in the expansions 
\BE \label{defcum}  \log \left( \E \left( e^{pX} \right) \right) = \sum_{m=1}^\infty c_m \frac{p^m}{n !} . \EE 
Following \cite{temperature,Castaing1993,dma01}, the cumulants of the log of the wavelet leaders are now introduced. 
Let $C_m(j) $ denote the $m$-th order cumulant of the random variables $\log (d_{j,k})$. Under  stationarity and ergodicity assumptions on this sequence,  $C_1(j) $ is the expectation of the $\log (d_{j,k})$, $C_2(j)$ their  variance,  etc. 
Assuming that the moments of order $p$ of the leader $ d_{j,k}$ exist, starting from the assumption $\E (d_{j,k}^p)  =  \E (d_{0,0}^p ) \cdot   2^{-j \zeta_f(p)}$, one obtains  that
\[ \log (\E (d_{j,k}^p) ) =  \log (\E (d_{0,0}^p) ) +  \zeta_f(p)  \log (2^{-j}) ;\]
using (\ref{defcum}), 
one obtains the following expansion for $p$ close to $0$:
 \[ \log (\E (d_{j,k}^p) ) =  \log (\E (e^{p \log d_{j,k}})) =  \sum_{m\geq 1} C_m(j) \frac{p^m}{ m!} . 
\]
Comparing these two expansions, we obtain that the $ C_m(j) $ are necessarily of the form
\begin{equation}
\label{equ-leadercum}
 C_m(j) = C^0_m + c_m \log (2^{-j} ) ,
\end{equation}
and that $\zeta_f(p)$ can be expanded around $0$ as
\begin{equation}
 \zeta_f(p)  =\sum_{m\geq 1} c_m \frac{p^m}{ m!}. 
\end{equation}
The concavity of $ \zeta_f$ implies that  $c_2 \leq 0$. 
Note that, even if the moment of order $p$ of $d_{j,k} $ is not finite, the cumulant of order $m$ of $\log (d_{j,k}) $ is likely to be finite and (\ref{equ-leadercum}) above is also likely to hold. 
Moreover, (\ref{equ-leadercum}) provides practitioners with a direct way to estimate the $c_m$ by means of linear regressions in $ C_m(j)$ versus $\log  (2^{-j} ) $ diagrams \cite{WENDT:2007:E}. 
The $c_m$ are hence not derived from an a posteriori expansion performed from the estimates of the $\zeta_f(p)$, but instead estimated directly.
For the sake of completeness, let us mention that  the polynomial expansion of $ \zeta_f(p)$ around $0$ can be translated, via the Legendre transform,  into an expansion of $L_f(h)$ around its maximum \cite{Wendt2008c}, on condition that $c_2 \neq 0$:
\begin{equation}
 L_f(h) = d + \frac{c_2}{2!}\left( \frac{ h-{c_1}}{c_2}\right)^2+\frac{-{c_3}}{3!}\left( \frac{ h-{c_1}}{{c_2}}\right)^3+\frac{-{c_4}+3{c_3}^2/{c_2}}{4!}\left( \frac{ h-{c_1}}{{c_2}}\right)^4+ \ldots
\end{equation} 

\noindent {\bf Mono vs. Multi-H\"older.} \quad The only general mathematical result concerning the  validity of the multifractal formalism is supplied by Theorem \ref{theomaj}, which, in general does not yield the spectrum. However, it does yield an interval in which the support of the spectrum is included: Indeed, let 
\[ \Hmin = \inf \{ h:  L_f(h)   \geq 0\}  \hspace{5mm} \mbox{ and}  \hspace{5mm} h_f^{max} = \sup \{ h:  L_f(h)   \geq 0\} . \]
Then, the support of the spectrum $D_f (h)$ is included in the interval $ [ \Hmin  , h_f^{max} ] $.  (It can indeed be checked that $\Hmin  $ thus defined  actually coincides with the uniform H\"older exponent, defined in Section \ref{secunif}.)  
An important particular case stems from the situation where $  \Hmin  $ and $ h_f^{max}$ coincide, in which case the spectrum is supported by the sole point $\Hmin  = h_f^{max}$, so that only one single  H\"older exponent exists in the data, which therefore constitute a monoh\"older function satisfying 
\[  \forall x, \hspace{6mm}  h_f (x) = \Hmin  =  h_f^{max}.  \]
This puts into light the necessity to use both positive and negative values of $p$ in the estimation: If the leader scaling function was determined only for positive $p$,  the upper bound supplied by Theorem \ref{theomaj} would yield an increasing  function, and therefore would, at best,  yield the increasing part of the spectrum.  
In particular, $h_f^{max}$ would not be obtained, hence preventing from drawing any conclusion with respect to data being monoh\"older or not. 

However, testing the theoretical equality $\Hmin  =  h_f^{max}$ is practically difficult as the estimation of these quantities turns out to be poor (cf. e.g., \cite{WENDT:2007:E,WENDT:2009:C}). 
Instead, one can estimate the parameter $c_2$ as defined above in (\ref{equ-leadercum}). 
When the estimated $c_2$ is found strictly negative, this unambiguously indicates multi-H\"older function. 
As discussed below, the converse is however not necessarily true: $c_2 = 0$ does not necessarily imply monoh\"older function, as may be incorrectly extrapolated from the Gaussian fBm case.

Let us also note that both approaches based on multifractal analysis can be regarded as non parametric: Testing for mono- or multi-H\"older can be achieved without recourse to any parametric setting,  in sharp contradistinction with other existing estimators (for example, many estimators evaluate the selfsimilarity index $H$ assuming that the data  are sample paths of an fBm).   

Estimated $c_2$ measured on real-world data are reported in Section~\ref{sec-App}:
Turbulence velocity is well-known as the emblematic real-world example where a  non vanishing parameter  $c_2$ is met, (cf. Section~\ref{sec-turb}, Fig.~\ref{fig-TurbFig}).
This is observed to be also the case for the Euro-USD price fluctuations (cf. Section~\ref{sec-finance} and Fig.~\ref{fig-finance}, bottom row, left plot) 
and for fetal heart rate variability (cf. Section~\ref{sec-fhbv} Fig.~\ref{fig-fECG}).
However, Internet traffic aggregated count time series are found to be characterized by $c_2=0$ at coarse scales and weakly departing from $0$ at fine scales (cf. Section~\ref{sec-Internet} and Figs.~\ref{fig-TraffFigb} and \ref{fig-TraffFigc}, bottom row, left plots). \\

\noindent {\bf Selfsimilarity?} \quad  Selfsimilar processes are often obtained as renormalized limits of random walks (see for instance the fBm, or stable L\' evy processes).  Because of the additive structure underlying their  nature, selfsimilar models are appealing for describing data.
From the examples shown in Section~\ref{secspec}, it can even be conjectured that scaling functions often are piecewise linear functions. 
In particular, combining  ergodicity and stationarity properties,  the Kolmogorov  scaling function  of selfsimilar processes with stationary increments should  behave linearly in $p$,  as $\eta_f(p) = pH $ on an interval of $p$ containing $0$  (where $H$ is the selfsimilarity parameter). Conjecturing that  this property also holds for   $\zeta_f (p)$, it follows  that 
 testing for selfsimilarity can be formulated as testing for $c_2 = 0$.
This has been conducted in a non parametric bootstrap framework, cf. e.g., \cite{ICASSP,WENDT:2007:E}. 
Note, however, that selfsimilarity with stationary increments together with $c_2=0$ does not imply mono-H\"older sample paths (cf. the example of Levy motion reported in  Section~\ref{secspec} above). 
Furthermore, $c_2 = 0$ indicates selfsimilarity but does not necessarily imply that data follow a fBm model.
This would in addition require to test for joint Gaussianity of all finite distribution functions. \\

\begin{figure}[h]
\centering
\includegraphics[width=1\linewidth]{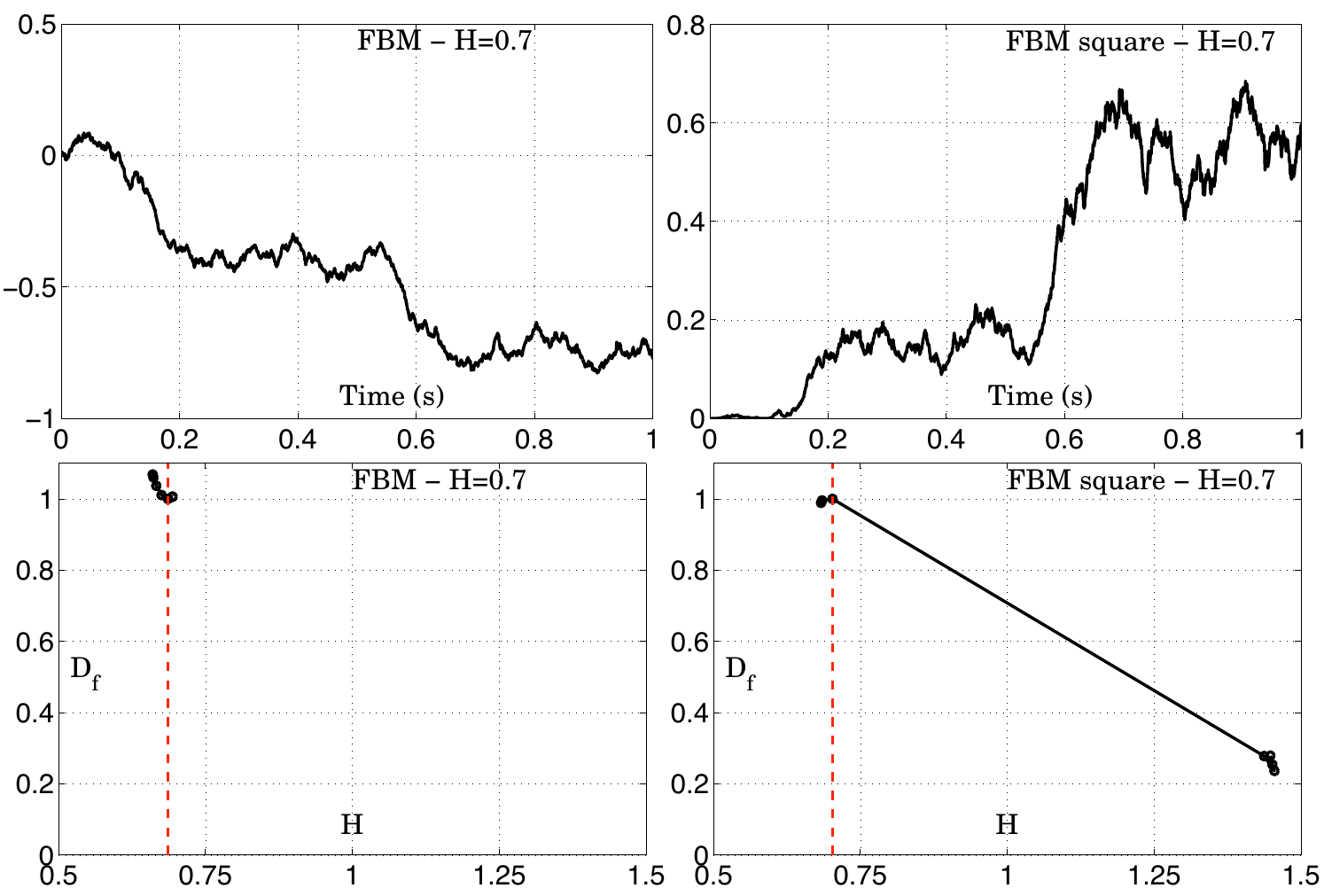}
\caption{\label{fig-FBMsquareFig} {\bf Non linear transform.} Top row: sample path of fBm, with $H = 0.7$  (left) and its squared transform (right). 
Bottom row: estimated wavelet leader based Legendre spectra.  
fBm is  a mono-H\"older process and its estimated Legendre spectrum collapses on the theoretical spectrum $D (h) =1$ for $h= H$. Its squared non linear point transform, however, is a bi-H\"older process, with theoretical multifractal spectrum given in \protect (\ref{equ-fbmsquared}): The estimated Legendre spectrum perfectly reproduces its corresponding concave hull.}
\end{figure}

\noindent {\bf Non linear point transform.} \quad Let us give some further insights into the relevance of the question of deciding whether data are mono or multi-H\"older. 
We aim at showing that the notion of monoH\"older function is, in some sense, unstable  since it can be altered under the action of the most simple smooth non-linear operator. 
To that end, let us consider  fBm $B_H(t)$, which is the simplest example of a mono-H\"older stochastic process: The H\"older exponent is everywhere equal to $h=H$. 
Let us now further consider its square transform $Y_H(t )= (B_H (t))^2$. 
On one hand, at points where the sample path of  fBm does not vanish, the action of the mapping $x\rightarrow x^2$ locally acts as a $C^\infty$ diffeomorphism, and the pointwise regularity is therefore preserved. 
On  other hand,  consider now  the (random) set $A$ of points where fBm vanishes.
 The  uniform modulus of continuity of fBm implies that 
 a.s., for $s$ small enough, 
 \[ \sup_t  | B_H(t+s ) -B_H(t) | \leq C |s|^{H}\sqrt{  \log (1/| s|) } . \]
 therefore, if $B$ vanishes at $t$, then $ B_H(t+s  )^2  \leq C |s|^{2H} \log (1/| s|) $,  so that $h_Y (t) \geq 2H$. The converse estimate follows from the fact that,
 for every $t$, 
 \[ \limsup_{ s \rightarrow 0} \frac{ | B_H(t+s ) -B_H(t) |}{  |s|^{H}  }  \geq 1,  \]
 so that, if $B_H(t) =0$, then 
  \[ \limsup_{ s \rightarrow 0} \frac{ ( B_H(t+s ) )^2}{  |s|^{2H}  }  \geq 1,  \]
  so that $h_Y (t) \leq 2H$.
Thus, at vanishing points of $B_H$,  the  action of the square is to shift the H\"older exponent from $h=H$ to $h=2H$.  
This set of points has been the subject of investigations by probabilists (cf. the flourishing literature on the {\em local time} either of Brownian motion, or of fBm);
in particular, it is known to be a fractal set of dimension $1-H$, cf \cite{MonPi}. 
It follows that  $Y_H(t )$ is a { \em bi-H\"older process}, whose spectrum is given by
\begin{equation}
 D_{Y_H} (h) = \left\{ \begin{array}{lll}
1 & \mbox{ if} & h = H , \\
1-H   & \mbox{ if}  & h = 2 H, \\ 
-\infty   & \mbox{ elsewhere.}  &  \end{array} \right. 
\label{equ-fbmsquared}
\end{equation}
Because the multifractal formalism yields the Legendre transform, i.e., the concave hull, of the multifractal spectrum, the leader based estimated multifractal spectrum 
would suggest practitioners to (wrongly) conclude that $Y_H(t )$ is a \emph{fully} multifractal function whose spectrum is supported by the interval $[ H, 2H ]$. 
This is illustrated in Fig~\ref{fig-FBMsquareFig} where it can be observed that the leader based measured multifractal spectrum yields a very precise estimation of the exact concave hull of $D_{Y_H} (h)$. 
Again, note that this is made possible only because negative values of $p$ can be used in the leader framework.

Further, the analysis of the $y=x^2$ point transform of fBm raises another disturbing issue. 
The multifractal properties of fBm are {\em homogeneous} (i.e., the spectrum measured from a restriction of the sample path to any interval $(a,b)$ on the real line is the same as that corresponding to the whole one).
This is no longer true for $(B_H)^2$: the spectrum measured on a restricted interval $(a,b)$ will vary depending on whether the interval includes or not a point where $B_H$ vanishes. 
This example shows that the multifractal properties are not preserved under  simple non linear transformations. 

\section{Real-World Applications}
\label{sec-App}

Reviewing the literature dedicated to real-world applications reveals that the concepts of scaling, scale invariance, selfsimilarity, fractal geometry and multifractal analysis have been and continue to be used to analyze data in numerous contexts of very different natures, ranging from natural phenomena
---  physics (hydrodynamic turbulence \cite{Frisch1995,Mandelbrot1974,my75}, astrophysics and stellar plasmas \cite{Starck1998}, statistical physics \cite{bk95,Keshner1982}, roughness of surfaces \cite{ard02}), biology (human heart rate  variabilities \cite{Ivanov2007,ivanov1999,kotani2003,lt05,lt96}, fMRI \cite{CIUCIU:2009:A,He2010}, physiological signals or images), geology (fault repartition, \cite{Foufoula94}) ---  to human activities --- computer network traffic \cite{pw00}, texture image analysis \cite{ard02}, population geographical repartition \cite{Dauphine2011,f97}, social behaviors  or financial markets \cite{CalvetFisher2008,Peters1994}.
B. Mandelbrot himself significantly contributed to the studies of many different applications, two of the most prominent possibly being hydrodynamic turbulence, the scientific field that actually gave birth to the essence of multifractal, (cf. e.g., \cite{Mandelbrot1974,Mandelbrot1990})  and finance (cf. e.g., \cite{Bacry2008,Bacry2010b,MandCal97,Mandelbrot1999}).

As a tribute to B. Mandelbrot's outstanding research work, a small number of applications
 from widely different origins have been selected to be revisited.
Their presentations are not intended as state-of-the-art reviews, but rather as illustrations, both of the theoretical arguments developed in the previous sections and of the benefits obtained from using multifractal analysis in these applications. 

\subsection{Hydrodynamic Turbulence}
\label{sec-turb}

\begin{figure}[h]
\centering
\includegraphics[width=1\linewidth]{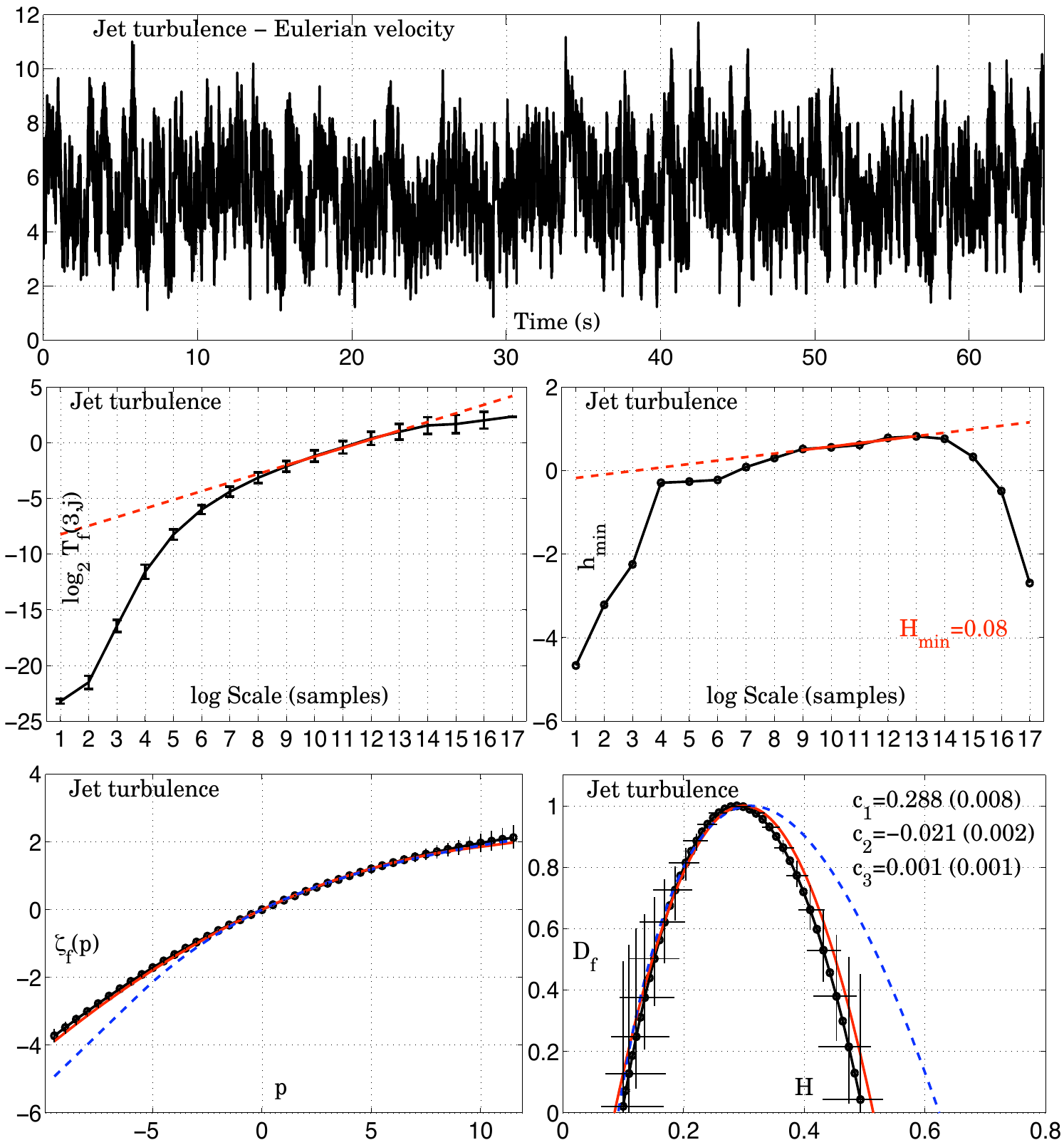}
\caption{\label{fig-TurbFig} {\bf Hydrodynamic Turbulence.} Top: longitudinal Eulerian turbulence velocity signal (hot-wire anemometry, Reynolds number: $\simeq 500$, Courtesy C. Baudet et al. \protect \cite{ruizetal95})). Middle: log-log plot for the wavelet leader based $T_f(3,j)$ (as in (\protect\ref{equ-WLSF})) and for the wavelet coefficients based $\Hmin$  (as in (\protect \ref{caracbeswav3hol})). Bottom: Wavelet leader scaling function (left) and multifractal spectrum (right) measured from a single run of data (black solid lines), compared to the log-normal (red solid lines) and log-Poisson (dashed blue lines) models. Data clearly select the log-normal model.}
\end{figure}

In the early $20$s, Richardson, combining empirical observations with the analysis of  Navier-Stokes equations, that mostly govern fluid flows, formulated a heuristic description of the apparently random movements observed in hydrodynamic turbulence flows, in terms of a cascade of energy: 
Between a coarse scale (where energy is injected by an external forcing) and a fine scale (where energy is dissipated by viscous friction), no characteristic scale can be
identified.
This seminal contribution  paved the road to tremendous amounts of works involving the concept of scaling to model turbulence flows (cf. \cite{Frisch1995,my75} for ancient and modern reviews, respectively). 
In 1941, Kolmogorov proposed a non explicit version of fBm to model velocity  (which implies that  $c_1=H$,  and $c_m \equiv 0$ for  $m \geq 2$)  \cite{k41}. 
However, a general consensus (a rare situation worth mentioning for the field of turbulence) amongst practitioners conducting real experiments rejected fBm as a valid model because the empirically observed scaling exponents (Kolomogorov scaling function) did not behave linearly in $p$, for $p \geq 0$. 
Framed in seminal contributions due to A. Yaglom \cite{Yaglom1966} and supported by physical arguments developed by 
Kolmogorov (and his collaborator A. Obukhov)   in 1962  \cite{kolmogorov62,obukhov62}, the energy transfer  from coarse to fine scales  has been modeled via split/multiply iterative random procedures that account for the physical intuitions beyond the vorticity stretching mechanisms implied by the Navier-Stokes equation. 
In 1974, B. Mandelbrot studied these various declinations of cascades and their properties \cite{Mandelbrot1974}, paving  the way for  their interpretation as a collection of singularities and hence to the first formulation of the multifractal formalism by G. Parisi and U. Frisch, based on the Kolmogorov scaling function \cite{ParFri85}.

This work triggered enourmous amounts of analysis of experimental and real-world turbulence data, aiming at measuring as precisely as possible  scaling exponents (i.e., the scaling function). 
The earliest contributions relied on the Kolmogorov scaling function (cf. e.g. \cite{arneodoetal96,arneodoetal98c,castaingetal90,cgm93,c89,pietropintoetal03}), and later work  initiated by A. Arneodo, E. Bacry and J.-F. Muzy  relied on the use of wavelets
 (cf. e.g. \cite{dma01,LASHERMES:2008:A,meneveau91,muzyetal91,muzyetal94,WENDT:2007:E}).

Fig.~\ref{fig-TurbFig} (top plot) shows an example of a turbulence longitudinal Eulerian velocity signal, measured with hot-wire anemometry techniques. 
It has been collected on a jet turbulence experiment \cite{ruizetal95}, at average Reynolds number $\simeq 580$, and has been made available to us by C. Baudet and collaborators (LEGI, Universit\'e Joseph Fourier, INPG, CNRS, Grenoble, France), who are gratefully acknowledged.
Middle row plots depict its scaling property in the so-called inertial range, a priori defined from the physical characteristics of the flow conditions: Left, the wavelet leader based $\log_2 T_f(j,q)$ versus $\log_2 2^j$ (for $p=3$), and $\Hmin$ on the right.
In the bottom row the measured wavelet leader scaling function and multifractal spectrum are plotted, clearly indicating multifractality.  

While the generic multifractality of turbulence data (understood as a strict departure from linearity in $p$ of the scaling function)  has never been called into question by any of these studies, a major open issue  remains: Which precise cascade model better describes turbulence? 
This question is not purely theoretical, since the precise ingredients entering the cascade that best match data may reveal physical mechanisms of importance for turbulence. 
Two such models are emblematic of this issue:  In 1962, Obukhov and Kolmogorov \cite{kolmogorov62,obukhov62} proposed a model mostly based on a law of large numbers argument and referred to as the log-normal multifractal model. It predicts that $ \zeta_f(p) = c_1p + c_2p^2/2$ and hence that $c_m \equiv 0$ for $m \geq 3$. 
In 1994, She and L\'ev\^eque \cite{sheleveque01} proposed an alternative construction based on the key ingredient that energy dissipation gradients must remain finite within turbulence flows. It is referred to as the log-Poisson model and yields a multifractal process with all non zero $c_m$s. 

Fig.~\ref{fig-TurbFig} (bottom plots) compares the scaling function (left) and multifractal spectrum (right) measured from data (solid black), together with bootstrap based confidence intervals (cf. \cite{WENDT:2007:E} for details) against those predicted by the log-normal (solid red) and log-Poisson (dashed blue) model.
Together with the fact that the measured $c_3$ can not be distinguished from $0$ ($c_3 = 0.001 \pm 0.001$), these plots clearly indicate that the log-normal model is preferred by turbulence data. 
This has been studied systematically and on larger data sets in \cite{LASHERMES:2008:A,WENDT:2007:E} and unambiguously confirms earlier findings reported in \cite{arneodoetal98c,dma01}.

\subsection{Finance Time Series: Euro-USD rate fluctuations}
\label{sec-finance}

\begin{figure}[h]
\includegraphics[width=1\linewidth]{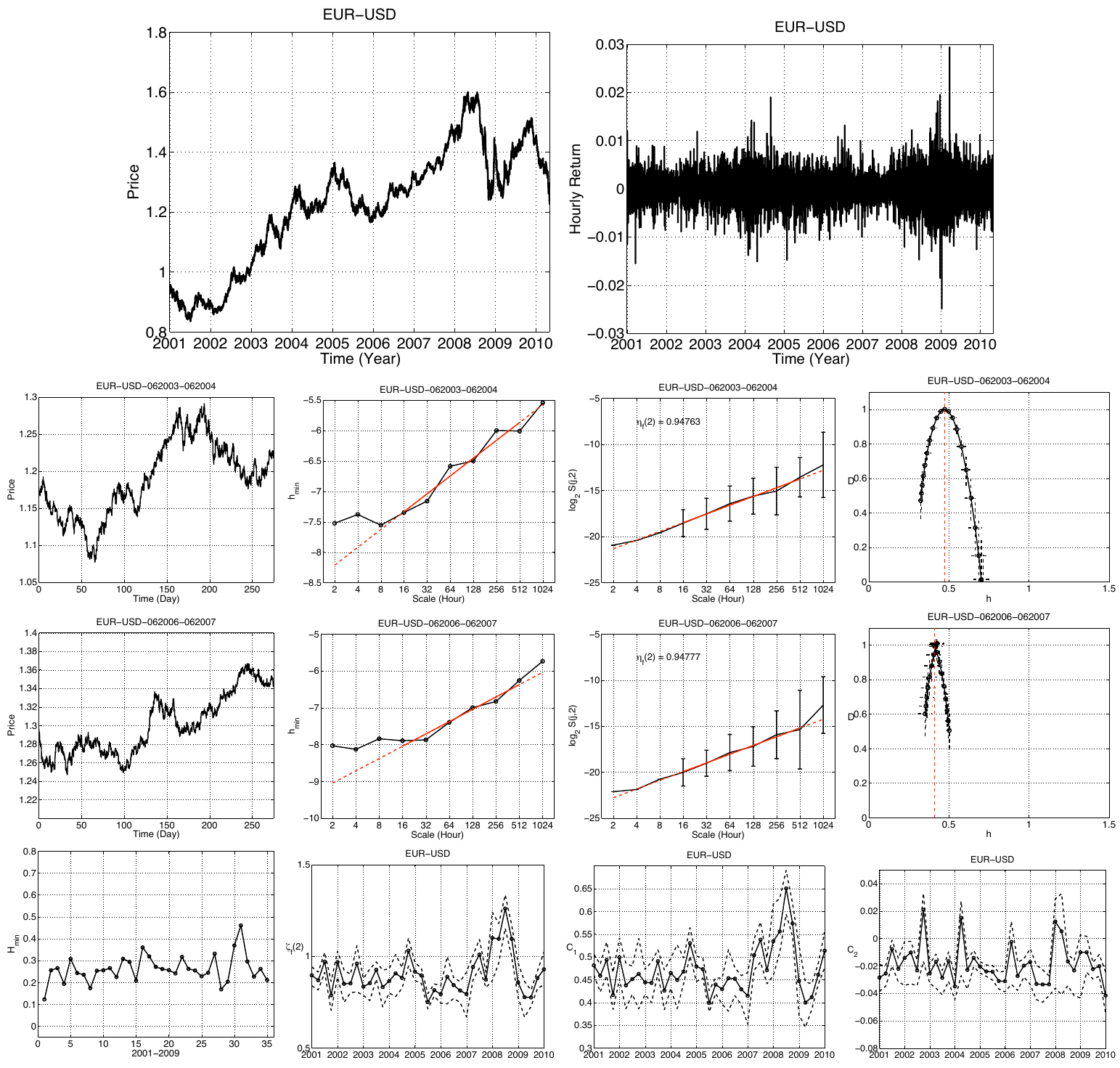}
\caption{\label{fig-finance} {\bf Euro-USD.} Top: Euro-USD rate fluctuations for 10 years (starting January 1st, 2001) (left) and hourly returns (right). 
Data Courtesy Vivienne Investissement, Lyon, France.
One year long time series (June 2003 to June 2004, 2nd row;  June 2006 to June 2007, 3rd row), with log-log plots for the wavelet coefficients based $\Hmin$  (cf. (\protect \ref{caracbeswav3hol})), log-log plots for the wavelet leader based $T_f(2,j)$ (given by (\protect \ref{equ-WLSF})), and 
multifractal spectra, measured using the wavelet leader multifractal formalism. 
Bottom row: Estimation in one year long time windows (with a $3$-month sliding shift) of $\Hmin$, $\eta_f(2)$, $c_1$, $c_2$ (together with confidence intervals).
These plots show that while $\Hmin >0$, $\zeta_f(2)$ often remains below the finite quadratic variations limit of $1$. 
Parameter $c_1$, $c_2$ variations indicate that the multifractal properties of the Euro-USD price fluctuations vary significantly along time.}
\end{figure}


Modeling stock market prices received considerable efforts and attention starting with  the seminal early contribution of Bachelier in 1900 who, in his thesis,  based his description of  stock market on Brownian motion, i.e., on an ordinary Gaussian random walk \cite{Bachelier1900}.
Essentially, this model implies that  daily \emph{returns}, i.e.,  increments of (or differences between) prices taken along two consecutive days, $X(t) = P(t)-P(t-1)$, consist of independent, identically distributed Gaussian zero-mean random variables. 
This hence leaves little room to earn money from stock market investments by pure statistical model based prediction. 
However, it has long been observed and/or expected that stock market prices significantly depart from the Brownian  model: 
First, returns may exhibit statistics that significantly depart from Gaussian distributions.
Second, it is often observed that volatility, a quantity essentially based on the squared (or absolute) value of the returns, display some form of dependency.
Notably, in their celebrated contribution, Granger and Joyeux indicated that volatility is well characterized by long memory, or long term dependence \cite{Granger1980,GrangerJoyeux1980}, a model that is closely related with selfsimilarity and fBm. 
This motivated substantial efforts to model stock market returns with the GARCH models family or FARIMA alternatives (cf. e.g., \cite{Beran1994,DittmanGranger2002,TEYSSIERE:2007:A}). 

B. Mandelbrot significantly contributed to the analysis of the statistical properties of stock market returns by performing an iconoclastic analysis of commonly used finance models, emphasizing heavy tails, as opposed to Gaussian distributions (the Noah effect), and long memory, as opposed to independence (the Joseph effect), cf. \cite{Mandelbrot1997} for a review, and several seminal contributions of B. Mandelbrot; see also e.g., \cite{E14,E20}. 
In particular, this called into question the Black and Scholes model strongly relying on Gaussianity,
and led B. Mandelbrot and collaborators to explicitly propose that stock market prices may have multifractal properties in  \cite{FisherCalvetMandelbrot1997} in contradistinction with most It\^o-martingale based models used so far.
B. Mandelbrot's work paved the way towards numerous efforts to assess multifractality empirically in stock market prices and exploit it (cf. e.g., \cite{Bacry2010,Bacry2008,CalvetFisher2002a,CalvetFisher2008,Peters1994}).

Almost one century after Bachelier, B. Mandelbrot and collaborators proposed a new model for stock market prices: fBm in multifractal time \cite{MandCal97,Mandelbrot1999}, a process that explicitly incorporates multifractal properties by a multifractal time jittering, obtained by a multiplicative cascade construction. 
In this model, a new generation of cascades is used, based on Compound Poisson distributions, following again an original idea of B. Mandelbrot \cite{BJM10,BaM02,Barral2002}. 
Let  $F$  denote the distribution function of the multiplicative cascade $\mu$ on $[0,1]$: $ F(t) = \mu ([0,t)).$
Then, fBm in multifractal time is defined as
\begin{equation} 
\label{equ-fbmmf}
X_t = B_H (F(t)), 
\end{equation}
where $B_H$ stands for fBm with selfsimilarity parameter $H$ (which is independent of $\mu$).  Note that, as  a consequence of (\ref{trandormdim}), the spectrum of $X$ is deduced from the spectrum of $F$ by a simple dilation along the $h$ axis
\[ÊD_X(h) = D_F (h/H);  \] 
see the results of R. Riedi and S. Seuret \cite{Riedi2003,Seur3} for additional discussions  concerning { \em multifractal subordination}, i.e. understanding multifractal properties of  composed processes).
Numerous  rich and interesting declinations of such constructions of cascades and multifractal processes have been proposed since (cf. e.g., 
\cite{Bacry2001,Bacry2010,Bacry2008,Bacry2002,Bacry2003,CalvetFisher2008,Chainais2005,Chainais2005a}).
Note that the idea of { \em multifractal subordination} (i.e., composition by a multifractal change of time) is extremely powerful, and can be applied in  many settings. 
Let us mention just one example: In \cite{Polya}, B. Mandelbrot and one of the authors   performed the multifractal analysis of the Polya function (a famous example of  continuous space-filling function)  by simply remarking that it can be \emph{factorized} as the composition of the repartition function of a simple dyadic cascade  with  a monoh\"older function.

As a tribute to B. Mandelbrot's contributions, though the authors have not personally contributed to the field, the multifractal properties of the Foreign exchange Euro-USD price fluctuations are now studied. 
Fig.~\ref{fig-finance} (top row) shows the hourly Euro-USD quotations (and hourly returns) since the birth of Euro (Data Courtesy Vivienne Investissement, Lyon, France, {\tt /www.vivienne-investissement.com/}). 
The return (or increment) time series clearly shows large fluctuations that are not consistent with a constant along time Gaussian distribution.
Scaling properties analyzed here are found to hold for scales ranging from the day to the month, and hence correspond to long term characteristics, as opposed to intra-day analysis, not conducted here.
Multifractal spectra measured over one year long time series (second and third rows, left plots), using the wavelet leader multifractal formalism, show that the Euro-USD quotations exhibit multifractal spectra that are in agreement with the fBm in multifractal time model, yet indicate that the multifractal properties may vary along the years.
To further analyze these variations, scaling attributes are measured in one year long time windows (with a $3$-month sliding time shift) and reported  (together with confidence intervals) in Fig.~\ref{fig-finance} (bottom row).
Interestingly, it is found that $\Hmin$ is permanently above $0$: $\Hmin >0$, which clearly rules out the use of jump processes to model ForEx time series.
Moreover, because $ \Hmin >0$, the finiteness of quadratic variations can be evaluated using the leader scaling function at $p=2$, $\zeta_f(2)$ (as in (\protect \ref{equ-WLSF}), cf. Section~\ref{secleade}). 
Fig.~\ref{fig-finance} (bottom row, 2nd plot) shows that $\zeta_f(2)$ fluctuates and often remains below the finite quadratic variations limit of $1$: $\zeta_f(2)<1$.
This result validates the hypothesis of  the finiteness of quadratic variations, underlying the use of It\^o-martingale processes as models to describe the Euro-USD variations, at least when long term variations are analyzed, as is the case here; 
this does not exclude that shorter term (or intra-day) variations might show different properties.
In Fig.~\ref{fig-finance} (bottom row, 3rd and 4th plots), parameters $c_1$, $c_2$ are plotted as functions of time.
These plots furthermore indicate that though a fBm in multifractal time-like process provides a valid model to describe the price fluctuations, the multifractal properties may vary along time.
Two different regimes are suggested by these plots: Before 2007, with a $c_1 \simeq 0.45 < 0.5$ and $c_2 \simeq -0.02$ ; and after 2008, with $c_1 \simeq 0.5$ and $c_2 \simeq -0.01$. 
Recalling that the ordinary Brownian motion (proposed by Bachelier) would correspond to  $c_1 = 0.5$ and $c_2 = 0$, it may be concluded that the (long term) fluctuations of the Euro-USD quotations show less structure for the later period than for the former, hence making predictions for investments more difficult (see also discussion in Section~\ref{homog} on homogeneous multifractal properties). 

Assuming that the fBm in multifractal time model studied by L. Calvet, A. Fisher and B. Mandelbrot is the correct description of prices, let us now show  how the parameter $H$  of the corresponding fBm can be inferred from the observed data using the wavelet scaling function (the heuristic  argument that follows can easily be turned into a mathematical proof).
Let us first return  to the definition of the Kolmogorov scaling function (\ref{kolmo}). 
Approximating the integral by a Riemann sum,  and taking $h = 2^{-j}$, one obtains
\[ 
\int |X_{t+h} - X_t|^p dt \sim 2^{-j} \sum_k  \left| X_{\frac{k+1}{2^j}} - X_{\frac{k}{2^j}}\right|^p  = 
2^{-j} \displaystyle\sum_k  \left| B_H \left(F\left( {\frac{k+1}{2^j}}\right)\right) - B_H \left(F\left({\frac{k}{2^j}}\right) \right) \right|^p \]
which, using the fact that $| B_H (t) -B_H (s) | \sim |t-s|^{H}$, is of the order of magnitude of 
\[ 2^{-j} \sum_k  \left| F\left( {\frac{k+1}{2^j}}\right) - F \left({\frac{k}{2^j}}\right)  \right|^{pH} .\] 
In particular, for $p= 1/H$, using the fact that $F$ is continuous and increasing, we obtain 
\[ 2^{-j} \sum_k  \left| F\left( {\frac{k+1}{2^j}}\right) - F \left({\frac{k}{2^j}}\right)  \right| = 2^{-j} \sum_k  \left( F\left( {\frac{k+1}{2^j}}\right) - F \left({\frac{k}{2^j}}\right) \right) = 2^{-j}  \] 
Since this quantity is, by definition, of the order of magnitude of  
$  (2^{-j})^{\eta_f(p)} $, we obtain that the exponent $H$ of the fBm  used in this  model (cf. \cite{MandCal97}) satisfies 
\BE  \label{CalMan} \eta_f\left(\frac{1}{H}\right) =1.  \EE 
Using the fact that the Kolmogorov and the wavelet scaling function coincide in the range of exponents involved here, (\ref{CalMan}) allows to recover the selfsimilarity  exponent  of the fBm  in this model from the wavelet scaling function.  In the example considered in Fig.~\ref{fig-finance}, one obtains that, in practice, the parameter $H$ thus estimated is very close to $c_1$. 

\subsection{Fetal Heart Rate Variability}
\label{sec-fhbv}

\begin{figure}[h]
\centering\includegraphics[width=0.8\linewidth]{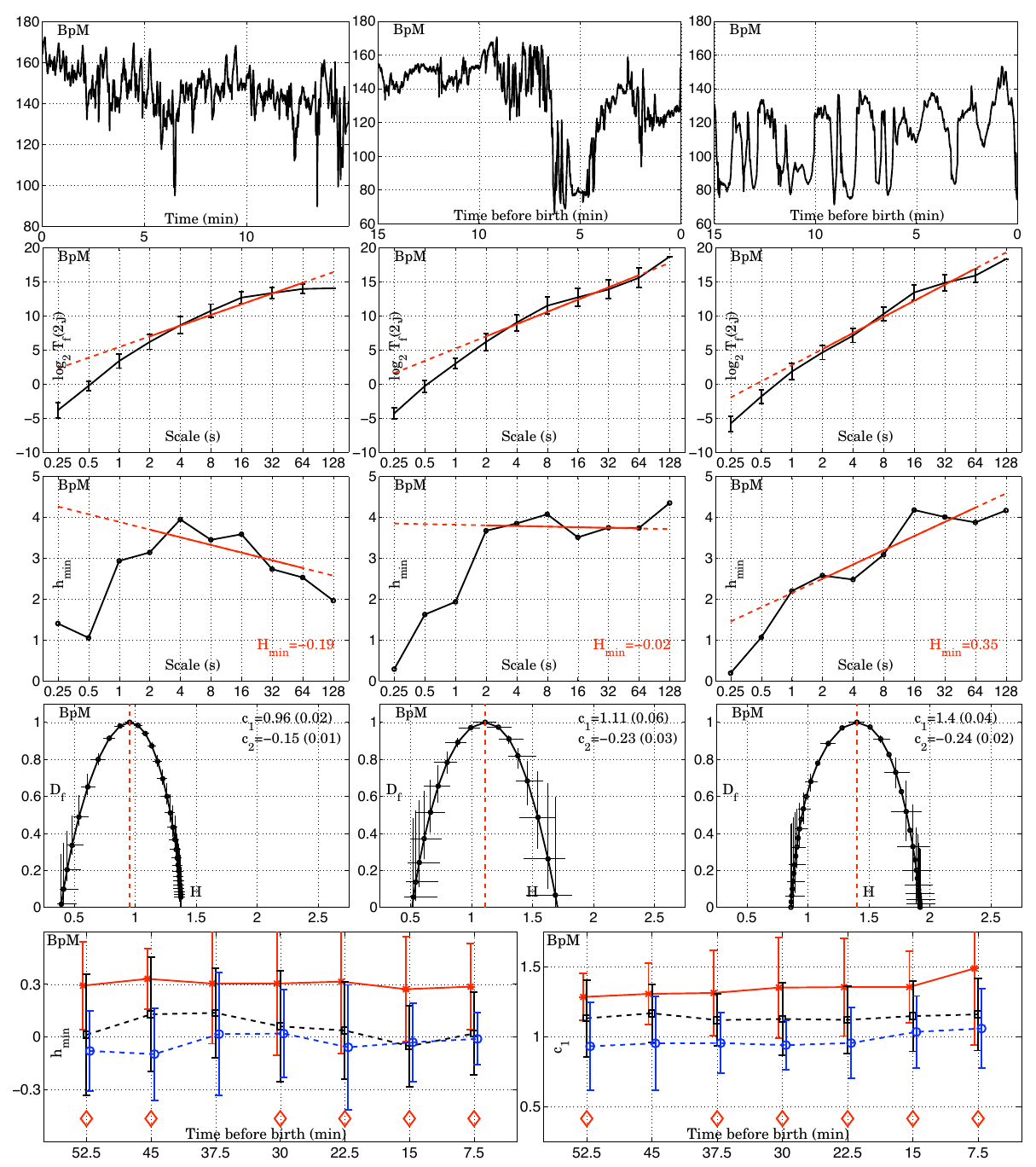}
\caption{\label{fig-fECG} {\bf Fetal-ECG.} Top: Examples of cardiac frequency time series (in Beats-per-Minute) for fetuses during delivery process - last $15$ min. before delivery, for a healthy subject, correctly diagnosed as such by the obstetrician (left) (True Negative), for a healthy subject, incorrectly diagnosed as suffering from per partum asphyxia (middle) (False Positive), for a fetus actually suffering from per partum asphyixia and diagnosed as such (True Positive). Data Courtesy M. Doret, Obstetrics Dept. of the University Hospital Femme-M\`ere-Enfant, Lyon, France.
Second Row: log-log plots for the wavelet leader based $T_f(2,j)$ (cf. (\protect \ref{equ-WLSF})).
Third Row: log-log plots for the wavelet coefficients based $\Hmin$  (as in (\protect \ref{caracbeswav3hol})).
Fourth Row: Multifractal spectra, measured using the wavelet leader multifractal formalism.
Bottom Row: Parameters $\Hmin$ (left) and $c_1$ (right) estimated over $15$ min. long sliding windows (with a $7.5$min overlap) over the last hour before delivery, values are averaged over the classes of $15$ patients each: True Positive (dashed blue line with 'o'), False Positive (dashed black line with squares, True Positive (solid red line with $\ast$). The red $\diamond$ correspond to windows where a Wilcoxon ranksum test between True Positive and True Negative gives a $p-$value below $5 \%$. This analysis shows that unhealthy subjects (True Positive) systematically display less HRV as compared to healthy ones (True Positive and False Positive) (cf. \cite{ABRY:2010:A,Doret2011}). 
}
\end{figure}

Heart rate variability (HRV), which refers to the short time fluctuations (within the minute) of  heart rate, has long been considered as a powerful tool to characterize the health status of a given subject: In a nutshell, the more variability, the healthier the subject is (cf. e.g., the seminal contribution \cite{Akselrod1987}). 
To characterize variability in heart rate time series, spectral analysis has long been used to measure the so called LF/HF balance, i.e., the ratio of energies measured in frequency bands attached to the sympathic and parasympathic activities of the autonomous central nervous system \cite{Akselrod1987}. 
Yet, it has long been observed that heart rate time series display \emph{fractal} properties and that the corresponding fractal exponent could be used as a non invasive tool for non healthy status detection (cf. \cite{Akselrod1987}). 
This opened the way to numerous analyses aiming at describing HRV data with fBm-like models and at correctly estimating the corresponding selfsimilarity  parameter $H$;  
later, multifractal analysis has naturally been considered as an analysis tool potentially improving the characterization of HRV scaling properties (cf. e.g., \cite{Ivanov2007,ivanov1999,kotani2003}).
Recently, multifractal analysis has been involved in the detection of per partum fetal asphyxia. 
Obstetricians are monitoring fetal-HRV during the delivery process to detect potential fetal asphyxia. 
Using a set of analysis criteria (the international FIGO rules), obstetricians measure the health status of the fetus      and decide on operative delivery or not.
A major difficulty they      are facing is the high number of False Positive (fetuses diagnosed as unhealthy during the delivery process but a posteriori evaluated as perfectly healthy), and hence of operative delivery that might  have been avoided.
This constitutes an important public health issue, as operative deliveries are accompanied by a significant number of serious posterior complications for both  newborns and their mothers.

Fig.~\ref{fig-fECG} (top row) shows examples of heart rate time series (in Beats-per-Minute) for a healthy subject, correctly diagnosed as such by the obstetrician (left) (True Negative), for a healthy subject, incorrectly diagnosed as suffering from per partum asphyxia (False Positive), and for a fetus actually suffering form per partum asphyixia and diagnosed as such (True Positive). Data Courtesy M. Doret, Obstetrics Dept. of the University Hospital Femme-M\`ere-Enfant, Lyon, France,  \cite{ABRY:2010:A,Doret2011}.
As shown in Fig.~\ref{fig-fECG} (second and third rows),  for all three classes (True Negative, True Positive, False Positive), HRV exhibits scaling properties for scales ranging from the second to the minute. 
The fourth row in Fig.~\ref{fig-fECG} shows multifractal spectra, estimated from $15$ min. long heart rate time series using wavelet leaders, which suggests that the healthier the subject, the larger its variability in the sense that its $c_1$ parameter is smaller (spectrum shifted to the left).
To confirm such an observation, multifractal parameters are  estimated within $15$ min. long sliding windows (with a $50 \%$ time overlap) for the last hour before delivery. 
The $15$ min. window duration is chosen after the obstetrician's request that asphyxia should be detected as early as possible. 
Multifractal parameters are then averaged within the three classes (each containing $15$ patients carefully selected as representative by obstetricians). 
A Wilcoxon ranksum test is conducted for each time window between the True Positive and False Positive classes.  
The results reported in Fig.~\ref{fig-fECG} (bottom row) clearly indicate a significant increase of $c_1$ and $\Hmin$ (hence a significant decrease in HRV) for the Unhealthy True Positives compared to the Healthy True Negatives and False Positives. 
The p-value of the Wilcoxon ranksum test is observed to be below $5 \%$ (cf. red $\diamond$  in Fig.~\ref{fig-fECG} (bottom row)) often $30$ min earlier than the time the decision to operate the delivery was actually taken. 
These promising results are reported in details in \cite{ABRY:2010:A,Doret2011} and will be further investigated.

It is worth noting that for heart rate time series (as for many other types of human rhythms signals, cf. e.g., fMRI applications \cite{CIUCIU:2009:A,He2010}), $\Hmin$ can be found either positive or negative depending on the subject, or its health status. 
This is consistent with the fact that practitioners, in applications, noticed that they often had to switch from fBm to fGn  (or conversely) to model data and understood this as a puzzling inconsistency (see for instance discussions in \cite{He2010}).
To some extent, the general framework developed here alleviates such inconsistency by introducing $\Hmin$ as one of the various parameters that can be used to analyze scaling without a priori  assumptions on the   data. 
In the analysis of fetal HRV described above, a fractional integration of order $1/2$ (larger than the absolute value of the most negative $\Hmin$ encountered in the whole database) has been applied to enable the use of the wavelet leader multifractal formalism in a consistent manner on the whole set of data.

\subsection{Aggregated Internet Traffic}
\label{sec-Internet}

\begin{figure}[h]
\includegraphics[width=1\linewidth]{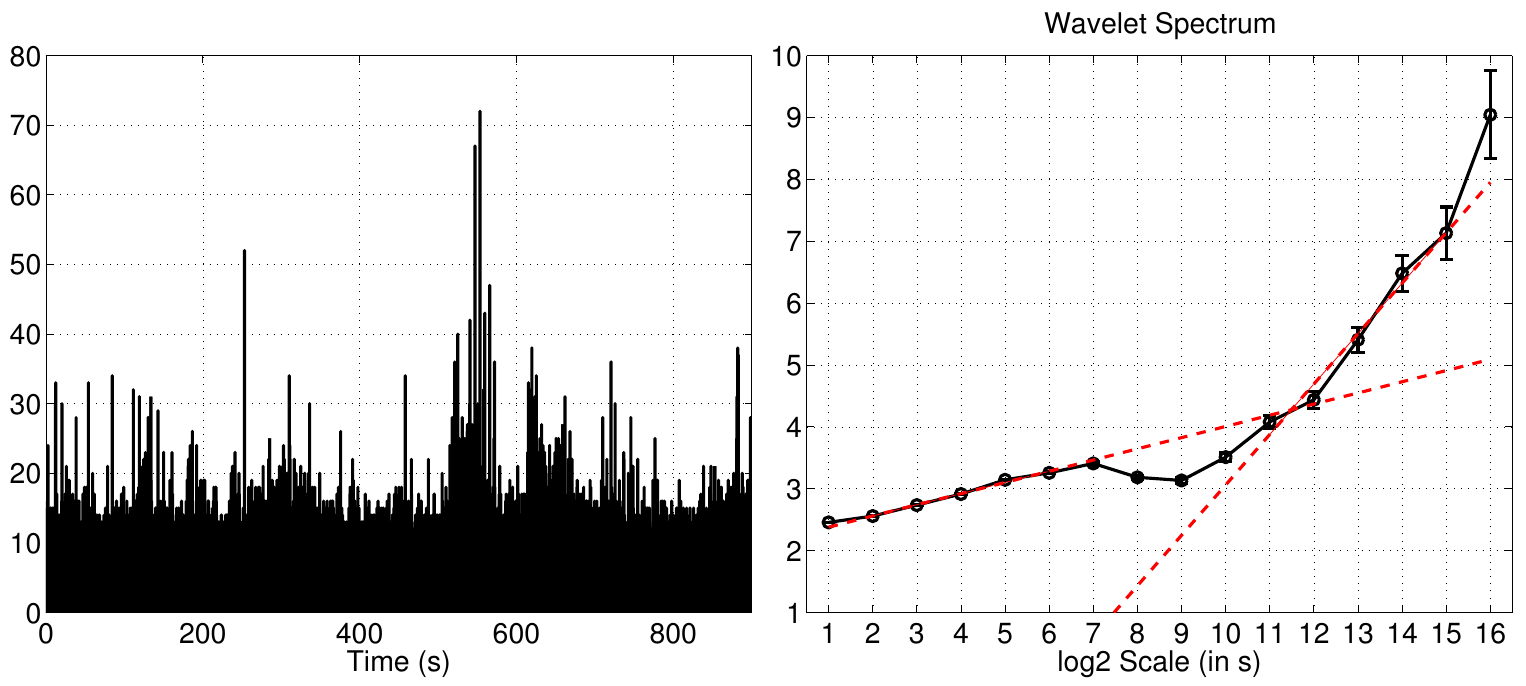}
\caption{\label{fig-TraffFiga} {\bf Internet Traffic: Aggregated time series and scaling.} Left: Aggregrated time series (Courtesy MAWI research project, \protect \cite{mawi}); Right: log-log plot for the wavelet coefficient based $S_f(2,j)$ (as in \protect (\ref{equ-WSF})).
It clearly displays two independent scaling ranges, below and above the second.}
\end{figure}

\begin{figure}[h]
\includegraphics[width=1\linewidth]{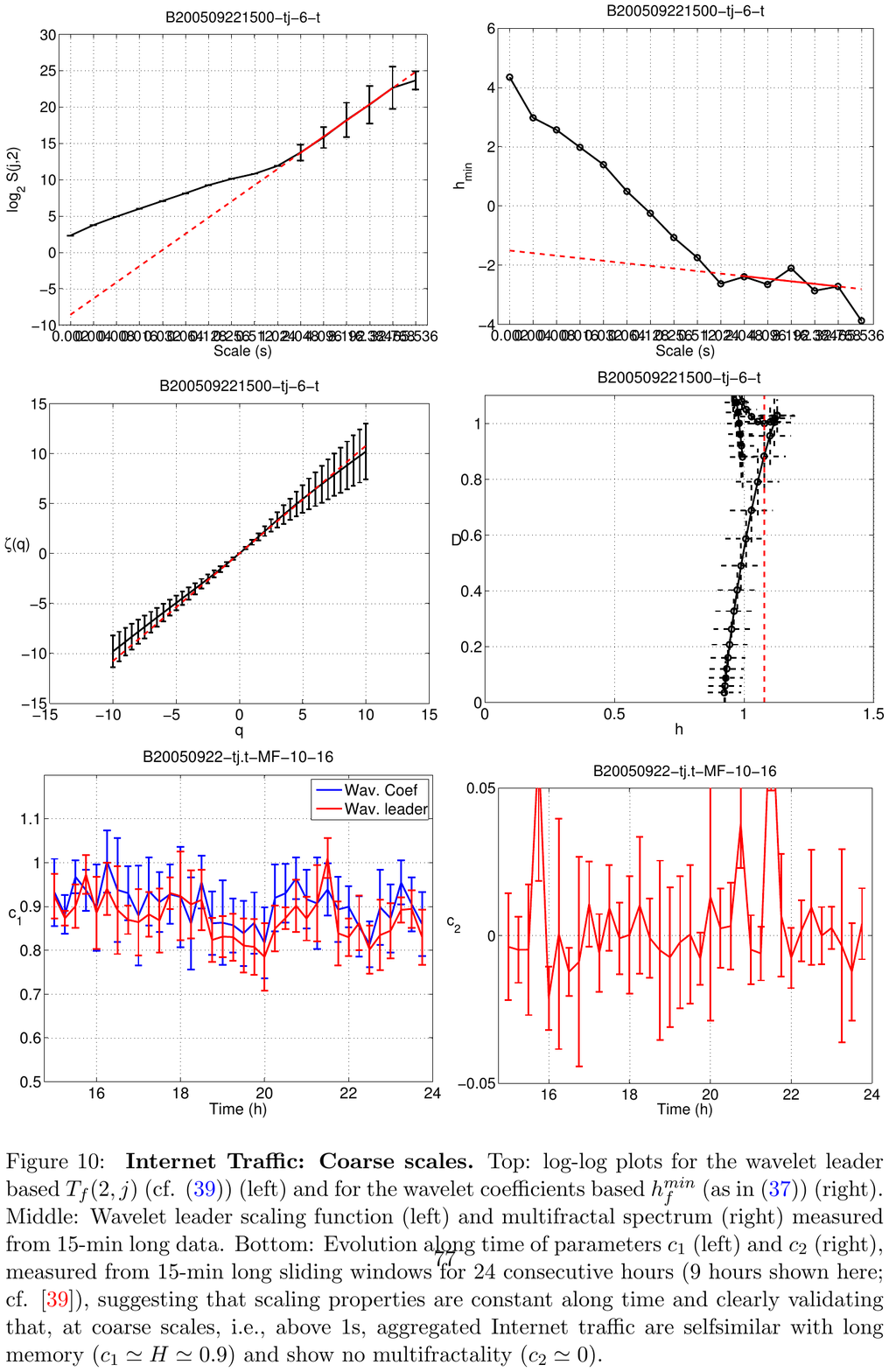}
\caption{\label{fig-TraffFigb} {\bf Internet Traffic: Coarse scales.} Top: log-log plots for the wavelet leader based $T_f(2,j)$ (cf. (\protect \ref{equ-WLSF})) (left) and for the wavelet coefficients based $\Hmin$  (as in (\protect \ref{caracbeswav3hol})) (right). 
Middle: Wavelet leader scaling function (left) and multifractal spectrum (right) measured from $15$-min long data.
Bottom: Evolution along time of parameters $c_1$ (left) and $c_2$ (right), measured from $15$-min long sliding windows for $24$ consecutive hours ($9$ hours shown here; cf. \protect \cite{BORGNAT:2009:A}), suggesting that scaling properties are constant along time and clearly validating that, at coarse scales, i.e., above $1$s, aggregated Internet traffic are selfsimilar with long memory ($c_1 \simeq H \simeq 0.9$) and show no multifractality ($c_2 \simeq 0$).}
\end{figure}

\begin{figure}[h]
\includegraphics[width=1\linewidth]{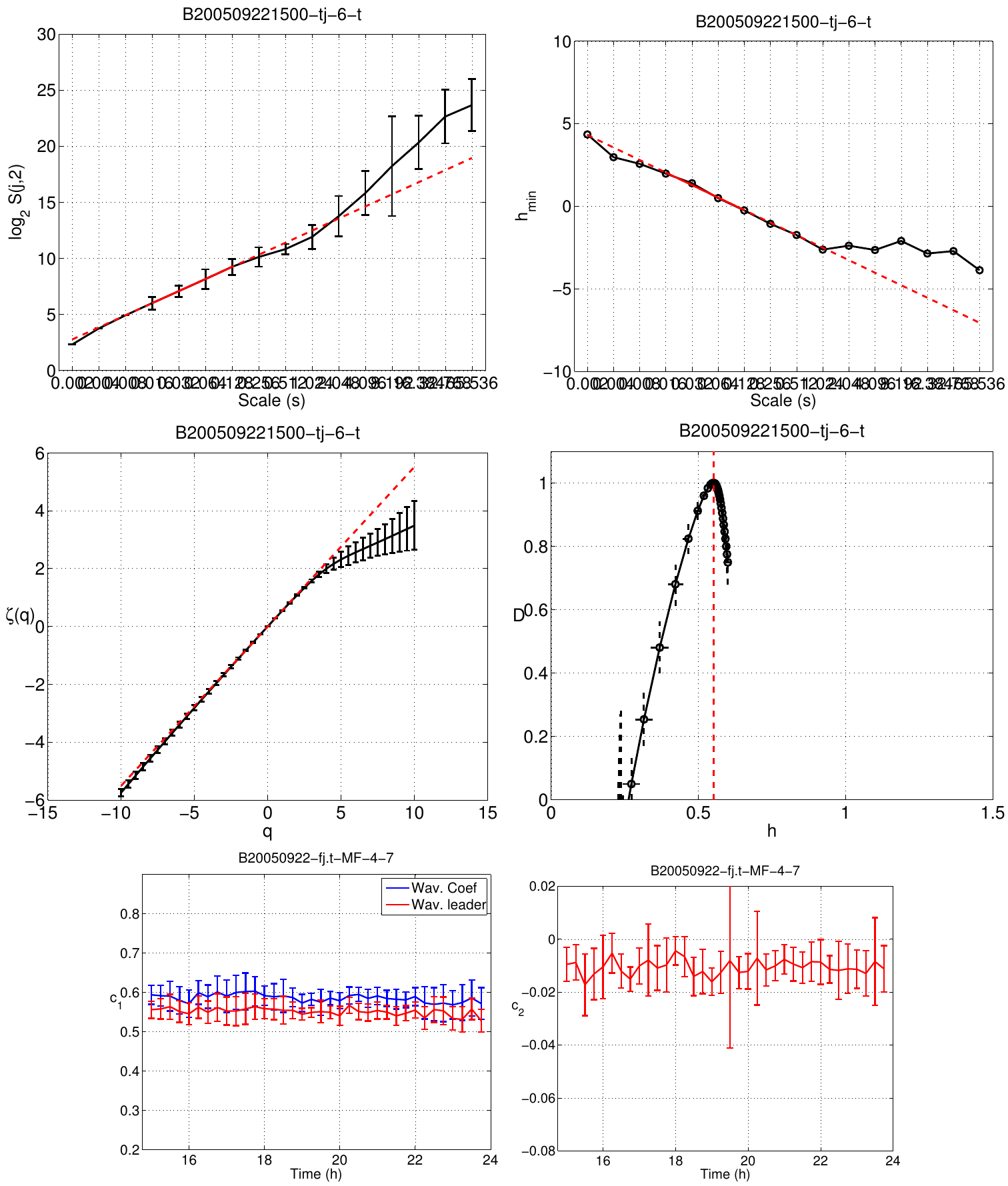}
\caption{\label{fig-TraffFigc} {\bf Internet Traffic: Fine Scales.} Top: log-log plots for the wavelet leader based $T_f(2,j)$ (as in (\protect \ref{equ-WLSF})) (left) and for the wavelet coefficients based $\Hmin$  (given by (\protect \ref{caracbeswav3hol})) (right). 
Middle: Wavelet leader scaling function (left) and multifractal spectrum (right) measured from $15$-min long data.
Bottom: Evolution along time of parameters $c_1$ (left) and $c_2$ (right), measured from $15$-min long sliding windows for $24$ consecutive hours ($9$ hours shown here). (cf. \protect \cite{BORGNAT:2009:A}), also suggest that scaling properties are constant along time and indicate that, at fine scales (i.e., below $1$s), aggregated Internet traffic are characterized by $c_1 \simeq H \simeq 0.6$ and $c_2 \simeq 0.005 \pm 0.005 $, revealing at best a very weak multifractality or no multifractality at all.}
\end{figure}

Internet traffic engineering  currently  is a challenging task involving a variety of issues ranging from economy development planning and security assessment to statistical data characterization. 
The earliest versions of statistical models used for Internet traffic modeling were stemming from the experience gained on the previous large communication network: Telephone communications. 
B. Mandelbrot himself made very early  contributions  to such modeling (cf. e.g., \cite{N6}), and this certainly was a major source of inspiration  for the later developments he made  in the theory of fractal stochastic processes. 
However, such telephone communications models were all based on the assumption that there exist one or several characteristic scales of time in data, while it has soon been recognized that Internet traffic is well characterized by long memory and self similarity (cf. \cite{Abry1998,leland94,leland93,Park2000,Paxson1995} for  seminal articles).  
This later raised the question of multifractality of Internet traffic (cf. e.g., \cite{abfrv02,Feldmann1998b,HOHN:2005:A,Riedi1999}) that we briefly illustrate here. 

Any activity performed on the net (emailing, web browsing, Video-On-Demand, data downloading, IP telephony,\ldots) amounts to establishing a connection between hosts (server-clients, peer to peer, \ldots), consisting of computers located all over the world, according to an agreed-on protocol, and exchanging a (very) large number of elementary quanta of information, the Internet Protocol packets (in short, IP Pkts).
Collecting Internet traffic data hence naturally consists of recording the IP Pkt time stamp and header (IP source, IP destination, Port Source, Port Destination). 
Obviously, analyzing the marked point process resulting of the tremendous number of IP Pkts exchanged is beyond feasible, and it is classical to instead monitor so-called aggregated (or count) time series: Practitioners select a given aggregation-level $\Delta$, whose choice depends on the application, and construct either the Pkt or Byte time series that counts the number of IP Pkts whose timestamp falls within $[k\Delta, (k+1)\Delta)$, or the corresponding volume in Byte. 

An example of such Pkt time series is plotted in Fig.~\ref{fig-TraffFiga} (left). 
It consists of traffic from the MAWI repository, made available by the Japanese WIDE research program \cite{mawi}.
The MAWI repository consists of data that have been collected everyday for $15$min at $2$ pm (Tokyo time) since January 2001 on major backbone links ($100$Mbps or $1$Gbps) 
connecting Japan and the United States (cf. \cite{BORGNAT:2009:A,Cho2000,DEWAELE:2007:A}. 

The left plot of the same figure displays the wavelet coefficient based $\log_2 (S_f(2,j))$ vs. $\log_2 (2^j) =j$ (as in \protect (\ref{equ-WSF})) and clearly shows that, for Internet traffic, there exist two ranges of scales for which scaling can be identified: Coarse scales, above $1$s, and fine scales, below $1$s. 
For each of those ranges of scales, the wavelet leader based multifractal formalism is applied. 

For the coarse scales, Fig.~\ref{fig-TraffFigb} clearly indicates on a $15$min long set of data aggregated at $1$ms that the leader scaling function is linear, and that the multifractal spectrum collapses on a single point, hence suggesting monofractality or that multifractality is at best very weak. 
To further investigate this point, the parameters $c_1$ and $c_2$ are estimated on sliding $15$min long windows for a trace collected during $24$ hours (on Sept., 22nd, 2005), for scales ranging from $1$s to $1$ min.
Results are reported for the last $8$ hours for simplicity and show, first, that the estimated parameters are remarkably constant along the day (while the changes in amplitude of the time series may be drastic) and, second, that $c_2 \simeq 0$. 
This indicates no multifractality at coarse scale and, because $c_1 \sim H \simeq 0.9$, a significant long memory in the data. 
These experimental results corroborate those reported in \cite{HOHN:2005:A} and are consistent with the main model due to Taqqu et al. (cf. \cite{taqqu97}) relating selfsimilarity and heavy-tails of the web file size distribution,   see also \cite{Abry2009,Hohn:2003:A,Loiseau2010}). 

For the fine scales, Fig.~\ref{fig-TraffFigc} shows on the same $15$min long set of data aggregated at $1$ms that the leader scaling function is linear for $p$ around $0$, hence indicating a very weak multifractality or no multifractality at all. 
Parameters $c_1$ and $c_2$, estimated on sliding $15$min long windows for the $24$ hour trace, again turn out to be constant along time, with $c_1 \simeq 0.55 \pm 0.05$ and $c_2 \simeq 0.007 \pm 0.005$. 
This hence validates that fine scales of aggregated traffic show at best little multifractality.
Furthermore, estimations of the selfsimilarity parameters (not shown here) being close to $ 0.5$ also indicate little correlation. 
These analyses tend to show that fine scales of aggregated traffic are characterized by the quasi-absence of dependence structure.
Moreover, the estimated parameters $c_1$  and $c_2$ remain remarkably constant along time for the whole day and across the various days analyzed, a highly unexpected result for Internet traffic, which is often described as widely varying, and which is continuously subject to numerous occurrences of anomalies and attacks.
This is essentially because data have been pre-processed using a random projection and a median procedure that robustify estimation against anomalous behaviors (cf. \cite{BORGNAT:2009:A,DEWAELE:2007:A}). In practice, the effect of this procedure is to eliminate  abnormal components in the traffic. 
The results reported here can hence be associated with the properties of regular and legitimate background Internet traffic, corresponding to the usual Internet user behaviors. 
Furthermore, this quasi-absence of dependence structure at fine scales  is consistent with the Cluster Point Process modeling of Internet traffic aggregated count time series proposed in \cite{HOHN:2005:A}, which implies that fine scales are actually not possessing strictly scaling properties. 
However, it also shows that the statistical properties of the fine scales are satisfactorily approximated by scaling behaviours.
Therefore, the fact that regular or legitimate traffic does not show multifractality does not imply however that multifractal parameters like $c_1$ and $c_2$ can not be used for traffic monitoring, such as anomaly detection for instance. 
Indeed, the occurrence of anomalies is likely to change the dependence structure of the fine scale statistics (as shown in e.g., \cite{DEWAELE:2007:A}) and consequently the empirically estimated parameters $c_1$ and $c_2$, which may thus prove very useful, even if scaling properties are measured that actually consist of approximations of the true statistical properties of the fine scales. 

Finally, let us mention that for aggregated count  time series $\Hmin <0$, as shown on the top left plots of  Figs.~\ref{fig-TraffFigb} and \ref{fig-TraffFigc}, hence a fractional integration of order  $1$ is systematically applied prior to performing the wavelet leader multifractal formalism.

\subsection{Texture classification: Vincent  Van Gogh meets Beno\^{\i}t Mandelbrot}
\label{sec-VG}

\begin{figure}[h]
\includegraphics[width=1\linewidth]{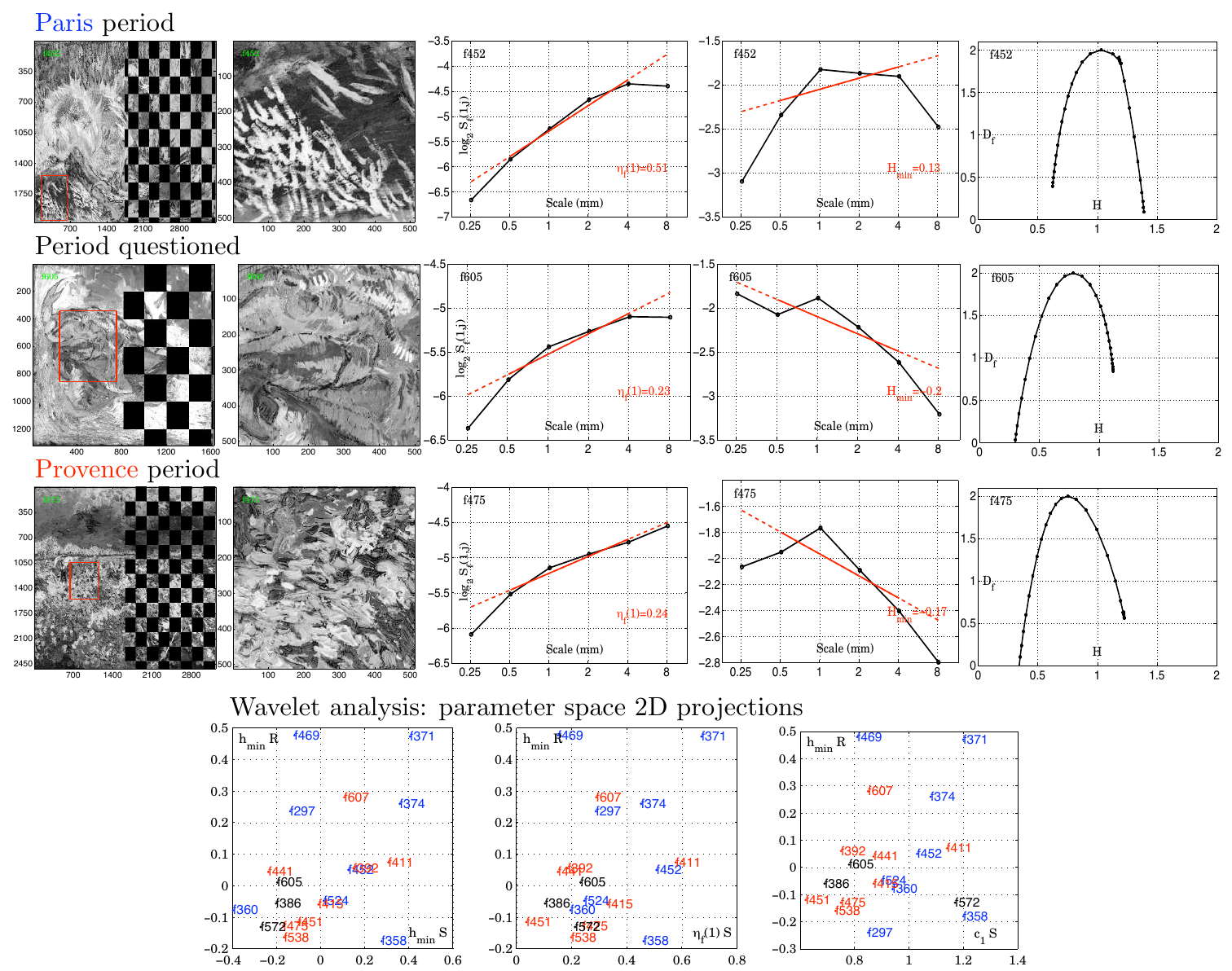}
\caption{\label{fig-VG} {\bf Van Gogh's Paintings.} Left columns: Representative examples of paintings of the Paris (top), Provence (bottom) and unknown (middle) periods, with selected homogeneous patches (second column),  Courtesy \emph{Image Processing for Art Investigation} research program, {cf. \tt www.digitalpaintinganalysis.org/}, and the Van Gogh and Kr\"oller-M\"uller Museums (The Netherlands).  Third column: log-log plots for the wavelet coefficient based $S_f(1,j)$ (as in (\protect \ref{equ-WSF})).  Fourth column: $\Hmin$. 
Right column: Wavelet leader based multifractal spectra, suggesting that the questioned painting has multifractal properties close to that of the Provence period example. Bottom row: 2D projections of multifractal attributes, indicating that paintings from the Paris period (blue) are overall globally more regular than paintings from the Provence period (red).}
\end{figure}

Applications described so far all implied the multifractal analysis of 1D signals only. 
However, B. Mandelbrot always advocated  that fractal analysis should  be applied to images, to characterize the roughness of textures (cf. e.g., \cite{H18,H17}). 
The wavelet coefficient based analysis of selfsimilarity analysis has been readily and immediately extended to 2D fields and images, thanks to the straightforward formulation of the 2D discrete wavelet transform (cf. e.g., \cite{Mallat1998}). 
The estimation of the corresponding $H$ parameter has then been used to characterize textures in images, often in the context of biomedical applications (cf. \cite{BR10,Bergou}). 
The extension of multifractal formalisms to images turned out to be significantly more intricate. 
Direct extension of the Kolmogorov scaling function were used to analyze rainfall repartitions by D. Schertzer and S. Lovejoy, see \cite{SL87}.
A second attempt based on the modulus maxima of the coefficients of the 2D continuous wavelet transform was proposed by A. Arneodo and his collaborators; it  employed in medical \cite{Arneodo2003a} or geophysical applications \cite{ard02}.
The wavelet leader formalism recently proposed and based on a 2D discrete wavelet transform enabled a significant step forward toward the use of multifractal analysis in higher dimensions. 
Its potential benefits have been assessed in \cite{WENDT:2009:E}, and it has recently been used to characterize texture in paintings \cite{Abry2012}. 

The paintings that Van Gogh performed while he was living in France are grouped into two periods by art experts and art historians: 
The Paris period, ending early 1888, and the later Provence period.
While art experts agree on the classification of most of Van Gogh's paintings, some are still not clearly attributed to either period. 
To perform such classification, art experts often make use of  material or stylometric features such as the density of brush strokes, their size or scale, the thickness of contour lines, the layers, the chemicals of the colors, \ldots.

Aiming at promoting the investigation of computer-based image processing procedures to assist art experts in their judgement, the \emph{Image Processing for Art Investigation} research program ({cf. \tt www.digitalpaintinganalysis.org/}), in collaboration with the Van Gogh and Kr\"oller-M\"uller Museums (The Netherlands), made available to research teams two sets of (partial and low resolution digitized versions of) $8$ Van Gogh's paintings each from the Paris and Provence period, and $3$ painting whose period remains undecided (nomenclature below corresponds to the Van Gogh museum catalog). 
The digitized copies are available through the usual RGB (Red, Blue, Green) channels and converted into the HSV (Hue, Saturation, Value) channels, hence enabling to process the gray-level light intensity information and the color information. 

Three examples, corresponding to the Paris period, the Provence period and unclassified paintings, respectively, are shown in Fig.~\ref{fig-VG} (left column) together with examples of homogenous patches selected for analysis (second column). 
The 2D wavelet leader based multifractal formalism has been applied, for each selected patch, independently to the six RGB and HSV channels. 
Examples of obtained multifractal spectra are shown in Fig.~\ref{fig-VG} (right column) and suggest, for example, that $f605$ is closer to $f475$ than to $f452$ and hence is likely to belong to the Provence period.
For each patch of the $2 \times 8+3$ paintings and each channel, multifractal attributes $\Hmin, \eta_f(1), c_1, \dots$ are estimated. 
Fig.~\ref{fig-VG} (bottom row) shows 2D projections of various pairs of such attributes, which all indicate that on average the Paris paintings (blue clusters) are globally more regular than the Provence paintings (red clusters). 
Also, the Saturation and Red channels are found to be the most discriminant.
This suggests that both a color and an intensity information (saturation and Red) are useful for discrimination. 
Interestingly, it has been reported by art experts that saturation is a key feature to distinguish between the Paris and Provence periods in Van Gogh's paintings. 
Finally, these 2D projections suggest that $f605$ and $f386$ belong to the lower left red cluster and hence to the Provence period, while conclusion is less straightforward for $f572$.
Such results are left to the appreciation of art experts. 
These analyses are detailed in \cite{Abry2012}.

\section{Conclusion}
\label{sec-conclusion}

A key issue for the future developments of multifractal analysis is  to disentangle (at least) four  different acceptations that can be associated with the key word \emph{multifractal}, and that are often misleadingly confused one with another, sometimes leading to potential misunderstandings and erroneous interpretations. \\

\noindent {\bf Theoretical multifractal analysis and local regularity.} \quad From the theoretical or mathematical side, multifractal analysis is deeply tied  with  the analysis of the fluctuations along time (or space) of the local regularity of the function or data of interest, these fluctuations being often coined as \emph{irregularity}. 
Usually, this local regularity is measured in terms of H\"older exponents $h_f(x)$.
The global description of the irregularity of $f$ is gathered in the multifractal spectrum $D_f(h)$, which consists of the Hausdorf dimensions of the geometrical iso-H\"older sets. 
For mathematics, multifractal is hence all about H\"older exponent and (fractal, or box or Haussdorf) dimensions of geometrical sets. 
Though listed first here, the mathematical formulation of multifractal analysis theory can be considered as the most recent piece of the history of multifractal, as it really started in  the  90s  only, see \cite{Bmp92}, as  a consequence of the seminal article of G. Parisi and U. Frisch  \cite{ParFri85}. 
Note however that precursors can be found, for instance in articles dealing with slow and fast points of Brownian motion: Recall that this specific process  displays exceptional points where the modulus of continuity is slightly  smaller (resp. bigger) than a.e.; the corresponding sets are random fractals. Following the seminal contribution of J.-P. Kahane, their dimensions have been  precisely determined,  respectively by S. Orey and  S. J. Taylor  for  fast points,  see  \cite{OrTay}, and by  E. Perkins  for slow points, see \cite{Per}. \\

\noindent {\bf Multifractal formalism and scale invariance.} \quad From the more practical perspective of real-world data analysis, scaling functions, consisting of time (or space) averages of the $p$-th power of scale dependent quantities (such as increments, oscillations, wavelet coefficients or wavelet leaders) --- also often referred to, in statistical physics, as partition or structure functions --- constitute the core aspect of the notion of multifractal. 
Importantly, practical multifractal analysis amounts in essence to assuming that scaling functions behave as power laws with  respect to the analysis scales, in the range  of  scales available in the data. 
The description of these power laws naturally leads to the notion of scaling exponents. 
Equivalently, these scaling exponents are termed scaling functions, with the practical and often implicit request that the limit underlying its definition actually exists (cf. (\ref{kolmo}) versus (\ref{scalKolmo1})). 

These power law behaviors --- or scale invariance, or scaling --- of scale dependent quantities, and the corresponding scaling exponents (or scaling functions), were at the heart of the intuition that B. Mandelbrot developed, connecting scaling properties to (fractal) dimensions, to selfsimilarity and/or to fluctuations of local regularity, and relating scaling to self similar random walks or multiplicative cascades. 

The multifractal formalism (in reference to the Thermodynamic formalism) therefore consists of the mathematical developments that relate scale invariance and scaling functions to local regularity and multifractal spectra.
In some sense, the multifractal spectrum constitutes the link between formal mathematics and physics or data analysis. \\

\noindent {\bf Multifractal processes.} \quad Multifractal processes do not constitute a well-defined class. 
In essence, this notion  refers to processes (or functions) that can be fruitfully used to model scale invariance properties in data. 
Following the seminal distinctions proposed by B. Mandelbrot, they can be thought of as falling into two major classes:
Selfsimilar random walks and multiplicative cascades (or martingales). 

The first one, now very classical after the seminal review of B. Mandelbrot  on fBm \cite{mvn68}, pertains to the class of additive models: 
A large step (or increment) can be obtained as the sum of many different small steps, yet sharing the same statistical properties. 
This class is hence deeply associated with selfsimilarity. 
Often, this class is confused with its emblematic but restrictive representative fBm, the only Gaussian selfsimilar process with stationary increments. 
The Gaussian nature of this process and its fully parametrized formulation lead to the formulation of various parametric estimation of its selfsimilarity exponent $H$. 
Such parametric estimation are sometimes incorrectly considered as a multifractal analysis of data. 
This is inaccurate since multifractal analysis aims at measuring much richer properties of data than the sole selfsimilarity and, above all, does not rely on the assumption that data are selfsimilar or Gaussian. 
Multifractal analysis thus has a much broader scope and can also be applied to  limits of other random walks,  such as L\'evy processes (and, more generally,  jump processes), see \cite{Panor} for overviews on the mathematical side.

The second class, multiplicative cascades, is usually considered as that of archetypal multifractal processes, in the sense that their (scaling functions and) multifractal spectra can be computed theoretically and found to often consist of \emph{continuous, smooth bell-shaped} functions. 
Therefore, for such processes, multifractal analysis unveils information that can not easily be obtained with other analysis tools (as opposed to the previous class of selfsimilar random walks). 

Furthermore, as mentioned earlier, and again following one of B. Mandelbrot's intuitions, this particular class can be extended by combinations (or subordinations) of selfsimilar process with multiplicative cascades (for instance fBm in multifractal time mentioned above, or L\'evy processes in multifractal time, as studied by J. Barral and S. Seuret  \cite{BaSe3,BaSe2}).

It is however clear that these examples do not even give  a vague idea of  the tremendous variety of all possible multifractal functions.  
A few  mathematical  results actually comfort this idea (and allow to transform it into a precise mathematical statement):  
Generically,  \emph{most} functions in a given function space are multifractal and satisfy the multifractal formalism, see \cite{JafFRa,JafBaire} for precise statements.  
Therefore, multifractality is not a rare property, and is certainly not restricted to the few examples which have been  investigated up to now. 

The classical opposition between monofractal versus multifractal processes, often used in practical data analysis, is not well grounded and somehow irrelevant. 
Confusion stems from the fact that Gaussian fBm is characterized by a single H\"older parameter $h_f (x) = H$ and is hence monoh\"older (its spectrum $D_f(h)$ is supported by  a single point). 
But selfsimilarity does not in general imply this monoh\"older property. 
Instead, the classification of data might be by opposing additive to multiplicative structure.
Indeed, the physics (or the biology or else) underlying data production may significantly differ depending on whether it is driven by additive mechanisms, or by  multiplicative ones. 
In that respect, multifractal analysis may help as the multifractal spectra of selfsimilar random walks generally differ from those of multiplicative constructions. 
However, in opposition to what is often used in applications, multifractal analysis per se, understood as the measurement of the scaling function and multifractal spectrum, is not sufficient for proving that there exist any additive or multiplicative structures underlying the data.
Multifractal analysis  does not aim at stating definite answers to fuzzy questions such as: {Are data following a cascade model or not?}
but, instead,  it provides quantitative answers to more restricted questions such as: Is $\Hmin >0$? Is $\zeta_f(2) > 1$?
In turns, these precise measurements can contribute to help practitioners to formulate hypotheses regarding data, or to classify them. \\

\noindent {\bf Practical multifractal analysis: signal and image processing.} \quad As stated in the previous paragraph, practical multifractal analysis amounts to measuring from data multifractal attributes such as $\eta_f(p), \Hmin, \zeta_f(p), c_m, L_f(h)$. 
These quantities can be computed from data without assuming any a priori data model. 
In other words, data to which multifractal analysis is  applied need not stem from any exact selfsimilar random walk or multiplicative construction. 
Multifractal analysis is hence, by nature,   a non parametric analysis tool that can be applied to any type of data, be they (multi)fractal or not. 
Whether 
$L_f(h)$ can be related with $D_f(h)$ or not does not prevent practitioners to use  these \emph{multifractal} attributes as tools for classification, diagnosis, detection, or else. 
This is notably the case in image processing, where it is very unlikely that databases containing images of widely different natures can be associated with specific construction models (such as random walks or cascades). 

Furthermore, in practice, one question is always central and was already raised by B. Mandelbrot: What multiresolution quantity should one start with? Increments, oscillations, wavelet coefficients, wavelet leaders?
In this respect, the present contribution aims at providing clear and general answers.
First, increments and oscillations only address the restricted case where $0\leq h_f(t) < 1, \forall t$; while wavelet based quantities enable to bypass that restriction by selecting smooth enough wavelets with a large enough number of vanishing moments. 
Therefore, wavelet coefficients generalize increments and the wavelet leaders constitute the generalized counterpart of oscillations. 
Second, wavelet coefficients and wavelet leaders should not be opposed. 
Wavelet coefficients enable to measure information in data such as the uniform H\"older exponent and a global selfsimilarity type scaling exponent;  they should hence always be applied first to data. 
Wavelet leaders should be applied next, when relevant ($\Hmin > 0$) and when seeking for additional and more refined analysis of scaling properties in data, via a full collection of scaling exponents $\zeta_f(p)$, or equivalently $c_m$, or equivalently $L_f(h)$.\\

\noindent {\bf  Multifractal toolbox.} \quad The practical solutions proposed here in terms of discrete wavelet transform coefficients and leaders are believed to be one of the methods enabling the best practical achievements in terms of combining firm mathematical grounding, satisfactory estimation performance, low implementation and computational costs, robustness and versatility. 
It can be accompanied with non parametric bootstrap procedures, providing confidence intervals for estimates and hypothesis test procedures. 
Technical details are reported in \cite{WENDT:2007:E} and {\sc Matlab} toolboxes implementing these tools are made publicly available by the authors.


\providecommand{\noopsort}[1]{}

\end{document}